\documentclass{article}

\usepackage{color}
\usepackage[colorlinks,citecolor=blue,linkcolor=red,urlcolor=cyan]{hyperref}
\usepackage{amssymb}
\usepackage{amsmath}
\usepackage{amsthm}

\usepackage{graphicx}
\usepackage{subfigure}
\usepackage{float}

\usepackage{footnote}
\makesavenoteenv{eqnarray}

\usepackage{enumerate}

\usepackage[a4paper,top=2cm,bottom=2cm,left=2cm,right=2cm,marginparwidth=1.8cm]{geometry}

\newtheorem{thm}{Theorem}[section]
\newtheorem{de}{Definition}[section]
\newtheorem{cor}{Corollary}[section]
\newtheorem{con}{Conjecture}[section]
\newtheorem{lem}{Lemma}[section]
\newtheorem{pr}{Proposition}[section]
\newtheorem{rem}{Remark}[section]
\newtheorem{cl}{Claim}[section]
\newtheorem{qu}{Question}[section]

\numberwithin{equation}{section}

\newcommand{\cp}[1]{\vcenter{\hbox{#1}}}

\newcommand{\re}{\mathrm{Re}}
\newcommand{\im}{\mathrm{Im}}
\newcommand{\D}{\mathbf{d}}
\newcommand{\li}{\mathrm{Li}_2}
\newcommand{\dd}{\mathrm{D}_2}
\newcommand{\vol}{\mathrm{Vol}}
\newcommand{\cs}{\mathrm{CS}}
\newcommand{\cv}{\mathrm{CV}}
\newcommand{\rt}{\mathrm{RT}}
\newcommand{\tv}{\mathrm{TV}}
\newcommand{\hess}{\mathrm{Hess}}

\title{On the Volume Conjecture for hyperbolic Dehn-filled
$3$-manifolds along the twist knots}

\author{Huabin Ge, Yunpeng Meng, Chuwen Wang and Yuxuan Yang}

\begin{document}

\maketitle

\begin{abstract}
    For a twist knot $\mathcal{K}_{p'}$, let $M$ be the closed $3$-manifold obtained by doing $(p, q)$ Dehn-filling along $\mathcal{K}_{p'}$. In this article, we prove that Chen-Yang's volume conjecture holds for sufficiently large $|p| + |q|$ and $|p'|$ for $M$. In the proof, we construct a new ideal triangulation of the Whitehead link complement which is different from Thurston's triangulation. Our triangulation has led to some new discoveries regarding symmetry, including insights into ``sister manifolds'' as introduced by Hodgson, Meyerhoff, and Weeks.
\end{abstract}

\tableofcontents

\section{Introduction}
The history of the ``volume conjecture" is not very long, but it is very important, and can even be regarded as one of the most important problems in the field of three dimensional geometric topology after Thurston's geometrization program and the virtual Haken conjecture were solved. Since the introduction of Jones polynomials, quantum topology has been widely studied by mathematicians and physicists, and many quantum invariants of knots and three-dimensional manifolds have been discovered. The volume conjecture began with
Kashaev \cite{Kashaev1,Kashaev2}, who introduced a complex valued invariant $\langle L\rangle_N$ of a link $L$, for each integer $N$. Based on the profound observation of some examples, he further conjectured the asymptotic exponential growth rate of $|\langle L\rangle_N|$ was given by the hyperbolic volume of $S^3\backslash L$. Later then, Murakami and Murakami \cite{Murakami} proved that Kashaev's invariant coincide with the values of the colored Jones polynomials at a certain root of unity. They reformulated Kashaev's conjecture as follows.

\begin{con}[Kashaev-Murakami-Murakami]
    For a hyperbolic link $L$ in $S^3$, denote $J_N(L,q)$ by its $N$-th colored Jones polynomial taking value at $q$. Then
    $$2 \pi \lim _{N \to \infty} \frac{\log \left|J_{N}(L, e^{\frac{2\pi\sqrt{-1}}{N}})\right|}{N}=\mathrm{Vol}\left(S^{3} \backslash L\right),$$
    where $\mathrm{Vol}(S^3\backslash L)$ denotes the hyperbolic volume of the complement of $L$ in $S^{3}$.
\end{con}

Surprisingly, the volume conjecture connects several completely different fields, such as TQFT and hyperbolic structures of three manifolds, which is why it is so magical and fascinating. As is well known, the geometric information of a three-dimensional manifold is given by Thurston's geometrization program, while its quantum topology is given by its quantum invariants, which are derived from the representation of quantum groups and appear to be irrelevant to its geometry at first glance. Therefore, the geometry and quantum topology of the three-dimensional manifold provide its geometric and topological information from different directions.

As soon as the volume conjecture was proposed, it attracted much attention, and many examples of knots or links were proven. For instance, the figure-eight knot by Ekholm, $5_2$ by Kashaev-Yokota, torus knots \cite{kashaev2000proof}, $(2, 2m)$-torus links \cite{hikami2004quantum}, Borromean rings \cite{garoufalidis2005volume}, twisted Whitehead links \cite{Zheng2007}, Whitehead doubles of torus knots \cite{Zheng2007}, Whitehead chains \cite{veen2008proof} and all hyperbolic knots with at most seven crossings \cite{ohtsuki2017asymptotic}\cite{oth7}. In addition, various extended versions of the volume conjecture have been proposed, and certain cases have been demonstrated. One version is to consider the growth rate of $J_N(L, q)$, rather than its absolute value, which should be given by the complex volume (see the following (\ref{complex-volume})); Another version is quite natural: for three-dimensional hyperbolic manifolds (closed, cusped, or those with totally geodesic boundary), they have very famous quantum invariants, such as the Reshetikhin-Turaev invariants \cite{RT1,RT2} and the Turaev-Viro invariants \cite{TV}. Is the asymptotic logarithmic growth rate of these invariants related to the hyperbolic volume? Witten's Asymptotic Expansion Conjecture (see \cite{ohtsuki2002}) predicts that, when evaluated at $q=e^{\frac{\pi\sqrt{-1}}{N}}$, the Witten-Reshetikhin-Turaev invariants of a 3-manifold only grow polynomially, with a growth rate related to classical invariants of the manifold such as the Chern-Simons invariant and the Reidemeister torsion. However, Chen-Yang \cite{chen-yang} investigated the asymptotic behavior of these invariants at $q=e^{\frac{2\pi\sqrt{-1}}{N}}$, and conjectured that they are of exponential order based on numerical computations. To be precise, consider a closed $3$-manifolds $M$, and let $\cs(M)$ be its Chern-Simons invariant. Denote
\begin{equation}
\label{complex-volume}
\cv(M)\equiv\vol(M) + \sqrt{-1}\cs(M) \mod\pi^2\sqrt{-1}
\end{equation}
by the complex volume of $M$. Let $\rt_r(M, q)$ be the $r$-th Reshetikhin-Turaev invariant, and $\tv_r(M, q)$ be the $r$-th Turaev-Viro invariant, both evaluated at the root $q=e^{\frac{2\pi\sqrt{-1}}{r}}$. Chen-Yang's volume conjecture for those invariants says:
\begin{con}[Chen-Yang]
\label{Conjecture-CY}
For a hyperbolic 3-manifold $M$. Then for $q=e^{\frac{2\pi\sqrt{-1}}{r}}$, and for $r$ running over all positive odd integers,
\begin{equation}
\label{equation-TV}
\lim_{r \to +\infty}\frac{2\pi}{r}\log(\tv_r(M, q))= \vol(M).
\end{equation}
If further assuming $M$ is closed and oriented, then for a suitable choice of the arguments,
\begin{equation}
\label{equation-RT}
\lim_{r \to +\infty}\frac{4\pi}{r}\log(\rt_r(M, q))\equiv \cv(M)\mod \pi^2\sqrt{-1}.
\end{equation}
\end{con}

Soon, Conjecture \ref{Conjecture-CY} were verified for the figure-8 knot, the Borromean rings complement \cite{DKY} by Detcherry-Kalfagianni-Yang and the fundamental shadow link complements \cite{BDKY-JDG} by Belletti-Detcherry-Kalfagianni-Yang. For the 3-manifold $M_p$ obtained from the figure-8 knot by a $|p|>4$ Dehn filling, Ohtsuki proved Conjecture \ref{Conjecture-CY} in \cite{Oh41}, where he also obtained a full asymptotic expansion of $\rt_r(M_p, q)$. Wong-Yang \cite{Wong-Yang41} further showed that Conjecture \ref{Conjecture-CY} hold for all hyperbolic 3-manifolds obtained by doing rational surgeries along the figure 8-knot. In the proof of \cite{BDKY-JDG}, an estimate on the growth of the quantum $6j$-symbol is essential. We refer Costantino \cite{Costantino}, Chen and Murakami \cite{Chen-Mura} for the asymptotic growth rate or the asymptotic expansion of quantum $6j$-symbols. 

Let $\mathcal{K}_{p'}$ be the twist knot with twist number $p'$ (Figure \ref{twist knot}), Chen and Zhu \cite{Chen-Zhu-1,Chen-Zhu-2} studied the asymptotic expansions of various quantum invariants related to the twist knot $\mathcal{K}_{p'}$ with $p'\geq6$. As a consequence, they asserted that the Kashaev-Murakami-Murakami volume conjecture is valid for $\mathcal{K}_{p'}$ with $p'\geq6$. This inspired us to consider Chen-Yang's volume conjecture for rational Dehn fillings of $\mathcal{K}_{p'}$. In this article, we prove Conjecture \ref{Conjecture-CY} for closed 3-manifolds obtained by doing $(p, q)$ \textbf{rational} surgery along $\mathcal{K}_{p'}$ with $p,q$ are coprime and $|p'|$, $|p|+|q|$ large enough (we note that, during the preparation of this paper, Chen-Zhu posed their results \cite{Chen-Zhu-3} on Chen-Yang's volume conjecture for closed hyperbolic 3-manifolds obtained by doing \textbf{integral} $p$-integral surgery along $\mathcal{K}_{p'}$ with $p'$, $p$ sufficiently large). Our first result says:
    \begin{figure}[H]
        \centering
        \includegraphics[width=0.25\linewidth]{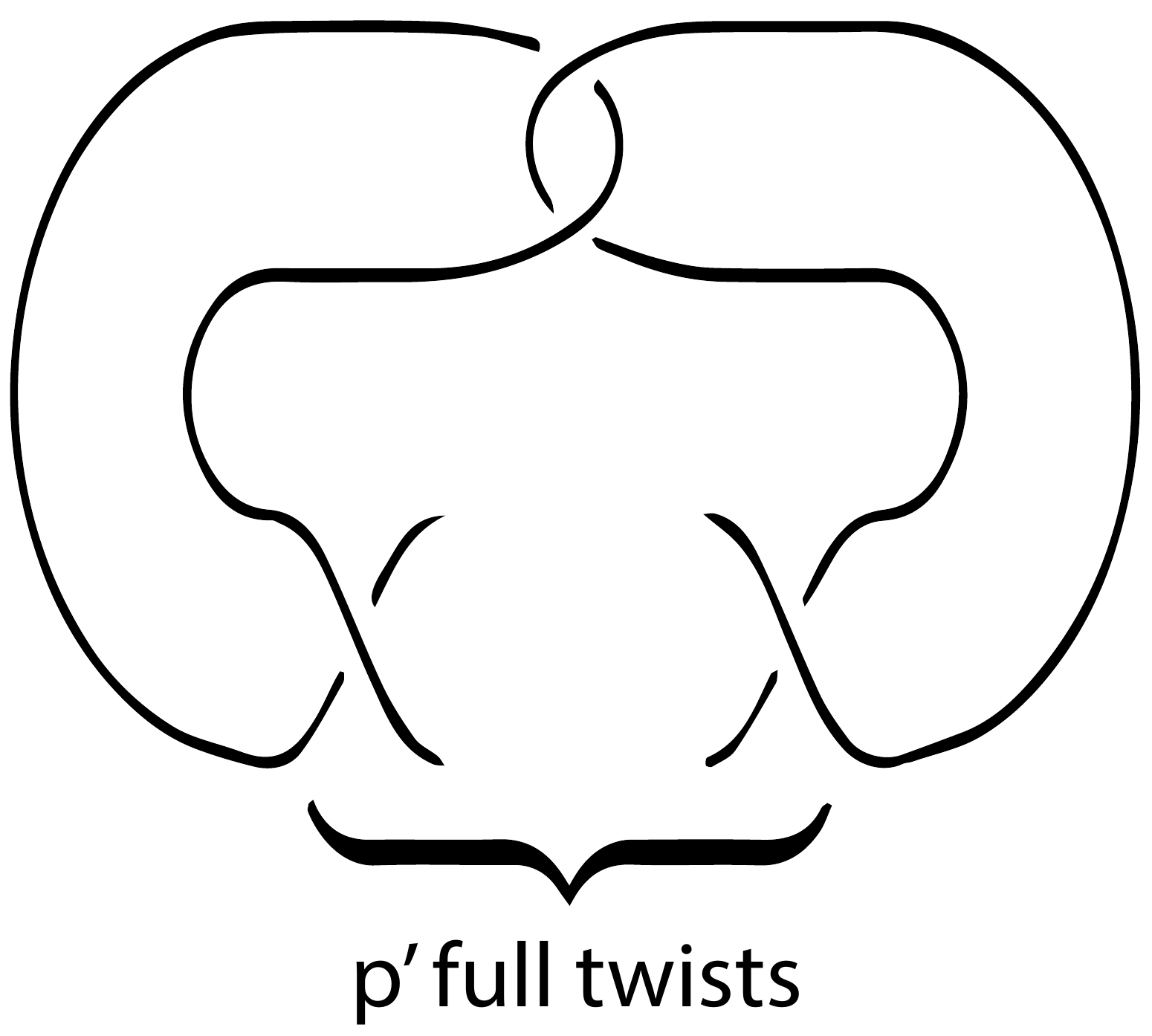}
        \caption{A twist knot $\mathcal{K}_{p'}$}
        \label{twist knot}
    \end{figure}

\begin{thm}
\label{thm1}
Let $M=\mathcal{K}_{p'}(p,q)$ be the closed 3-manifold obtained by doing the $(p,q)$ Dehn surgery along the twist knot $\mathcal{K}_{p'}$, then for $|p|+|q|$, $|p'|$ sufficiently large, $r$ positive odd integer,
\begin{equation}
\rt_r(M, e^{\frac{2\pi\sqrt{-1}}{r}})=C(r)t(M)e^{\frac{r}{4\pi}\cv(M)}\Big(1 + O\Big(\frac{1}{r}\Big)\Big),
\end{equation}
where $ C(r)$ is of norm $1$ independent of the geometry structure of $M$, $t(M)$ is a invariant of $\mathcal{K}_{p'}(p,q)$.
\end{thm}

At the root $q=e^{\frac{2\pi\sqrt{-1}}{r}}$, by \cite{DKY}, \cite{Turaev1}, \cite{Witten1}, \cite{Roberts1}, \cite{Wong-Yang41}, the Turaev-Viro invariant is up to a scalar the square of the norm of the Reshetikhin-Turaev invariant for a closed oriented $3$-manifold, i.e.
$$\tv_r(M)=2^{b_2(M)-b_0(M)+2}|\rt_r(M)|^2.$$
where $b_0(M)$, $b_2(M)$ is the zeroth and second $\mathbb{Z}_2$-Betti numbers of $M$. As a consequence, we have

\begin{thm}
\label{thmtv}
Let $M=\mathcal{K}_{p'}(p,q)$ be the closed 3-manifold obtained by doing the $(p,q)$ Dehn surgery along the twist knot $\mathcal{K}_{p'}$, then for $|p|+|q|$, $|p'|$ sufficiently large, $r$ positive odd integer,
\begin{equation}
\tv_r(M, e^{\frac{2\pi\sqrt{-1}}{r}}) =2^{1+b_2(M)}|t(M)|^2e^{\frac{r}{2\pi}\vol(M)} \Big(1 + O\Big(\frac{1}{r}\Big)\Big).
\end{equation}
\end{thm}


Using the asymptotic formulas in Theorem \ref{thm1} and Theorem \ref{thmtv}, we derive the following corollary.
\begin{cor}
Let $M=\mathcal{K}_{p'}(p,q)$ be the closed 3-manifold obtained by doing the $(p,q)$ Dehn surgery along the twist knot $\mathcal{K}_{p'}$, then for $|p|+|q|$, $|p'|$ sufficiently large, $r$ positive odd integer, Conjecture \ref{Conjecture-CY} is true, i.e.
    \begin{eqnarray*}
      \lim_{r \to \infty} \frac{4\pi}{r} \log \rt_r(M, e^{\frac{2\pi\sqrt{-1}}{r}})&\equiv&\cv(M)\mod \pi^2\sqrt{-1},\\
      \lim_{r \to \infty} \frac{2\pi}{r} \log \tv_r(M, e^{\frac{2\pi\sqrt{-1}}{r}})&=&\vol(M).
    \end{eqnarray*}
\end{cor}

Overall, we use Ohtsuki \cite{Oh41} and Wong-Yang's \cite{Wong-Yang41} framework to approach our main results. However, in the process of proof, there are many essential differences in our article. The main differences arise from the fact that the canonical ideal triangulation of the Whitehead link complement does not align with our potential function. Therefore, we manually constructed an ideal triangulation of the Whitehead link complement space, consisting of five ideal tetrahedra. Our triangulation can also be derived from Thurston's canonical triangulation through a series of Pachner moves.

We then encountered a difficulty that appeared for the first time in the proof of the volume conjecture: the Dehn filling coefficients corresponding to the potential function differed from what we wanted. Inspired by Hodgson-Meyerhoff-Weeks' conclusion \cite{HMW}, we proposed a theory called the \textbf{sister potential function}, which fixes this discrepancy in the Dehn coefficients. The theory provides an analytical explanation for Hodgson-Meyerhoff-Weeks' conclusion. In addition to this, it also provides inspiration for the proof of the volume conjecture for other cases.


Before beginning the calculations, we use Habiro's renowned cyclotomic expansions  \cite{habiro}(or both two forms of the Masbaum's formula \cite{masbaum})
to greatly simplify the subsequent analysis. Additionally, before estimating, we observed a symmetry in the Fourier coefficients, which generalizes the equality of the leading Fourier terms and the \textbf{Big Cancellation} discovered by Chen and Zhu  \cite{Chen-Zhu-1}. This discovery made the proof  more concise. Additionally, before using the saddle point method, we adjusted the integration region to satisfy Ohtsuki's convexity conditions.

Moreover, we have discovered an interesting lemma (Lemma \ref{lem2}) where we can directly get the relationship between the imaginary part of the critical value of the potential function and the volume without using the Neumann-Zagier-Yoshida Theory \cite{NZ}\cite{Yoshida}. This simple correspondence has an inspiration that the shapes of the idea tetrahedra can be directly observed from the form of the potential function. This greatly helps us find good idea triangulations in the future.

~

\textbf{Outline of the proof.} The proof follows the guideline of Wong-Yang and Ohtsuki's strategy. In Section \ref{sec2}, we provide the necessary concepts. In Section \ref{sec3}, we compute $\mathrm{RT}_r(\mathcal{K}_{p'}(p,q))$ and get the potential function $V$.

In Section \ref{sec4}, we utilize the discrete Fourier transform method to express $\mathrm{RT}_r(\mathcal{K}_{p'}(p,q))$ as a sum of integrals. We discover the symmetric of the Fourier coefficients, which is a generalization of the \textbf{Big Cancellation} as showed in \cite{Chen-Zhu-1}. As a corollary, the key terms which  exponential growth rate greater than the volume cancel out directly.

\begin{figure}[H]
        \centering
        \includegraphics[width=0.3\linewidth]{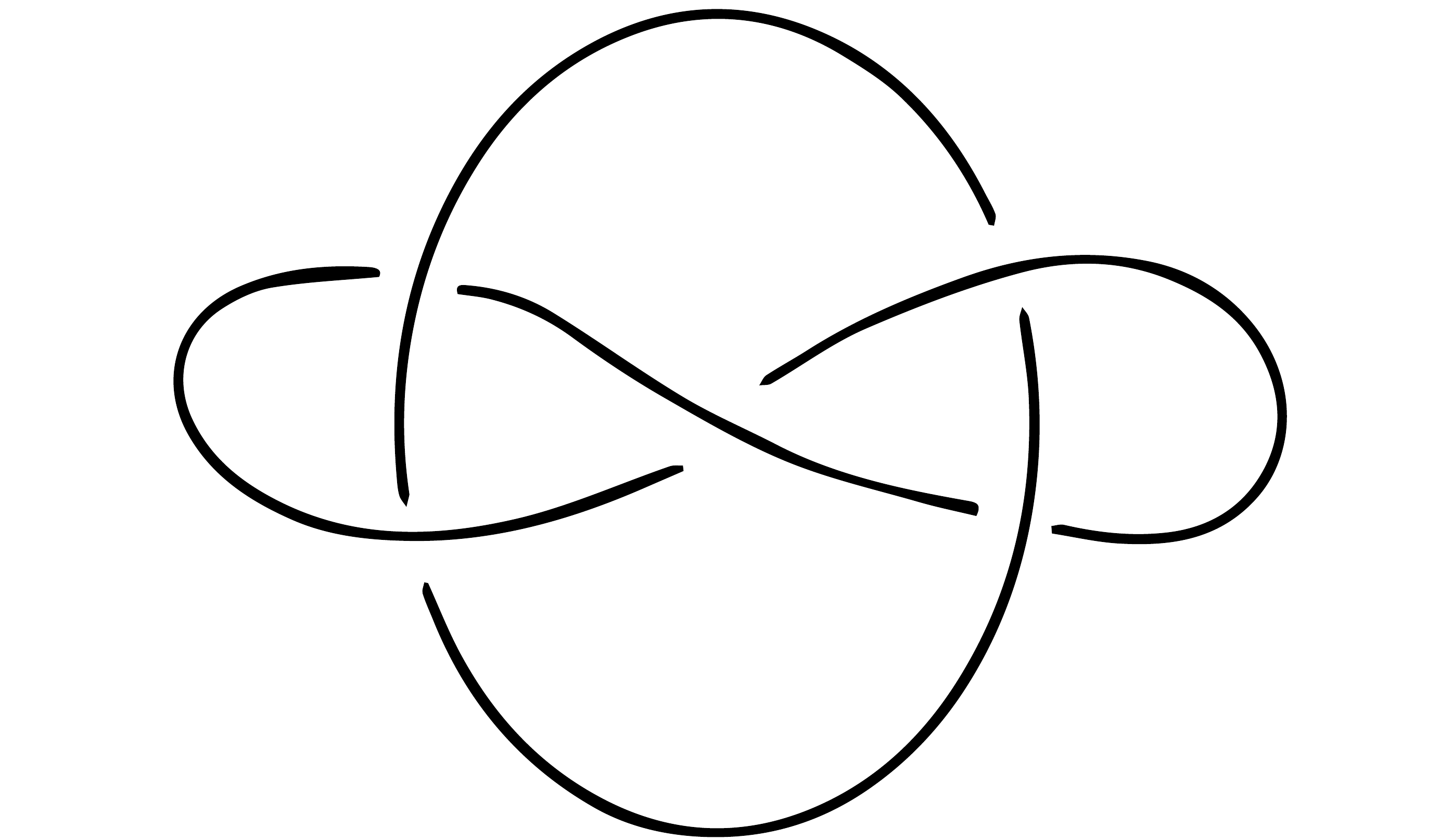}
        \caption{Whitehead link}
        \label{Whiteheadlink}
    \end{figure}

In Section \ref{sec5}, we find the relationship between $V$ and the geometry of $\mathcal{K}_{p'}(p,q)$. We find $V$ is corresponding to the geometry to $\mathcal{W}((p+4q,-q),(1,p'-\frac{1}{2}))$, where $\mathcal{W}$ is the Whitehead link (Figure \ref{Whiteheadlink}). However, $\mathcal{W}((p+4q,-q),(1,p'-\frac{1}{2}))$ is a $4\pi$ cone angle hyperbolic manifold, which is not a hyperbolic manifold in the common sense. To overcome this difficulty, we define the sister potential function $W$ and establish a connection between $W$ and the geometry of $\mathcal{W}((p,q),(1,-p'))$.  More precisely, Hodgson, Meyerhoff, and Weeks \cite{HMW} gives a geometric explanation (Figure \ref{sister manifolds}) of the intriguing phenomenon between the potential function $V$ and sister potential function $W$, and the critical points of the $W$ are equivalent to the hyperbolic structure equations of $\mathcal{W}((p,q),(1,-p'))$.

 In Section \ref{sec6}, we employed the saddle point method to demonstrate that the exponential growth rates of the remaining four terms are equal to the volume, and they will not cancel each other out, unlike the other terms, which have exponential growth rates that are subdominant with respect to the volume. Additionally, we provided the proof of Theorem \ref{thm1} and Theorem \ref{thmtv}.

 In Section \ref{sec7}, we discuss some interesting problems related to our work. Finally, we present a program for proving the volume conjecture.

 ~

\textbf{Acknowledgments.} The authors would like to thank Professor Tian Yang for helpful discussions, and Qing Lan for computer skills. Ge and Wang would like to thank Professor Qingtao Chen for his concern and helpful conversations. Ge wants to thank Professor Ruifeng Qiu for his encouragement. Meng is thankful for the encouragement and support from his mentor, Xuezhi Zhao. Yang would like to express his deepest gratitude to his advisor, Professor Yi Liu, for his invaluable guidance and unwavering support throughout his research journey. The paper is supported by NSFC, no.12341102, no.12122119.

\section{Preliminaries}\label{sec2}
\subsection{Reshetikhin-Turaev invariants}
    Reshetikhin-Turaev invariants is a quantum invariant of closed $3$-manifold. We will defined it by skein theory in this section.

    In a $3$-sphere $S^3$, a framed link $L$ refers to a collection of smoothly embedded disjoint thicken circles into $S^3$. A framed link diagram $D(L)$ in $S^3$ is obtained by taking a regular neighborhood of a link diagram in the plane. The following is an example of obtaining a framed link diagram from a link diagram.
\begin{figure}[H]
    \centering
    \includegraphics[width=0.5\linewidth]{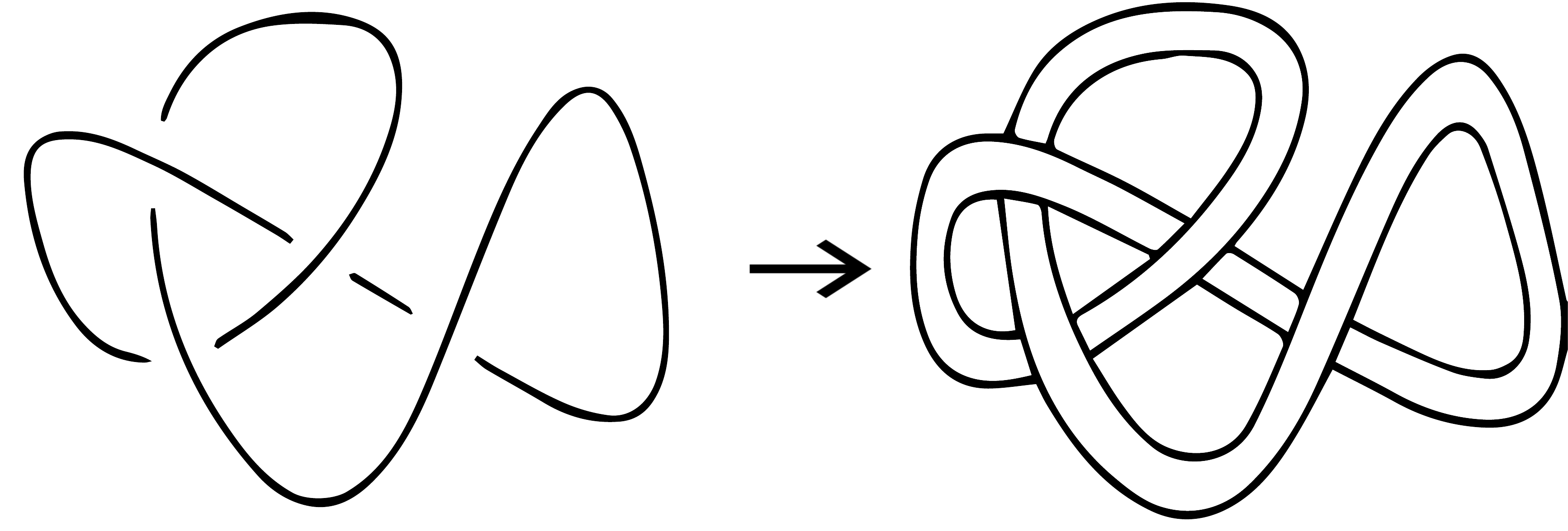}
    \caption{framed link diagram}\label{fig:framed link diagram}
\end{figure}
    It is not difficult to prove that any framed link can be represented by an appropriate framed link diagram through suitable ambient isotopy.

    The Kauffman bracket skein module  $\mathrm{K}_{r}(S^3) $  of  $S^3$   is the $ \mathbb{C} $-module generated by  framed link diagrams in $S^3$ modulo the following two relations
    \begin{enumerate}[(1)]
        \item  \emph{Kauffman Bracket Skein Relation:} \ $\cp{\includegraphics[width=1cm]{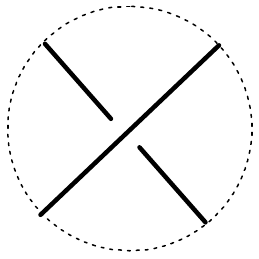}}\ =\ e^{\frac{\pi\sqrt{-1}}{r}}\ \cp{\includegraphics[width=1cm]{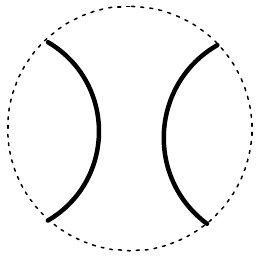}}\  +\ e^{-\frac{\pi\sqrt{-1}}{r}}\ \cp{\includegraphics[width=1cm]{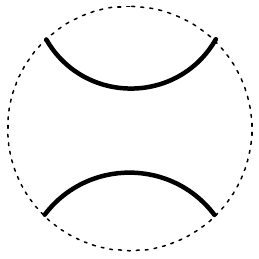}}.$

        \item \emph{Framing Relation:} \ $L \cup \cp{\includegraphics[width=0.8cm]{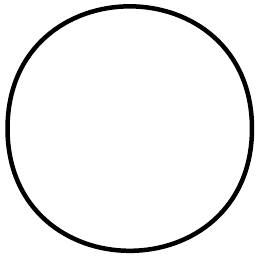}}=(-e^{\frac{2\pi\sqrt{-1}}{r}}-e^{-\frac{2\pi\sqrt{-1}}{r}})\ L.$
    \end{enumerate}

    Similarly, we can define the Kauffman bracket skein module $\mathrm{K}_{r}(A\times[0,1])$ on the product of an annulus $A$ with a closed interval, it is  the $ \mathbb{C}$-module generated by all framed link diagrams on the $A$ , and modulo by the two types of relations mentioned above. There is a commutative multiplication induced by the juxtaposition of  $A$, making it a  $\mathbb{C}$-algebra. It is not hard to prove that  $\mathrm{K}_{r}(A\times[0,1])$ is the polynomial algebra $\mathbb{C}[z]$, where $z$ is the core cricle of $A$. For any link diagram $D$ in $\mathbb R^2$ with $k$ ordered components and $b_1, \dots, b_k\in \mathrm K_r(A\times [0,1]),$ let
    $$\langle b_1,\dots, b_k\rangle_D$$
    be the complex number obtained by cabling $b_1,\dots, b_k$ along the components of $D$ considered as an element of $K_r(\mathrm S^3)$ then taking the Kauffman bracket $\langle\ \rangle$ as defined in \cite{lick}.

    Suppose $M$ is obtained by doing surgery on a three-dimensional sphere along a frame link diagram $D(L)$, $\sigma(L)$ denotes the signature of the linking matrix of $L$ . We can define the Reshetikin-Turaev invariant of M as follows. Let $e_n(z)$ be the n-th Chebyshev polynomial, which is defined by the recursive relation $e_n(z)=ze_{n-1}(z)-e_{n-2}(z)$ , where $e_0(z)=1$ and $e_1(z)=z$. The Kirby coloring  $\omega_{r} \in \mathrm{K}_{r}(A \times[0,1]) $ is then defined by
    $$\omega_{r}=\sum_{n=0}^{r-2}(-1)^{n}[n+1] e_{n},$$
    where $[n]$ is the quantum interger defined by
    $$[n]=\frac{e^{\frac{2 n \pi \sqrt{-1}}{r}}-e^{-\frac{2 n \pi \sqrt{-1}}{r}}}{e^{\frac{2 \pi \sqrt{-1}}{r}}-e^{-\frac{2 \pi \sqrt{-1}}{r}}}.$$
    Set
    $$\mu_{r}=\frac{\sin \frac{2 \pi}{r}}{\sqrt{r}}.$$
\begin{de}
    We can define the $r$-th Reshetikin-Turaev invariant as
    $$\mathrm{RT}_{r}(M)=\mu_{r}\left\langle\mu_{r} \omega_{r}, \ldots, \mu_{r} \omega_{r}\right\rangle_{D(L)}\left\langle\mu_{r} \omega_{r}\right\rangle_{U_{+}}^{-\sigma(L)},$$
    where  $U_+ $ is the diagram of  unknot with framing $1$.
\end{de}

\subsection{Dilogarithm, Lobachevsky functions and Bloch-Wigner functions}
Dilogarithm functions are closed related to the hyperbolic volume of $3$-manifolds. We will discuss some analytical property of dilogarithm functions in this section.

  Let  $  \log : \mathbb{C} \backslash(-\infty, 0] \to \mathbb{C}$  be the standard logarithm function defined by
  $$\log z=\log |z|+\sqrt{-1} \arg z$$
  with $ -\pi<\arg z<\pi $.

  The dilogarithm function $\li: \mathbb{C} \backslash(1, \infty) \to \mathbb{C}$  is defined by$$\li(z)=-\int_{0}^{z} \frac{\log (1-t)}{t} d t,$$ where the integral is along any path in  $\mathbb{C} \backslash(1, \infty)$  connecting 0 and  $z$, which is holomorphic in  $\mathbb{C} \backslash[1, \infty)$  and continuous in $\mathbb{C} \backslash(1, \infty)$. In the unit disk $ \{z \in \mathbb{C}||z|\le1\},   \li(z)=\sum_{n=1}^{\infty}\limits \frac{z^{n}}{n^{2}}$.

   The dilogarithm function satisfies the follow properties (see eg. Zagier \cite{dilog}),
   \begin{equation}\label{e1}
     \li\left(\frac{1}{z}\right)=-\li(z)-\frac{\pi^{2}}{6}-\frac{1}{2}(\log (-z))^{2}.
   \end{equation}
Moreover, on the unit circle
   $\{z=e^{2 \sqrt{-1} \theta} \mid 0 \leqslant \theta \leqslant \pi\}$,
    \begin{equation}\label{e2}
      \li\left(e^{2 \sqrt{-1} \theta}\right)=\frac{\pi^{2}}{6}+\theta(\theta-\pi)+2 \sqrt{-1} \Lambda(\theta),
    \end{equation}
    where  $\Lambda: \mathbb{R} \to \mathbb{R}$ is the Lobachevsky function (see eg. Thurston's notes \cite{Thurston} Chapter 7) defined by
    $$\Lambda(\theta)=-\int_{0}^{\theta} \log |2 \sin t| d t.$$
    The Lobachevsky function is an odd function of period  $\pi$. It achieves the absolute maximums at $ k \pi+\frac{\pi}{6}, k \in \mathbb{Z}$, and the absolute minimums at  $k \pi+\frac{5 \pi}{6}, k \in \mathbb{Z}$.

We recall the famous $5$-term formula for dilogarithm. This corresponds to this Pachner $2$-$3$ move in geometry.
    \begin{pr}[\cite{dilog}]\label{5terms1}
    \label{5terms}
         $$\li(x)+\li(y)-\li\Big(\frac{x}{1-y}\Big)-\li\Big(\frac{y}{1-z}\Big)+\li\Big(\frac{xy}{(1-x)(1-y)}\Big)=-\log(1-x)\log(1-y).$$
    \end{pr}

\begin{lem}\label{lemLambda}
    For $x, \theta\in\mathbb{R}$ where $|a|<\frac{1}{2}$, there is a constant $C$ such that
    $$|\Lambda(x+a)-\Lambda(x)|<C|a\log|a||.$$
\end{lem}

\begin{proof}[Proof of Lemma \ref{lemLambda}]
    Because of $\Lambda(x)$ is $C^1$ for $x\ne k\pi, k\in\mathbb{Z}$, it suffices to prove the case $x=0$, i.e.
    $$|\Lambda(a)|=O(|a|\log a)$$
    for small $a>0$. We have

    \begin{eqnarray*}
      |\Lambda(a)|&=&\frac{1}{2}|\im\li(e^{2\sqrt{-1}a})|\\
      &=&\frac{1}{2}|\sum_{n=1}^{\infty}\frac{\sin 2na}{n^2}|\\
      &\le&\frac{1}{2}(|\sum_{1\le n\le\frac{\pi}{2a}}\frac{\sin 2na}{n^2}|+|\sum_{n>\frac{\pi}{2a}}\frac{\sin 2na}{n^2}|)\\
      &<&\frac{1}{2}(\sum_{1\le n\le\frac{\pi}{2a}}\frac{2n|a|}{n^2}+\sum_{n>\frac{\pi}{2a}}\frac{1}{n^2})\\
      &=&|a|\sum_{1\le n\le\frac{\pi}{2a}}\frac{1}{n}+\frac{1}{2}\sum_{n>\frac{\pi}{2a}}\frac{1}{n^2}\\
      &=&O(|a\log a|)+O(a)\\
      &=&O(|a\log a|).
    \end{eqnarray*}
\end{proof}

    \begin{de}
    The Bloch-Wigner function $\dd(z)$ is related to the dilogarithm function by $\dd(0)=\dd(1)=0$ and
        $$\dd(z)=\im(\li(z))+\log|z|\cdot\im\log(1-z), \;\text{if }z\in\mathbb{C}\backslash\{0,1\}.$$
    \end{de}

    \begin{pr}[\cite{dilog}]\label{prD}
        $\dd(z)$ enjoys the following remarkable properties :
        \begin{enumerate}
          \item $\dd(z)$ is real analytic on $\mathbb{C}\backslash\{0,1\}$;
          \item $D(\overline{z})=-D(z)$;
          \item $D(z)$ is positive on the upper plane $\{z\in\mathbb{C}|\im z>0\}$;
          \item $\dd(e^{2\sqrt{-1}\theta})=\li(e^{2\sqrt{-1}\theta})=2\Lambda(\theta)$ for $\theta\in\mathbb{R}$;
          \item $\dd(z)=\dd(1-\frac{1}{z})=\dd(\frac{1}{1-z})=-\dd(\frac{1}{z})=-\dd(1-z)=-\dd(\frac{z}{z-1})$ for $z\ne0,1$;
          \item $5$-term formula for $\dd$:
          $$\dd(x)+\dd(y)-\dd\Big(\frac{x}{1-y}\Big)-\dd\Big(\frac{y}{1-z}\Big)+\dd\Big(\frac{xy}{(1-x)(1-y)}\Big)=0.$$
        \end{enumerate}
    \end{pr}

    The following lemma means that the imaginary part of the critical value of the potential function equal to the volume, which will be showed in Section \ref{sec5.4}.

    \begin{lem}\label{lem2}
        Suppose $\vec{\mathbf{z}}=\left(\begin{matrix}
                                     z_1 \\
                                     \vdots \\
                                     z_n
                                   \end{matrix}\right)\in\mathbb{C}^n,\vec{\mathbf{r}}=(r_1,\cdots,r_m)\in\mathbb{R}^m,L=\left(\begin{matrix}
                                     \vec{l_1} \\
                                     \vdots \\
                                     \vec{l_m}
                                   \end{matrix}\right)\in\mathbb{R}^{m\times n}$, $p$ is a $n$-variables real polynomial with degree$\le2$ and
        $$f(\vec{\mathbf{z}})=\vec{\mathbf{r}}\left(\begin{matrix}
                                     \li(e^{\sqrt{-1}\vec{l_1}\vec{\mathbf{z}}}) \\
                                     \vdots \\
                                     \li(e^{\sqrt{-1}\vec{l_m}\vec{\mathbf{z}}})
                                   \end{matrix}\right)+p(\vec{\mathbf{z}}),$$
        then we have
        $$\im f(\vec{\mathbf{z}})=\vec{\mathbf{r}}\left(\begin{matrix}
                                     \dd(e^{\sqrt{-1}\vec{l_1}\vec{\mathbf{z}}}) \\
                                     \vdots \\
                                     \dd(e^{\sqrt{-1}\vec{l_m}\vec{\mathbf{z}}})
                                   \end{matrix}\right)+(\re\nabla f(\vec{\mathbf{z}}))\im\vec{\mathbf{z}},$$
        particularly, if $\vec{\mathbf{z_0}}$ is a critical point\footnote{$\nabla f(\vec{\mathbf{z_0}})=\vec{0}$} of $f$, then we have
        $$\im f(\vec{\mathbf{z_0}})=\vec{\mathbf{r}}\left(\begin{matrix}
                                     \dd(e^{\sqrt{-1}\vec{l_1}\vec{\mathbf{z_0}}}) \\
                                     \vdots \\
                                     \dd(e^{\sqrt{-1}\vec{l_m}\vec{\mathbf{z_0}}})
                                   \end{matrix}\right).$$
    \end{lem}

    \begin{proof}[Proof of Lemma \ref{lem2}]
        We have
        $$\im p(\vec{\mathbf{z}})=(\re\nabla p(\vec{\mathbf{z}}))\im\vec{\mathbf{z}}$$
        and for $k\in\{1,\cdots,m\}$, by definition
        $$ \dd(e^{\sqrt{-1}\vec{l_k}\vec{\mathbf{z}}})= \im(\li(e^{\sqrt{-1}\vec{l_k}\vec{\mathbf{z}}}))+\log|e^{\sqrt{-1}\vec{l_k}\vec{\mathbf{z}}}|\im\log(1-e^{\sqrt{-1}\vec{l_k}\vec{\mathbf{z}}}),$$
        so it suffices to prove
        \begin{equation}\label{e3}
          \log|e^{\sqrt{-1}\vec{l_k}\vec{\mathbf{z}}}|\im\log(1-e^{\sqrt{-1}\vec{l_k}\vec{\mathbf{z}}})+(\re\nabla\li(e^{\sqrt{-1}\vec{l_k}\vec{\mathbf{z}}}) )\im(\vec{\mathbf{z}})=0.
        \end{equation}

        We have $$\log|e^{\sqrt{-1}\vec{l_k}\vec{\mathbf{z}}}|=\re(\sqrt{-1}\vec{l_k}\vec{\mathbf{z}})=-\vec{l_k}\im(\vec{\mathbf{z}})$$ and $\nabla\li(e^{\sqrt{-1}\vec{l_k}\vec{\mathbf{z}}})=(\cdots,\frac{\D\li(e^{\sqrt{-1}\vec{l_k}\vec{\mathbf{z}}})}{\D z_i},\cdots)=(\cdots,-\sqrt{-1}\log(1-e^{\sqrt{-1}\vec{l_k}\vec{\mathbf{z}}})l_{ki},\cdots)=-\sqrt{-1}\log(1-e^{\sqrt{-1}\vec{l_k}\vec{\mathbf{z}}})\vec{l_k}$, so we have
        $$\re\nabla\li(e^{\sqrt{-1}\vec{l_k}\vec{\mathbf{z}}})=\im\log(1-e^{\sqrt{-1}\vec{l_k}\vec{\mathbf{z}}})\vec{l_k},$$
        which proves (\ref{e3}).
    \end{proof}

    The following part of this section will help us estimate as Section \ref{subsecy}.

        Define
        \begin{equation}\label{eF}
          F(t,X)=\begin{cases}
               0, & \mbox{if }X\ge0; \\
               2(2t-\pi)X, & \mbox{if }X\le0.
             \end{cases}
        \end{equation}
    we have  the following estimation
\begin{lem}\label{lem1}
    let t be a real number  with $0<t<\pi$, if $|X|\ge X_0\ge0$, We have
    \begin{equation}\label{eim0}
      |\im \li(e^{2 \sqrt{-1}(t+X \sqrt{-1})})-F(t,X)|\le m(X_0),
    \end{equation}
    where
    $$m(x)= \max_{b\in \partial D(0,e^{-2x})}\im\li(b)$$
    and $D(0,r)=\{z\in\mathbb{C}||z|\le r\}$.
\end{lem}

\begin{proof}[Proof of Lemma \ref{lem1}]
    If $X\ge0$, we have $|\im \li(e^{2 \sqrt{-1}(t+X \sqrt{-1})})-F(t,X)|=|\im \li(e^{2 \sqrt{-1}(t+X \sqrt{-1})})|$. Let $b=e^{2 \sqrt{-1}(t+X \sqrt{-1})}$, then $b\in  D(0,e^{-2x})$. $\im\li(b)$ is a harmonic function of $b$, so we have
            $$\max_{b\in  D(0,e^{-2x})}\im\li(b)=\max_{b\in \partial  D(0,e^{-2x})}\im\li(b).$$

            If $X\le0$, according to (\ref{e1}), we have
    $$\im \li(e^{2 \sqrt{-1}(t+X \sqrt{-1})})+\im \li(e^{2 \sqrt{-1}(-t-X \sqrt{-1})})=2(2t-\pi)X,$$
    so we have $|\im \li(e^{2 \sqrt{-1}(t+X \sqrt{-1})})-F(t,X)|=|\im \li(e^{2 \sqrt{-1}(-t-X \sqrt{-1})})|$, this translates to what we proved before.
\end{proof}

\begin{cor}\label{cor1.1}
    \begin{enumerate}[(1)]
        \item
            \begin{equation}\label{eim1}
                \lim _{X \to \pm\infty} |\im \li(e^{2 \sqrt{-1}(t+X \sqrt{-1})})-F(t,X)|=0.
            \end{equation}
        \item
            \begin{equation}\label{eim3}
                |\im \li(e^{2 \sqrt{-1}(t+X \sqrt{-1})})-F(t,X)|\le v_3,
            \end{equation}
            where $v_3=1.01494\cdots$ is the volume of the hyperbolic regular ideal tetrahedron.
        \item
            Let $X_0 =\frac{\log5}{4}$, if $|X|\ge X_0$ we have
            $$|\im \li(e^{2 \sqrt{-1}(t+X \sqrt{-1})})-F(t,X)|<1/2.$$
    \end{enumerate}
\end{cor}

\begin{proof}[Proof of Corollary \ref{cor1.1}]
    \begin{enumerate}[(1)]
        \item this is because of
        $$\lim_{x\to+\infty}m(x)=0.$$

        \item  $m(x)$ is a decreasing function, so we have
            $$m(x)\le m(0)=\max_{b\in \partial D(0,1)}\im\li(b)=\max_{b\in \partial D(0,1)}\dd(b)\le\max_{b\in \mathbb{H}}\dd(b)=v_3,$$
            the inequality is equal if and only if $b=\frac{1+\sqrt{-3}}{2}$.

        \item WLOG, we just need to consider the case of $X=X_0$. Let $f(t)=\im \li(e^{2 \sqrt{-1}(t+X_0 \sqrt{-1})})$. We have
    \begin{figure}[H]
        \label{f(t)}
        \centering
        \includegraphics[width=0.6\linewidth]{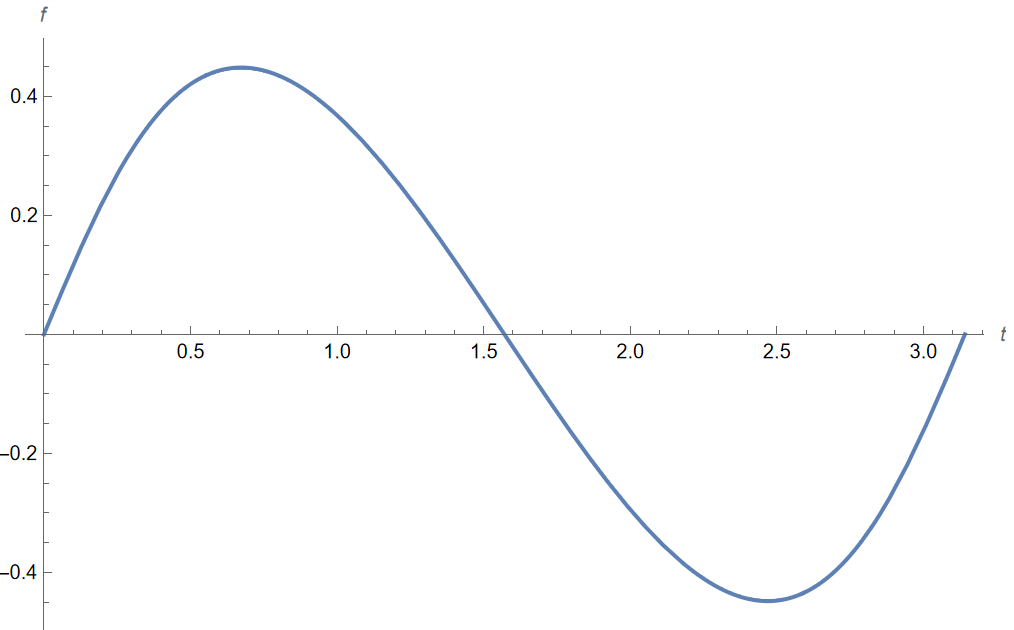}
        \caption{The function $f(t)$}
    \end{figure}
    \begin{eqnarray*}
        f'(t)&=&\im\frac{\D \li(e^{2 \sqrt{-1}z})}{\D z}|_{z=t+X_0\sqrt{-1}}\\
        &=&\im(-2\sqrt{-1}\log(1-e^{2 \sqrt{-1}(t+X_0 \sqrt{-1})}))\\
        &=&-\log(1-2e^{-2X_0}\cos 2t+e^{-4X_0})\\
        &=&\log\frac{5}{6-2\sqrt{5}\cos 2t}.
    \end{eqnarray*}
    $f'(t)$ is positive on $[0, t_0)$ and $(\pi-t_0,\pi]$ and negative on $(t_0,\pi-t_0)$, where $t_0=\frac{1}{2}\arccos \frac{1}{2\sqrt{5}}$. Hence $f(t)$ achieves its maximum value $0.448473\cdots$ at $t=t_0$, and its minimum value $-0.448473\cdots$ at $t=\pi-t_0$.
    \end{enumerate}
\end{proof}

\subsection{Quantum dilogarithm functions }
Quantum dilogarithm functions which is close related to the dilogarithm functions is an important tool to transform quantum factorials into analytic functions.

    let $r\ge3$ be an odd number, we recall the Faddev's quantum dilogarithm $\varphi_{r}(z)$ for $\{z \in \mathbb{C}|-\frac{\pi}{r}<\re z<\pi+\frac{\pi}{r}\}$, defined by the following contour integral
    \begin{equation}
    \label{Faddev-quant-dilog}
    \varphi_{r}(z)=\frac{4 \pi \sqrt{-1}}{r} \int_{\Omega} \frac{e^{(2 z-\pi) x}}{4 x \sinh (\pi x) \sinh (\frac{2 \pi x}{r})} d x,
    \end{equation}
    where the contour is $ \Omega=(-\infty,-\epsilon] \cup\{z \in \mathbb{C}||z| =\epsilon, \im z \ge 0\} \cup[\epsilon, \infty)$.

    Let $t=e^{\frac{4\pi\sqrt{-1}}{r}}$ and $(t)_{n}=\prod_{k=1}^{n}\limits\left(1-t^{k}\right) $, we recall the following lemma.
\begin{lem}[\cite{Wong-Yang41}, Lemma 2.2]
\label{lem3}
    $$(t)_{n}=2^{\epsilon_n}e^{\frac{r}{4 \pi \sqrt{-1}}(\varphi_{r}(\frac{\pi}{r})-\varphi_{r}(\frac{2\pi n+\pi}{r}-\epsilon_n \pi))},$$
    where $\epsilon_n=\begin{cases}
                                                  0, & \mbox{if } 0\le n\le\frac{r-1}{2} ,\\
                                                  1, & \mbox{if } \frac{r+1}{2}\le n\le r-1.
                                                \end{cases}$
\end{lem}

   $ \varphi_{r}(z) $ and $\li$ are closely related as follows.
\begin{lem}[\cite{52}, Proposition A.1]\label{converge}
    \

    \begin{enumerate}[(1)]
        \item For every $z$ with $0<\mathrm{Re}z<\pi,$
            \begin{equation}\label{conv1}
                \varphi_r(z)=\mathrm{Li}_2(e^{2\sqrt{-1}z})+\frac{2\pi^2e^{2\sqrt{-1}z}}{3(1-e^{2\sqrt{-1}z})}\frac{1}{r^2}+O\Big(\frac{1}{r^4}\Big).
            \end{equation}
        \item For every $z$ with $0<\mathrm{Re}z<\pi,$
            \begin{equation}\label{conv2}
                \varphi_r'(z)=-2\sqrt{-1}\log(1-e^{2\sqrt{-1}z})+O\Big(\frac{1}{r^2}\Big).
            \end{equation}
        \item
            As $r\to \infty,$ $\varphi_r(z)$ uniformly converges to $\mathrm{Li}_2(e^{2\sqrt{-1}z})$ and $\varphi_r'(z)$ uniformly converges to $-2\sqrt{-1}\log(1-e^{2\sqrt{-1}z})$ on a compact subset of $\{z\in \mathbb C\ |\ 0<\mathrm{Re}z<\pi\}.$
    \end{enumerate}
\end{lem}

Denote $\{n\}=t^\frac{n}{2}-t^{-\frac{n}{2}}$ and $\{n\}!=\prod_{k=1}^{n}\limits\{k\}$. We have
\begin{equation}\label{e4}
  \{n\}!=(-1)^n t^{-\frac{n(n+1)}{4}}(t)_n.
\end{equation}

\begin{lem}[\cite{DK} Proposition 4.1]\label{converge2}
    For $1\le n\le r-1$ we have
    $$\log |\{n\}!|=-\frac{r}{2\pi}\Lambda(\frac{2n\pi}{r})+O(\log r).$$
    Moreover, in this estimate $O(\log r)$ is uniform: there exists a constant $C_3$ independent of $n$ and $r$, such that $O(\log r)\le C_3\log r$.
\end{lem}

\begin{lem}[\cite {52} Proposition A.2, Proposition A.3]\label{lem4}
    $$\varphi_r(z)+\varphi_r(\pi-z)=2z(z-\pi)+\frac{\pi^2}{3}-\frac{2\pi^2}{3r^2}.$$
    $$\varphi_r(\frac{\pi}{r})=\frac{\pi^2}{6}+2\pi\sqrt{-1}\frac{\log r}{r}-\big(\pi^2+2\pi\sqrt{-1}\log2\big)\frac{1}{r}+\frac{2\pi^2}{3r^2}.$$
\end{lem}

\subsection{Geometry of hyperbolic ideal tetrahedron}
This section we will introduce the shape of ideal tetrahedron and show formulas of the volume of ideal tetrahedron by Bloch-Wigner functions.

    Let $\mathbb{H}^3$ be the $3$-dimension hyperbolic space. We use the half-space model, in which $\mathbb{H}^3$ is represented by $\mathbb{C}\times\mathbb{R}_+$, with the standard hyperbolic metric in which the geodesics are either vertical lines or semicircles in vertical planes with endpoints in $\mathbb{C}\times\{0\}$, and the geodesic planes are either vertical planes or else hemispheres with boundary in $\mathbb{C}\times\{0\}$.

    \begin{figure}[H]
      \centering
      \includegraphics[width=0.2\linewidth]{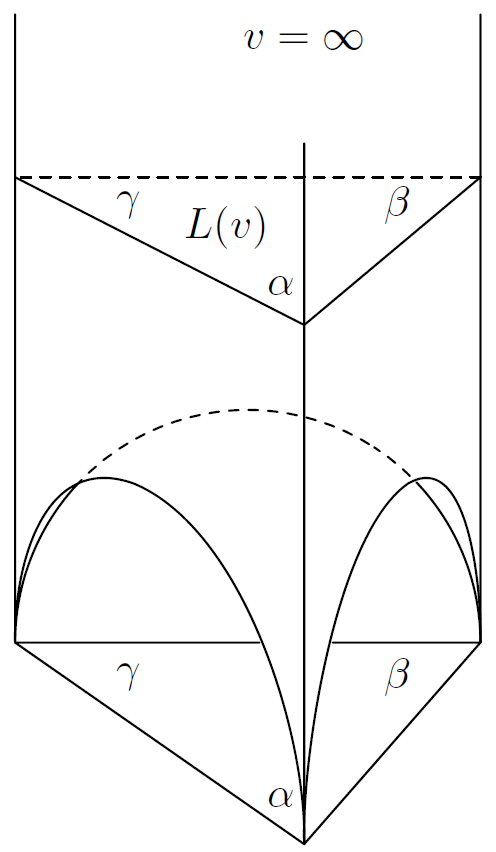}
      \caption{An ideal tetrahedron in $\mathbb{H}^3$}\label{tetrahedron}
    \end{figure}

    An ideal tetrahedron is a tetrahedron whose vertices are all in $\partial\mathbb{H}^3=\mathbb{C}\cup\{\infty\}=\mathbb{P}^1(\mathbb{C})$. We will introduce the properties of ideal tetrahedron (see \cite{ratcliffe}). Let $T$ be such a tetrahedron and let $\Sigma$ be a horosphere based at an ideal vertex $v$ of $T$ that does not meet the opposite side of $T$. Then $L(v)=\Sigma\cap T$ is a Euclidean triangle, called the link of $v$ in $T$. The (orientation preserving ) similarity class of the link $L(v)$ of a vertex $v$ of an ideal tetrahedron $T$ in $\mathbb{H}^3$ determines $T$ up to (orientation preserving ) congruence, so the shape space of ideal tetrahedron is identified with the  space of  Euclidean triangle. The dihedral angles of opposite edges of $T$ are equal. Let $T_{\alpha,\beta,\gamma}$ be an ideal tetrahedron in $\mathbb{H}^3$ with dihedral angles $\alpha,\beta,\gamma$. We have $\alpha+\beta+\gamma=\pi$.

    Let $\triangle(u,v,w)$ be a Euclidean triangle in the complex plane $\mathbb{C}$ with vertices $u,v,w$ labeled counterclockwise around $\triangle$. To each vertex of $\triangle$ we associate the ratio of the sides adjacent to the vertex in the following manner:
    $$z(u)=\frac{w-u}{v-u}, \; z(v)=\frac{u-v}{w-v}, \; z(w)=\frac{v-w}{u-w}.$$
    $z(u),z(v),z(w)$ are called the vertex invariants of $\triangle(u,v,w)$. We have $\im z(u)>0$ and $\im\log(z(u))$ is the angle of $\triangle(u,v,w)$ at $u$. Denote\footnote{We have $z''=(z')'$ and $((z')')'=z$.}
    $$z'=1-\frac{1}{z}, \; z''=\frac{1}{1-z},$$
    we have $z(w)=z(u)', \; z(v)=z(u)''$ and $z(u)z(v)z(w)=-1$.

    \begin{figure}[H]
      \centering
        \begin{minipage}[b]{0.4\linewidth}
          \includegraphics[width=1\linewidth]{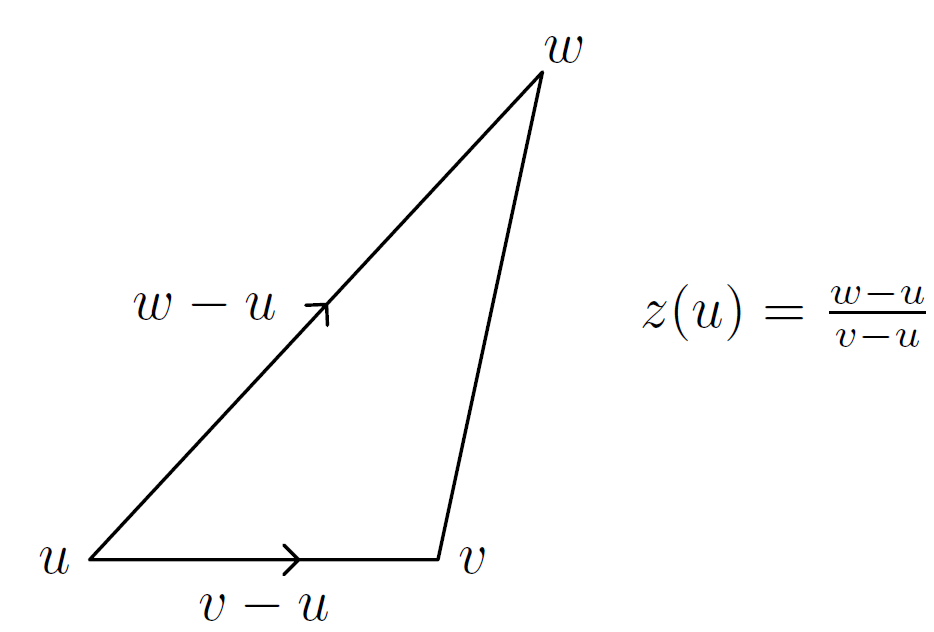}
            \label{triangle}
        \end{minipage}
        \begin{minipage}[b]{0.27\linewidth}
          \includegraphics[width=1\linewidth]{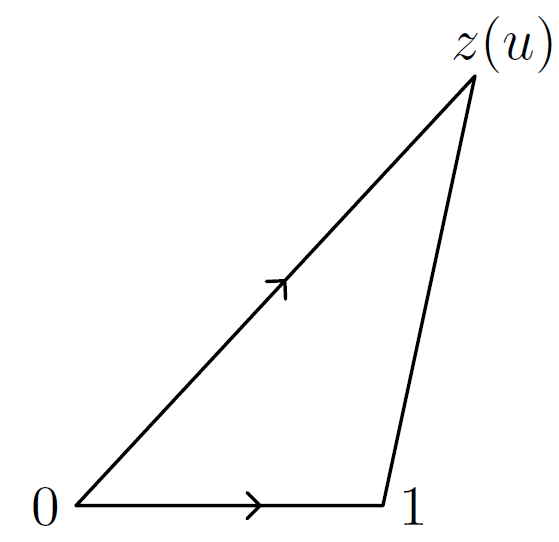}
          \label{triangle01}
        \end{minipage}
        \caption{The vertex invariants}
    \end{figure}

\begin{figure}[H]
      \centering
      \includegraphics[width=0.36\linewidth]{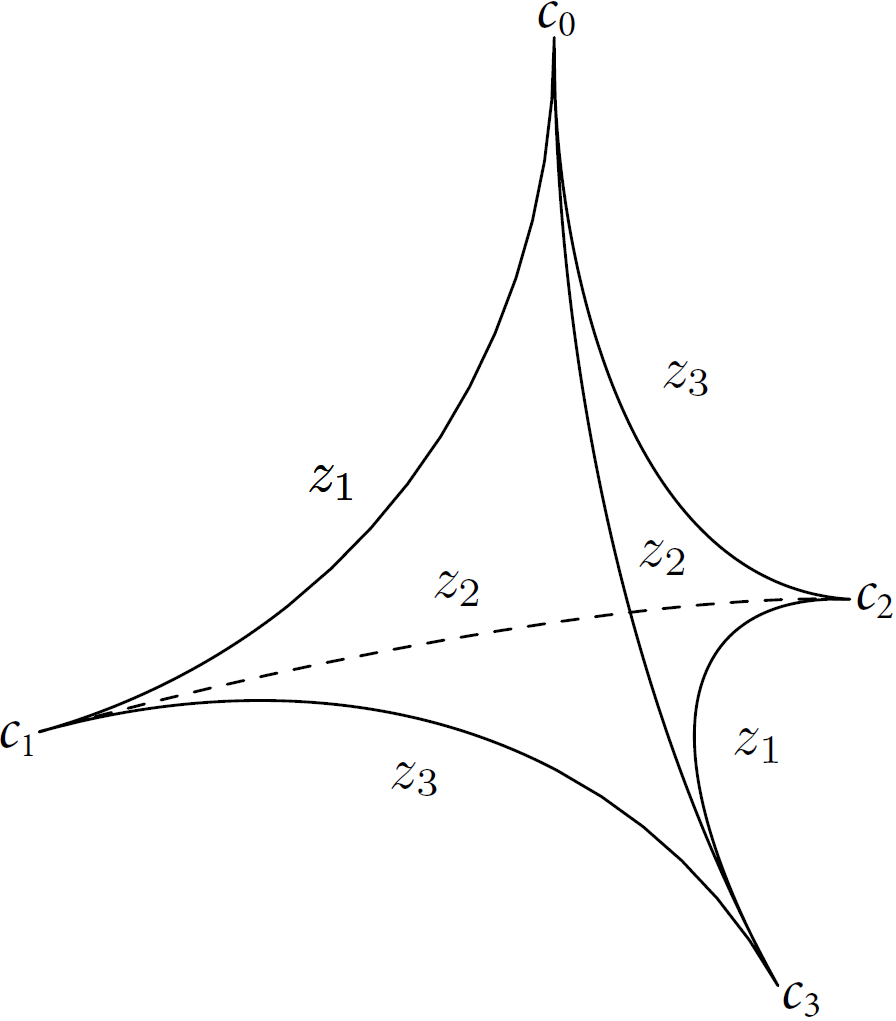}
      \caption{The edge invariants of an ideal tetrahedron}\label{ideal tetrahedron}
    \end{figure}

    For a ideal tetrahedron $T$, let $v$ be a vertex of $T$. We label the edges of $T$, incident with $v$, with the corresponding vertex invariants $z_1,z_2,z_3$ of the link of $v$. Then opposite edges of $T$ have the same label. The three parameters $z_1,z_2,z_3$ are indexed according to the right-hand rule with your thumb pointing towards a vertex of $T$ (see Figure \ref{ideal tetrahedron}). The complex parameters $z_1,z_2,z_3$ are called the edge invariants of $T$.

    The edge invariants of $T$ equal to the cross-ratios of the $4$ vertexes of $T$.
    We define the cross-ratio $[\cdot,\cdot\,;\,\cdot,\cdot]$ as below.
    $$[c_0,c_1;c_2,c_3]=\frac{(c_2-c_1)(c_3-c_0)}{(c_2-c_0)(c_3-c_1)}.$$
    Let the $4$ vertexes of $T$ be $c_0,c_1,c_2,c_3$ as showed in Figure \ref{ideal tetrahedron}. We have
    $$z_1=[c_0,c_1;c_2,c_3], \; z_2=[c_0,c_3;c_1,c_2], \; z_3=[c_0,c_2;c_3,c_1].$$

    \begin{pr}\label{prvol}
        The volume of $T_{\alpha,\beta,\gamma}$ is given by
        \begin{enumerate}
          \item  $\vol(T)=\Lambda(\alpha)+\Lambda(\beta)+\Lambda(\gamma)$
          \item $\vol(T)=\dd(z_1)=\dd(z_2)=\dd(z_3)$
        \end{enumerate}
    \end{pr}

\subsection{Dehn filling along twist knots}
This section we will discuss the topology of $\mathcal{K}_{p'}(p,q)$.

    We firstly recall the definition of Dehn-filling in $S^{3}$ along a link. Let $L$ be a link with $k$ components in $S^3$. We first remove a tubular neighborhood of $L$, resulting in $k$  boundaries $T_i$ ($1\leq i\leq k$). Then we attach a solid torus $S^1\times\mathbb{D}^2$ along each $T_i$ according to a diffeomorphism $\varphi_i:\partial(S^1\times\mathbb{D}^2)\to T_i$.  We can choose the canonical meridian $m_i$ and longitude $l_i$ for each $T_i$ according to the Rolfson's book \cite{rolfsen}: the meridian bounds a disk, the longitude and the core circle together form the boundary of a framing link with $0$-framing, and we assign a direction to the meridian and longitude such that their linking number is one. $\varphi_i$ induces the isomorphism $\varphi_{i*} :H_1(\partial(S^1\times\mathbb{D}^2))\to H_1(T_i^2)$.
    Then we can write the curve $\varphi(x\times\partial\mathbb{D}^2)$ in terms of their basis
    $$\varphi_{i*}(\{x\}\times\partial\mathbb{D}^2)=a_i m_i+b_i l_i$$
    with an ambiguity of a $\pm$ sign depending on how one wishes to orient the curve $x\times\partial\mathbb{D}^2$.  We denote the resulting closed manifold by $M_{((a_1,b_1),...,(a_k,b_k))}$.

    For convenience of exposition, we introduce generalized Dehn filling parameters (see \cite{martelli}). Let a generalised Dehn filling parameter $\vec{\mathbf{s}}=\left(s_{1}, \ldots, s_{k}\right) $ be a sequence where each $s_{i}$ $(i\in\{1,\cdots,k\})$ is either the symbol  $\infty$  or a rationally related pair of real numbers $ (p, q)=\left(k p', k q'\right)=k\left(p', q'\right)$  where  $k>0$  is real and  $\left(p', q'\right)$  are coprime integers. A generalised Dehn filling parameter  $\vec{\mathbf{s}}$  determines a Dehn filling  ${L^{\text {fill }} }$ of $ L$ as follows: for every $i$, if $s_{i}=k\left(p', q'\right)$  we fill  $T_{i}$  by killing the slope  $p' m_{i}+q' l_{i}$, while if $ s_{i}=\infty $ we do nothing. We also mark the cores of the filled solid tori with the label  $\alpha_{i}=\frac{2 \pi}{k}>0$. If the system of Thurston gluing equations
    $$p_i m_i+q_il_i=2\pi \sqrt{-1}$$
    have solutions, then it can be realized as a $3$-dimensional  hyperbolic cone-manifold with cone angle $\alpha_{i}$.

    The twist knots $\mathcal{K}_{p'}$ are defined in Figure \ref{twist knot}. For example, $\mathcal{K}_{-1}=4_{1}$, $\mathcal{K}_{0}$ is unknot, $\mathcal{K}_{1}=3_{1}$, $\mathcal{K}_{2}=5_{2}$.

    The Whitehead link $W$ and twisted knots are closely related:
    $$\mathcal{W}((p,q),(1,-p'))\cong\mathcal{K}_{p'}(p,q),$$
     where $p'$ is an integer, and $(p,q)$ are coprime integers.

\begin{figure}[H]
        \centering
        \includegraphics[width=0.8\linewidth]{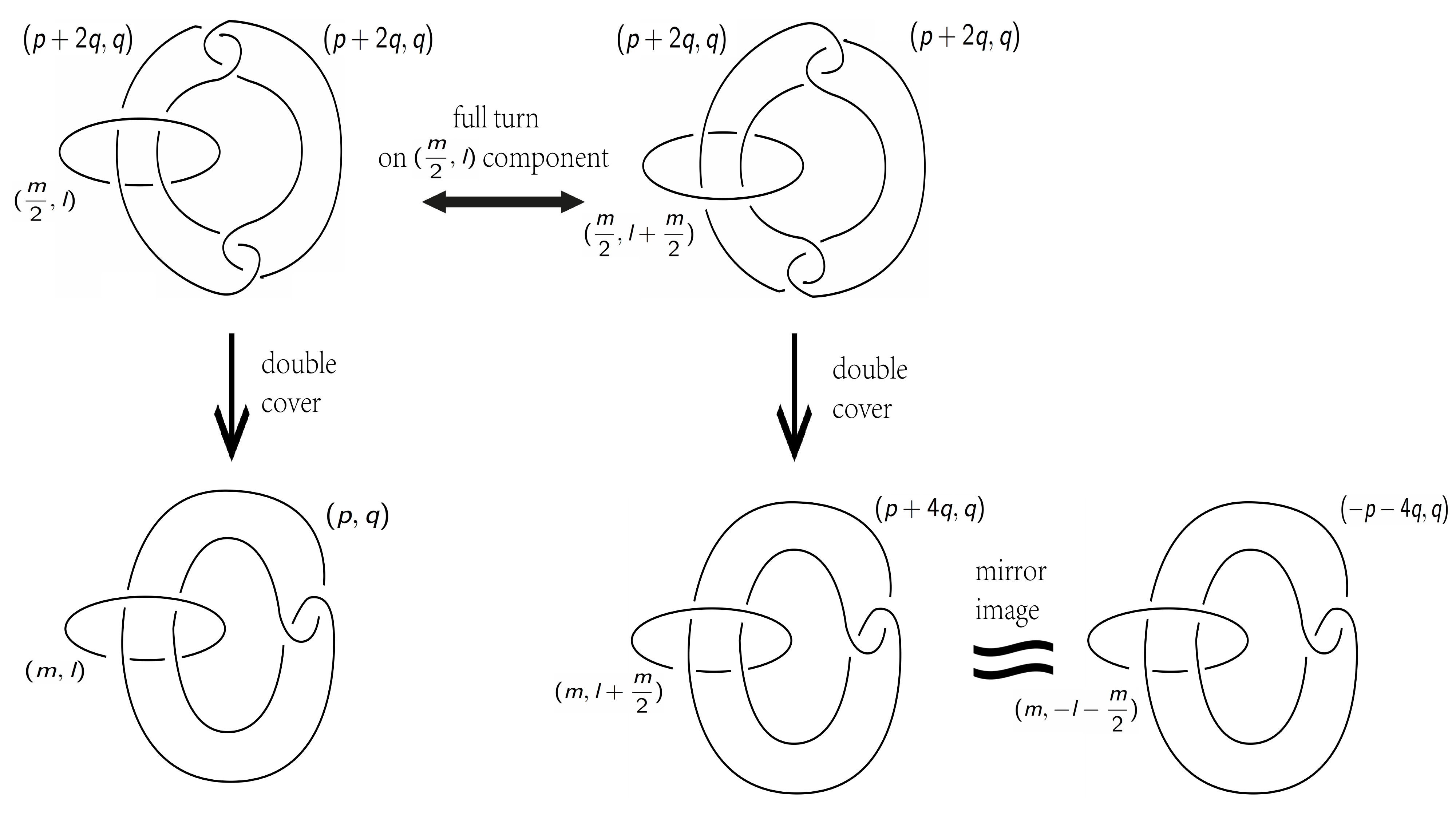}
        \caption{Sister manifolds}\label{sister manifolds}
    \end{figure}

      There is an interesting fact about the closed hyperbolic  manifolds resulting from Dehn filling on the Whitehead link: different coefficients in the Dehn filling can produce three-dimensional hyperbolic manifolds with the same volume. In \cite{HMW}, Hodgson, Meyerhoff and Weeks proved that $\mathcal{W}((m,l),(p,q))$\footnote{It's worth noting that the two component of $\mathcal{W}$ are topologically symmetric, so we have $\mathcal{W}((m,l),(p,q))=\mathcal{W}((p,q),(m,l))$}  and the mirror image of $\mathcal{W}((m,-l-\frac{m}{2}),(-p-4q,q))$ have a common degree $2$ branched covering space (see Figure \ref{sister manifolds}). So if they were hyperbolic manifold, they have same complex volume.

\subsection{Continued fractions}\label{secfrac}
Wong and Yang \cite{Wong-Yang41} discover an elegant expression of Reshetikhin-Turaev invariants of the Dehn-filled $3$-manifold along knots. Their formula involves continued fractions, we will discuss these tools in this section.

    Let
$$S=\left[\begin{array}{cc}
0 & -1 \\
1 & 0
\end{array}\right], \quad T=\left[\begin{array}{ll}
1 & 1 \\
0 & 1
\end{array}\right]$$

Given a pair of co-prime integer $(p, q)$, let  $U=\left[\begin{array}{cc}
p & -\tilde{q} \\
q & \tilde{p}
\end{array}\right]\in S L(2, \mathbb{Z}),$
where the integer $\tilde{p}$, $\tilde{q}$ are chosen so that $p\tilde{p} + q\tilde{q}= 1$. For suitable chosen $\tilde{p}$, $\tilde{q}$, we can always express\footnote{The expression of $U$ is not unique. However, each such expression of $U$ is enough for the latter use.} $U$ as the following
$$U=T^{a_{k}} S \ldots T^{a_{1}} S$$
by the standard theory from Modular groups (see \cite{Apostol} for example), with $a_i\in\mathbb{Z}$. For each $1\leq i \leq k$, denote
$$U_i=\left[\begin{array}{cc}
A_{i} & B_{i} \\
C_{i} & D_{i}
\end{array}\right]=T^{a_{i}} S \cdots T^{a_{1}} S.$$
As a result, we have that:
\begin{equation}
\left[\begin{array}{l}A_{k} \\C_{k}\end{array}\right]=\left[\begin{array}{l}p \\q\end{array}\right]
\end{equation}
and
\begin{equation}
\left[\begin{array}{l}
A_{k-1} \\
C_{k-1}
\end{array}\right]=\left[\begin{array}{c}
q \\
-p+a_{k} q
\end{array}\right] .
\end{equation}
we have the following lemma:

\begin{lem}[\cite{Wong-Yang41}]
We have the continued fraction expression \( \dfrac{A_{i}}{C_{i}} = a_{i} - \dfrac{1}{a_{i-1} - \dfrac{1}{\ddots - \dfrac{1}{a_{1}}}} \), $1\leq i\leq k$.
\end{lem}
So the rational number $\frac{p}{q}$ can be expressed as a continued fraction in the following manner:
$$\frac{p}{q} = a_{k} - \frac{1}{a_{k-1} - \frac{1}{\ddots - \frac{1}{a_{1}}}}.$$

In this section, we will recall some notations pertaining  to the  continued fraction expansion of $\frac{p}{q}$, as discussed in \cite{Wong-Yang41}. Let
$$K_{i} = \frac{(-1)^{i+1} \sum_{j=1}^{i}\limits a_{j} C_{j}}{C_{i}}.$$

Define the map \( I: \{0, \ldots, |q|-1\} \rightarrow \mathbb{Z} \) as
\[ I(s) = -C_{k-1}\left(2 s+1+K_{k-1}\right).\]

Let \( J: \{0, \ldots, |q|-1\} \rightarrow \mathbb{Q} \) be the map defined by
$$ J(s) = \frac{2 s+1}{q} + (-1)^{k} \sum_{i=1}^{k-1} \frac{(-1)^{i+1} K_{i}}{C_{i+1}}.$$
Additionally, let  $K:\{0, \ldots,|q|-1\} \rightarrow \mathbb{Q}$ be the map defined by

$$ K(s)=\frac{C_{k-1}\left(2 s+1+K_{k-1}\right)^{2}}{q}+\sum_{i=1}^{k-2} \frac{C_{i} K_{i}^{2}}{C_{i+1}}.$$

We summarize the following lemma in \cite{Wong-Yang41}  which is crucial in the computation of Reshetikhin-Turaev invariants.

\begin{lem}[\cite{Wong-Yang41} Lemma 3.3]\label{IJK}
    \

    \begin{enumerate}[(1)]
        \item $$I(s)\equiv 1-q\mod2.$$
      \item there exists a unique $s^+,s^- \in\{0, \ldots,|q|-1\}$ and $m^+,m^-\in\mathbb{Z}$ such that
        \begin{eqnarray*}
          I(s^+)&=&1-q+2 m^+ q,\\
          I(s^-)&=&-1-q+2 m^- q.
        \end{eqnarray*}

      \item We have the following congruences:
        $$J(s^+) \equiv \frac{\tilde{p}}{q} \mod \mathbb{Z}, \;J(s^-) \equiv-\frac{\tilde{p}}{q} \mod\mathbb{Z}.$$
        Furthermore
        $$J(s^+) \equiv-J(s^-) \mod 2\mathbb{Z}.$$

      \item
        $$K(s^+)=K(s^-)\equiv-\frac{\tilde{p}}{q}\mod\mathbb{Z}.$$
    \end{enumerate}
\end{lem}

Define
\begin{equation}\label{ek}
  k(s,m)=\frac{I(s)}{q}+1-2m.
\end{equation}
By Lemma \ref{IJK} we have
$$qk(s,m)\equiv1\mod2$$
and
$$k(s^+,m^+)=\frac{1}{q}, \;k(s^-,m^-)=-\frac{1}{q}.$$

\begin{lem}\label{IJK2}
    For $(s,m)$ and $(s',m')$ such that $k(s,m)+k(s',m')=0$, we have
    \begin{enumerate}[(1)]
      \item $J(s)\in\frac{\mathbb{Z}}{q}$ and $J(s)+J(s')\equiv0\mod2$.
      \item One of $s-s',m-m'$ is odd.
      \item $K(s)-K(s')\equiv 0\mod4$ and $1+m'-m+\frac{K(s)-K(s')}{4}\equiv 0\mod2$.
    \end{enumerate}
\end{lem}

\begin{proof}[Proof of Lemma \ref{IJK2}]
    \begin{enumerate}[(1)]
      \item Wong and Yang proved the following two properties in \cite{Wong-Yang41} proof of Lemma 3.3.
        \begin{enumerate}
          \item For $i\in\{2,\cdots,k\}$, define
            $$E_i=C_i\sum_{j=1}^{i-1}\frac{(-1)^{j+1}K_j}{C_{j+1}}$$
            then $E_i\in\mathbb{Z}$;

          \item $$pqJ(s)=I(s)+q(a_k(2s+1)+a_k(-1)^kE_k-(-1)^kE_{k-1})$$
        \end{enumerate}

        By (a) and $C_k=q$ we have $J(s)=\frac{2s+1+(-1)^kE_k}{q}\in\frac{\mathbb{Z}}{q}$.
        $$\text{(b)}\Leftrightarrow pJ(s)=k(s,m)+2m-1+a_k(2s+1)+a_k(-1)^kE_k-(-1)^kE_{k-1}$$
        so we have
        \begin{equation}\label{eJ1}
          p(J(s)+J(s'))\equiv0\mod2
        \end{equation}

        On the another hand
        \begin{eqnarray}
          &&J(s)-J(s')=\frac{2(s-s')}{q}\nonumber\\
          &\Rightarrow&q(J(s)-J(s'))\equiv0\mod2\nonumber\\
          &\Leftrightarrow&q(J(s)+J(s'))\equiv0\mod2\label{eJ2}
        \end{eqnarray}

        Combine (\ref{eJ1}), (\ref{eJ2}) and $(p,q)=1$, we conclude the proof.

      \item  We have
    \begin{eqnarray*}
      &&k(s,m)-k(s',m')=-\frac{2C_{k-1}(s-s')}{q}-2(m-m')\\
      &\Leftrightarrow&C_{k-1}(s-s')+q(m-m')=-\frac{q}{2}(k(s,m)-k(s',m'))\\
      &\Leftrightarrow&C_{k-1}(s-s')+q(m-m')=-qk(s,m)\\
      &\Rightarrow&C_{k-1}(s-s')+q(m-m')\equiv1\mod2
    \end{eqnarray*}

    So one of $s-s',m-m'$ is odd.
      \item We have
    \begin{eqnarray*}
      &&k(s,m)+k(s',m')=0\\
      &\Leftrightarrow&\frac{-2C_{k-1}(s+s'+1+K_{k-1})}{q}+2-2m-2m'=0\\
      &\Leftrightarrow&\frac{C_{k-1}(s+s'+1+K_{k-1})}{q}=1-m-m'\\
      &&\frac{K(s)-K(s')}{4}\\
      &=&\frac{C_{k-1}(s+s'+1+K_{k-1})(s-s')}{q}\\
      &=&(1-m-m')(s-s')
    \end{eqnarray*}

    So we have
    \begin{eqnarray*}
      &&1+m'-m+\frac{K(s)-K(s')}{4}\\
      &\equiv&1-m-m'+(1-m-m')(s-s')\\
      &=&(1-m-m')(s-s'+1)\mod 2
    \end{eqnarray*}
    By (2) we have $(1-m-m')(s-s'+1)$ is even.
    \end{enumerate}
\end{proof}

\section{Computation of the Reshetikhin-Turaev invariants}\label{sec3}
To compute the Reshetikhin-Turaev invariants, we discuss the colored Jones polynomial at first.

For a knot $K$ in $S^3$, due to the work of Habiro \cite{habiro} and Masbaum \cite{masbaum}, its colored Jones polynomial can be expressed as
$$J_N(K,q)=\sum_{n=0}^{N-1} f_K(n)\prod_{k=1}^{n}\{N+k\}\{N-k\},$$
where $f_K(n)$ is a Laurent polynomial in $\mathbb{Z}[t^{\pm\frac{1}{2}}]$ determined by $K$ and $n$. For example, $f_{4_1}(n)\equiv1$. Particularly, for the twist knot $\mathcal{K}_{p'}$, we have
\begin{pr}
    \begin{equation}\label{e4.1}
       f_{\mathcal{K}_{p'}}(n)=t^{\frac{n(n+3)}{4}}\sum_{l=0}^n(-1)^lt^{p'(l^2+l)}\{2l+1\}\frac{\{n\}!}{\{n-l\}!\{n+l+1\}!}
    \end{equation}
\end{pr}

\begin{rem}
    It's not obvious that $f_{\mathcal{K}_{p'}}(n)$ is a Laurent polynomial. \cite{masbaum} gives two forms of $f_{\mathcal{K}_{p'}}(n)$, the one is as before, the another one is a complex multiple summation expressions which a Laurent polynomial.

    This fact makes the formulation of the Proposition \ref{pr1} much simpler, which greatly simplifies the calculation. If we ignore this fact, we have to deal with the terms in (\ref{e5}) whose both the numerator and denominator are $0$. This means that the summing term (\ref{e5}) will be a lot more, and the region $D$ showed in Section \ref{subsec4.1} will be bigger. In this case $D$ will have two extra parts which corresponding to the large volume terms, whose volume may be $\vol(\mathcal{W})=3.66\cdots$ which is larger than $\vol(\mathcal{K}_{p'}(p,q))$!

    However, the terms corresponding to the two extra parts will completely cancel out each other, this is because of the symmetry of the Fourier coefficients which is similar to Section \ref{subsec4.3}. This inspires us to conjecture that Habiro's cyclotomic expansions \cite{habiro} implies some symmetry of the Fourier coefficients in a sence.
\end{rem}

In knot theory, it is possible to transform any rational Dehn-filling procedure performed on a knot into an integer surgery on a framed link $L$ consisting of $k$ components. This is visually represented in the following figure(we choose $K=\mathcal{K}_{p'}(p,q)$), where the surgeries coefficients  $a_1, \ldots, a_k$ are integer number.
\begin{figure}[H]
    \centering
    \includegraphics[width=0.6\linewidth]{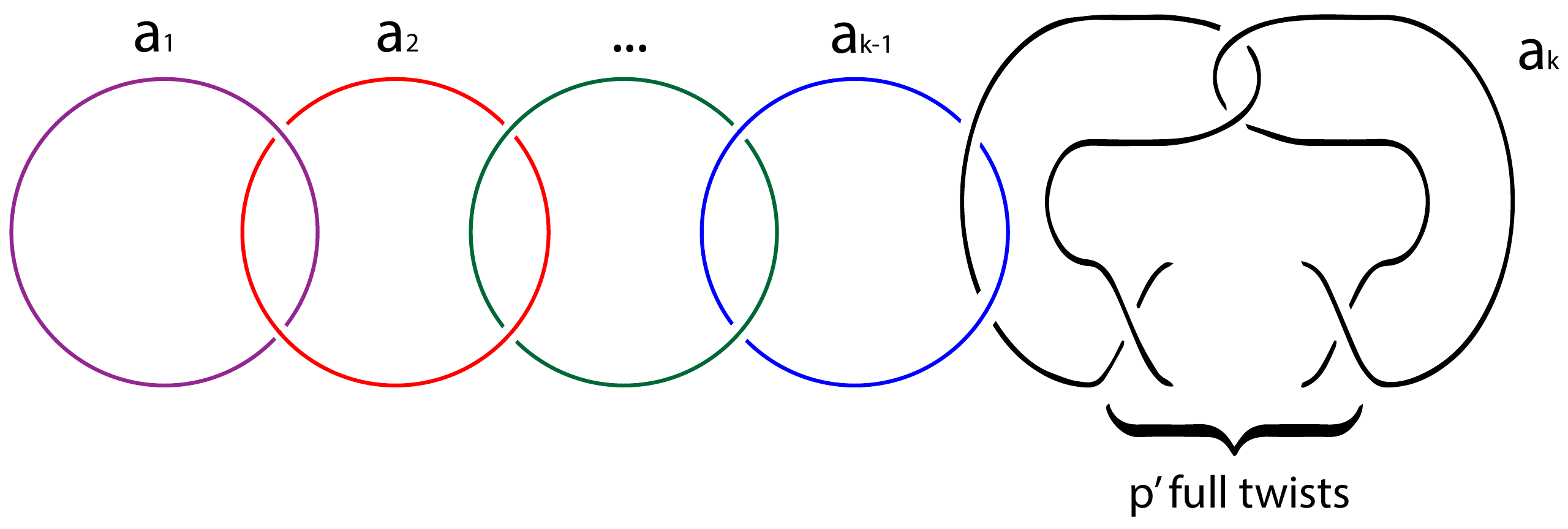}
    \label{figL}\caption{Surgery on $L$}
\end{figure}

Wong and Yang \cite{Wong-Yang41} simplified the Reshetikhin-Turaev invariants for general $K\subset S^3$ as follows, which shows the relation between the Reshetikhin-Turaev invariants and the (Habiro's cyclotomic expansions of) colored Jones polynomial.

\begin{pr}\label{pr1}
    $$\rt_r(K(p,q))=\kappa_r\sum_{s=0}^{|q|-1}\sum_{m=1-\frac{r}{2}}^{\frac{r}{2}-1}\sum_{n=|m|}^{\frac{r}{2}-1}\sin(\frac{x}{q}-J(s)\pi)\left(e^{-\sqrt{-1}x+\frac{r}{4\pi\sqrt{-1}}(-\frac{p}{q}x^2+\frac{2\pi I(s)}{q}x+4xy+K(s)\pi^2)}\right)$$
    $$\cdot\frac{(t)_{r-m-n-1}}{(t)_{n-m}}f_{K}(\frac{r}{2}-n-1)$$
    where
    $$\kappa_r=\frac{(-1)^{\frac{3(k+1)}{4}+\sum_{i=1}^k\limits a_i}e^{\frac{\pi\sqrt{-1}}{r}\big(3\sigma(L)-\sum_{i=1}^k\limits a_i-\sum_{i=2}^k\limits \frac{1}{C_{i-1}C_i}\big)+\frac{\pi\sqrt{-1}r}{4}\big(\sigma(L)+3a_k\big)}}{2r\sqrt{q}}$$
    $m,n$ are half-integers, $\sigma(L)$ is the signature of linking matrix of the link $L$ and
    $$x=\frac{2\pi m}{r},y=\frac{2\pi n}{r}$$
\end{pr}
\begin{rem}
Proposition \ref{pr1} is in fact the formula (3.9) of Propositon 3.4 in \cite{Wong-Yang41}. The $f_K(\cdots)$ in Proposition \ref{pr1} corresponds to $f_{4_1}=1$ in \cite{Wong-Yang41}.
\end{rem}

By (\ref{e4}) and (\ref{e4.1}), we have
$$f_{\mathcal{K}_{p'}}(n)=\sum_{l=0}^n(-1)^{n+l+1}t^{(p'+\frac{1}{2})(l^2+l)+\frac{n^2+3n+1}{2}}\{2l+1\}\frac{(t)_n}{(t)_{n-l}(t)_{n+l+1}}$$
and
\begin{eqnarray*}
  f_{\mathcal{K}_{p'}}(\frac{r}{2}-n-1)&=&\sum_{l=0}^{\frac{r}{2}-n-1}(-1)^{\frac{r}{2}-n+l}t^{(p'+\frac{1}{2})(l^2+l)+\frac{(\frac{r}{2}-n-1)^2+3(\frac{r}{2}-n-1)+1}{2}}\{2l+1\}\frac{(t)_{\frac{r}{2}-n-1}}{(t)_{\frac{r}{2}-n-l-1}(t)_{\frac{r}{2}-n+l}}\\
  &\stackrel{l'=l+\frac{1}{2}}{=}&\sum_{l'=\frac{1}{2}}^{\frac{r-1}{2}-n}(-1)^{\frac{r-1}{2}-n+l'}t^{(p'+\frac{1}{2})(l'^2-\frac{1}{4})+\frac{n^2}{2}-\frac{(r+1)n}{2}+\frac{r^2}{8}+\frac{r}{4}-\frac{1}{2}}\{2l'\}\frac{(t)_{\frac{r}{2}-n-1}}{(t)_{\frac{r-1}{2}-n-l'}(t)_{\frac{r-1}{2}-n+l'}}
\end{eqnarray*}

If $t=e^{\frac{4\pi\sqrt{-1}}{r}}$, let $z=\frac{2\pi l'}{r}$, we have
$$f_{\mathcal{K}_{p'}}(\frac{r}{2}-n-1)=2(-1)^{\frac{r+1}{2}}e^{\frac{\pi\sqrt{-1}}{r}(-p'-\frac{5}{2})}\sum_{l'=\frac{1}{2}}^{\frac{r-1}{2}-n}\sin (2z)\left(e^{-\sqrt{-1}y+\frac{r}{4\pi\sqrt{-1}}(-(4p'+2)z^2-2y^2-2\pi y-2\pi z)}\right)$$
$$\cdot\frac{(t)_{\frac{r}{2}-n-1}}{(t)_{\frac{r-1}{2}-n-l'}(t)_{\frac{r-1}{2}-n+l'}}$$

Now we can start the calculation of $\rt_r(\mathcal{K}_{p'}(p,q))$. Substitute the above formula into Proposition \ref{pr1}, we have
    \begin{eqnarray}\label{e5}
      \rt_r(\mathcal{K}_{p'}(p,q))=c_r\sum_{s=0}^{|q|-1}\sum_{m=1-\frac{r}{2}}^{\frac{r}{2}-1}\sum_{n=|m|}^{\frac{r}{2}-1}\sum_{l'=\frac{1}{2}}^{\frac{r-1}{2}-n}\sin(\frac{x}{q}-J(s)\pi)\sin(2z)e^{-\sqrt{-1}(x+y)} \\
    \cdot e^{\frac{r}{4\pi\sqrt{-1}}(-\frac{p}{q}x^2+4xy-2y^2-(4p'+2)z^2+\frac{2\pi I(s)}{q}x-2\pi y-2\pi z+K(s)\pi^2)}\frac{(t)_{r-m-n-1}(t)_{\frac{r}{2}-n-1}}{(t)_{n-m}(t)_{\frac{r-1}{2}-n-l'}(t)_{\frac{r-1}{2}-n+l'}}\nonumber
    \end{eqnarray}
where
\begin{equation}\label{e6}
  c_r=\frac{(-1)^{\frac{3(k+1)}{4}+\sum_{i=1}^k\limits a_i+\frac{1}{2}}e^{\frac{\pi\sqrt{-1}}{r}\big(3\sigma(L)-\sum_{i=1}^k\limits a_i-\sum_{i=2}^k\limits \frac{1}{C_{i-1}C_i}-p'-\frac{5}{2}\big)+\frac{\pi\sqrt{-1}r}{4}\big(\sigma(L)+3a_k+2\big)}}{r\sqrt{q}}
\end{equation}
We have
$$x\in(-\pi,\pi),y\in[|x|,\pi),z\in(0,\pi-y)$$
and by Lemma \ref{lem3},
$$\frac{(t)_{r-m-n-1}(t)_{\frac{r}{2}-n-1}}{(t)_{n-m}(t)_{\frac{r-1}{2}-n-l'}(t)_{\frac{r-1}{2}-n+l'}}=2^{\epsilon_1+\epsilon_2+\epsilon_3}$$
$$\cdot e^{\frac{r}{4\pi\sqrt{-1}}(\varphi_r(y-x+\frac{\pi}{r}-\epsilon_1\pi)+\varphi_r(\pi-y-z)+\varphi_r(\pi-y+z-\epsilon_2\pi)-\varphi_r(2\pi-x-y-\frac{\pi}{r}-\epsilon_3\pi)-\varphi_r(\pi-y-\frac{\pi}{r})-\varphi_r(\frac{\pi}{r}))}$$
where $\epsilon_1=\begin{cases}
    0, & \mbox{if }y-x\in[0,\pi)  \\
    1, & \mbox{if }y-x>\pi
    \end{cases}$, $\epsilon_2=\begin{cases}
                    0, & \mbox{if }y-z\ge0  \\
                    1, & \mbox{if }y-z<0
                  \end{cases}$, $\epsilon_3=\begin{cases}
                    0, & \mbox{if }x+y>\pi \\
                    1, & \mbox{if }x+y\in[0,\pi)
                  \end{cases}$.

Let's define the potential function $V_r(x,y,z,s)$ (we will define potential function $V_r(x,y,z,s,m,n,l)$ and $V(x,y,z,s,m,n,l)$ later). the "potential" here means that it is the key point of the proof of volume conjecture and we hope the potential function is closed related to the hyperbolic volume of $3$-manifold (we will show this in Section \ref{sec5}).
\begin{eqnarray}
    \label{e7}
    V_r(x,y,z,s)=\varphi_r(y-x+\frac{\pi}{r}-\epsilon_1\pi)+\varphi_r(\pi-y-z)+\varphi_r(\pi-y+z-\epsilon_2\pi)\\
    -\varphi_r(2\pi-x-y-\frac{\pi}{r}-\epsilon_3\pi)-\varphi_r(\pi-y-\frac{\pi}{r})-\varphi_r(\frac{\pi}{r}) -\frac{p}{q}x^2\nonumber\\
   +4xy-2y^2-(4p'+2)z^2+\frac{2\pi I(s)}{q}x-2\pi y-2\pi z+K(s)\pi^2\nonumber
\end{eqnarray}
 and $\epsilon(x,y,z)=2^{\epsilon_1+\epsilon_2+\epsilon_3}$, we have
\begin{thm}\label{thm2}
        \begin{equation}\label{e8}
          \rt_r(\mathcal{K}_{p'}(p,q))=c_r\sum_{s=0}^{|q|-1}\sum_{m=1-\frac{r}{2}}^{\frac{r}{2}-1}\sum_{n=|m|}^{\frac{r}{2}-1}\sum_{l'=\frac{1}{2}}^{\frac{r-1}{2}-n}g_r(s,x,y,z)
        \end{equation}
        where

        \begin{equation}\label{e8.1}
          g_r(s,x,y,z)=\epsilon(x,y,z)\sin(\frac{x}{q}-J(s)\pi)\sin(2z)e^{-\sqrt{-1}(x+y)+\frac{r}{4\pi\sqrt{-1}}V_r(x,y,z,s)}
        \end{equation}
        and $c_r,x,y,z,\epsilon_1,\epsilon_2,\epsilon_3,I(s),J(s),K(s),V_r$ are defined as before.
\end{thm}

\section{Discrete Fourier transform}\label{sec4}
    In this section, we use the discrete Fourier transform method to express $\mathrm{RT}_r(M)$ as a sum of integrals, which facilitate our estimation (in Section \ref{sec6}).

    In Section \ref{subsec4.1} we study the properties of the potential function $V_r$ and $V$.

    In Section \ref{subsec4.2} we use the Poisson summation formula to express $\mathrm{RT}_r(M)$ as a sum of integrals.

    In Section \ref{subsec4.3} we discover the symmetry of the Fourier coefficients, which is a generalization of the \textbf{Big Cancellation} as showed in \cite{Chen-Zhu-1}. As a corollary, the troublesome terms cancel out directly.

\subsection{Potential function $V_r$ and $V$}\label{subsec4.1}

We will discuss two regions $D$ and $D_\epsilon$ at first. This process simplifies the definition of potential functions. $D$ is the whole area over which we sum or integrate. The definition of $D_\epsilon$ is technical, which is to satisfy Lemma \ref{lemim<3.5} and Lemma \ref{lemyconvex}.

\begin{figure}[H]
        \centering
        \includegraphics[width=0.6\linewidth]{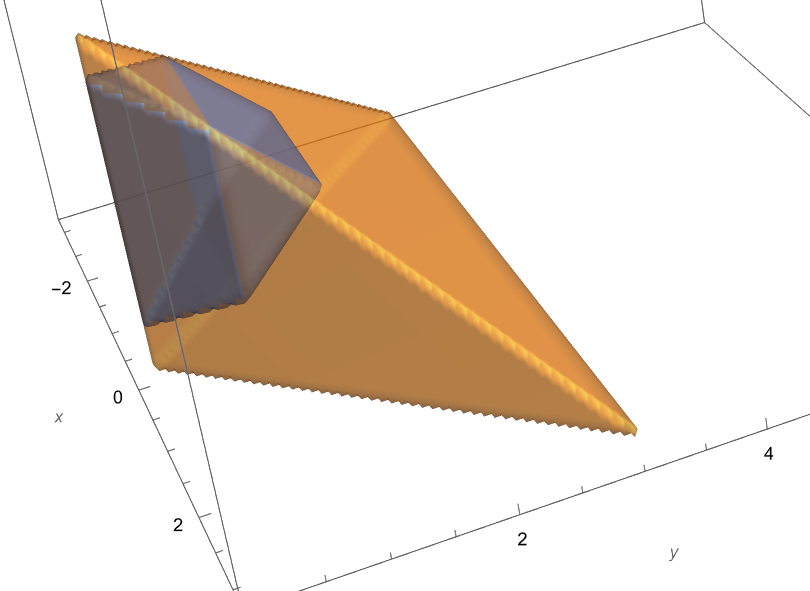}
        \caption{The orange region is $D$, and the blue region is $D_0$}\label{D&D1}
    \end{figure}

For $\epsilon\ge0$, let
$$\begin{aligned}
    D=\{&(x,y,z)\in\mathbb{R}^3|y\in[|x|,\pi),z\in(0,\pi-y)\}\\
    D_\epsilon=\{&(x,y,z)\in\mathbb{R}^3|-\frac{\pi}{4}+\epsilon<x<\frac{\pi}{4}-\epsilon,y-x>\epsilon,x+y>\epsilon,z-y>\epsilon,\\
    &y+z<\pi-\epsilon,y+2z<\frac{7}{4}\pi-\epsilon,2z-y>\frac{\pi}{4}-\epsilon\}
\end{aligned}$$
Then for $\epsilon>0$. we have $(x,y,z)\in D$ and $D_\epsilon\subset D_0\subset D$. Figure \ref{D&D1} shows the shape of $D$ and $D_0$.

 If $(x,y,z)\in D_0$, we have $\epsilon_1=0,\epsilon_2=\epsilon_3=1,$
 $$\epsilon(x,y,z)=4$$
 and
\begin{eqnarray}
    \label{e9.9}
    &&V_r(x,y,z,s)=\varphi_r(y-x+\frac{\pi}{r})+\varphi_r(\pi-y-z)+\varphi_r(z-y) -\varphi_r(\pi-x-y-\frac{\pi}{r})\\
   &&-\varphi_r(\pi-y-\frac{\pi}{r})-\varphi_r(\frac{\pi}{r}) -\frac{p}{q}x^2+4xy-2y^2-(4p'+2)z^2+\frac{2\pi I(s)}{q}x-2\pi y-2\pi z+K(s)\pi^2\nonumber
\end{eqnarray}

By Lemma \ref{lem4}, we have
\begin{eqnarray*}
  &&-\varphi_r(\pi-x-y-\frac{\pi}{r}) \\
  &=&\varphi_r(x+y+\frac{\pi}{r})+2(x+y+\frac{\pi}{r})(\pi-x-y-\frac{\pi}{r})-\frac{\pi^2}{3}+\frac{2\pi^2}{3r^2}\\
  &=&\varphi_r(x+y+\frac{\pi}{r})-2(x+y)^2+2\pi(x+y)-\frac{\pi^2}{3}-\frac{4\pi}{r}(x+y)+\frac{2\pi^2}{r}-\frac{\pi^2}{3r^2}
\end{eqnarray*}
so we have
\begin{eqnarray}\label{e10}
    &&V_r(x,y,z,s)=\varphi_r(y+x+\frac{\pi}{r})+\varphi_r(y-x+\frac{\pi}{r})+\varphi_r(-y+z)+\varphi_r(-y+\pi-z) \\ &&-\varphi_r(\pi-y-\frac{\pi}{r})-\frac{p+2q}{q}x^2-4y^2-(4p'+2)z^2+2\pi(\frac{ I(s)}{q}+1)x-2\pi z+K(s)\pi^2\nonumber\\
    &&-\varphi_r(\frac{\pi}{r}) -\frac{\pi^2}{3}-\frac{4\pi}{r}(x+y)+\frac{2\pi^2}{r}-\frac{\pi^2}{3r^2}\nonumber
\end{eqnarray}

Let
\begin{eqnarray}
  V_r(x,y,z,s,m,n,l)&=&V_r(x,y,z,s)-4\pi m x-4\pi n y-4\pi l z\label{e19}
\end{eqnarray}

\begin{de}\label{defV}
    The potential function $V$ is define as below. For $x,y,z$ satisfied $\re(x,y,z)\in D$,
    \begin{eqnarray*}
      &&V(x,y,z,s,m,n,l)=\li(e^{2\sqrt{-1}(y+x)})+\li(e^{2\sqrt{-1}(y-x)})+\li(e^{2\sqrt{-1}(-y+z)})+\li(e^{2\sqrt{-1}(-y-z)})\\
      &&-\li(e^{-2\sqrt{-1}y})-\frac{p+2q}{q}x^2-4y^2-(4p'+2)z^2+2\pi k(s,m)x-4\pi ny-2\pi(2l+1) z+K(s)\pi^2-\frac{\pi^2}{2}\nonumber
    \end{eqnarray*}
where $k(s,m)=\frac{I(s)}{q}+1-2m$.
\end{de}

 It is worth mentioning that $V$ is an analytic function for $\re(x,y,z)\in D_0$.

By simple quadratic function calculation, we have the following identical equation which shows the symmetries of the potential functions.
\begin{pr}\label{pr3}
    For $(s,m)$ and $(s',m')$ such that $k(s,m)+k(s',m')=0$, we have

    \begin{enumerate}[(1)]
        \item
            \begin{equation}\label{ec1}
              \overline{V(x,y,z,s,m,n,l)}=V(-\overline{x},-\overline{y},-\overline{z},s',m',-n,-l-1)+(K(s)-K(s'))\pi^2
            \end{equation}
      \item
        \begin{eqnarray}\label{e20}
            V_r(x,y,z,s,m,n,l)&=&V_r(-x,y,z,s',m',n,l)-\frac{8\pi}{r}x+(K(s)-K(s'))\pi^2\\
            V(x,y,z,s,m,n,l)&=&V(-x,y,z,s',m',n,l)+(K(s)-K(s'))\pi^2\nonumber
        \end{eqnarray}
      \item
        \begin{eqnarray}\label{e21}
          V_r(x,y,z,s,m,n,l)&=&V_r(x,y,\pi-z,s,m,n,-2p'-2-l)-4(p'+l+1)\pi^2\\
          V(x,y,z,s,m,n,l)&=&V(x,y,\pi-z,s,m,n,-2p'-2-l)-4(p'+l+1)\pi^2\nonumber
        \end{eqnarray}
      \item
            \begin{equation}\label{eVz}
                V(x,y,z,s,m,n,l)=V(x,y,\pi+z,s,m,n,l-2p'-1)+4(l-p')\pi^2
            \end{equation}
    \end{enumerate}
\end{pr}

The following lemma shows that $V_r$ converges uniformly to $V$.

\begin{lem}\label{convergeV}
    For $(x,y,z)\in D_{\frac{\epsilon}{2}}$, we have
    \begin{eqnarray*}
      V_r(x,y,z,s)=&&V(x,y,z,s)-(\log(1-e^{2\sqrt{-1}(y+x)})+\log(1-e^{2\sqrt{-1}(y-x)})+\log(1-e^{-2\sqrt{-1}y})\\
      &&-2\sqrt{-1}(x+y)+\log r+\frac{3\pi\sqrt{-1}}{2}-\log2)\frac{2\pi\sqrt{-1}}{r}+\frac{v_r(x,y,z)}{r^2}
    \end{eqnarray*}
    with $|v_r(x,y,z)|$ bounded from above by a constant independent of $r$.
\end{lem}

\begin{proof}[Proof of Lemma \ref{convergeV}]
  By Lemma \ref{converge} and Lemma \ref{lem4}, we have
  \begin{eqnarray*}
    \varphi_r(y+x+\frac{\pi}{r})&=&\varphi_r(y+x)+\varphi_r'(y+x)\frac{\pi}{r}+O(\frac{1}{r^2})\\
    &=&\li(e^{2\sqrt{-1}(y+x)})-\frac{2\pi\sqrt{-1}}{r}\log(1-e^{2\sqrt{-1}(y+x)})+O(\frac{1}{r^2})\\
    \varphi_r(y-x+\frac{\pi}{r})&=&\varphi_r(y-x)+\varphi_r'(y-x)\frac{\pi}{r}+O(\frac{1}{r^2})\\
    &=&\li(e^{2\sqrt{-1}(y-x)})-\frac{2\pi\sqrt{-1}}{r}\log(1-e^{2\sqrt{-1}(y-x)})+O(\frac{1}{r^2})\\
    \varphi_r(\pi-y-z)&=&\li(e^{-2\sqrt{-1}(y+z)})+O(\frac{1}{r^2})\\
    \varphi_r(z-y)&=&\li(e^{2\sqrt{-1}(z-y)})+O(\frac{1}{r^2})\\
    \varphi_r(\pi-y-\frac{\pi}{r})&=&\varphi_r(\pi-y)-\varphi_r'(\pi-y)\frac{\pi}{r}+O(\frac{1}{r^2})\\
    &=&\li(e^{-2\sqrt{-1}y})+\frac{2\pi\sqrt{-1}}{r}\log(1-e^{-2\sqrt{-1}y})+O(\frac{1}{r^2})\\
    \varphi_r(\frac{\pi}{r})&=&\frac{\pi^2}{6}+2\pi\sqrt{-1}\frac{\log r}{r}-(\pi^2+2\pi\sqrt{-1}\log2)\frac{1}{r}+O(\frac{1}{r^2})
  \end{eqnarray*}
  Substitute the above formulas, we are done.
\end{proof}

The following parts of this section shows that the terms outside $D_\epsilon$ are neglectable. This process simplifies our analysis, because of the terms in $D_\epsilon$ have convexity (Lemma \ref{lemyconvex}).

\begin{lem}\label{lemim<3.5}
     There exists a sufficiently small $\epsilon>0$, such that for $(x,y,z)\in \overline{D\backslash D_\epsilon}$, we have
    $$\im V(x,y,z,s,m,n,l)<3.5$$
\end{lem}

\begin{figure}[H]
  \centering
  \includegraphics[width=0.4\textwidth]{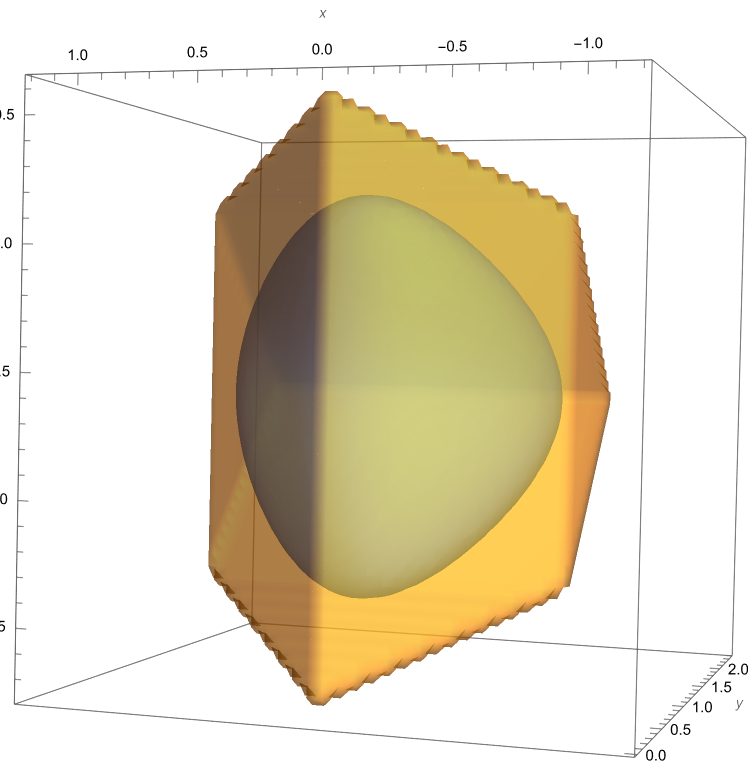}
  \caption{$D_0$ and the region $\im V(x,y,z)>3.5$}\label{D0imV}
\end{figure}
The above picture explains Lemma \ref{lemim<3.5}. The orange region is $D_0$ and the blue region is $\{(x,y,z)\in D|\im V(x,y,z,s,m,n,l)>3.5\}$.

The proof of Lemma \ref{lemim<3.5} is in Appendix.

\begin{thm}\label{thm4}
     There exists a sufficiently small $\epsilon,\delta>0$, such that for $(x,y,z)\in D\backslash D_\epsilon$ and $r$ sufficiently large,
    $$|g_r(s,x,y,z)|<e^{\frac{r}{4\pi}\cdot(3.5-\delta)}$$
\end{thm}

\begin{proof}[Proof of Theorem \ref{thm4}]
  We have
  \begin{eqnarray*}
    &&|g_r(s,x,y,z)|\\
    &\stackrel{(\ref{e5})}{=}&|\sin(\frac{x}{q}-J(s)\pi)\sin(2z)\frac{(t)_{r-m-n-1}(t)_{\frac{r}{2}-n-1}}{(t)_{n-m}(t)_{\frac{r-1}{2}-n-l'}(t)_{\frac{r-1}{2}-n+l'}}|\\
    &\le&|\frac{(t)_{r-m-n-1}(t)_{\frac{r}{2}-n-1}}{(t)_{n-m}(t)_{\frac{r-1}{2}-n-l'}(t)_{\frac{r-1}{2}-n+l'}}|\\
    &=&|\frac{\{r-m-n-1\}!\{\frac{r}{2}-n-1\}!}{\{n-m\}!\{\frac{r-1}{2}-n-l'\}!\{\frac{r-1}{2}-n+l'\}!}|
  \end{eqnarray*}

  By Lemma \ref{converge2} we have
  \begin{eqnarray*}
    &&\log|\frac{\{r-m-n-1\}!\{\frac{r}{2}-n-1\}!}{\{n-m\}!\{\frac{r-1}{2}-n-l'\}!\{\frac{r-1}{2}-n+l'\}!}|\\
    &=&-\frac{r}{2\pi}(\Lambda(-x-y-\frac{2\pi}{r})+\Lambda(\pi-y-\frac{2\pi}{r})-\Lambda(y-x)-\Lambda(-\frac{\pi}{r}-y-z)-\Lambda(-\frac{\pi}{r}-y+z))+O(\log r)\\
    &\stackrel{Lemma\ \ref{lemLambda}}{=}&-\frac{r}{2\pi}(\Lambda(-x-y)+\Lambda(-y)-\Lambda(y-x)-\Lambda(-y-z)-\Lambda(-y+z))+O(\log r)\\
    &\stackrel{Proposition\ \ref{prD}}{=}&-\frac{r}{4\pi}(\dd(e^{2\sqrt{-1}(-x-y)})+\dd(e^{2\sqrt{-1}(-y)})-\dd(e^{2\sqrt{-1}(y-x)})-\dd(e^{2\sqrt{-1}(-y-z)})-\dd(e^{2\sqrt{-1}(-y+z)}))\\
    &&+O(\log r)\\
    &\stackrel{Proposition\ \ref{prD}}{=}&\frac{r}{4\pi}(\dd(e^{2\sqrt{-1}(y+x)})+\dd(e^{2\sqrt{-1}(y-x)})+\dd(e^{2\sqrt{-1}(-y+z)})+\dd(e^{2\sqrt{-1}(-y-z)})-\dd(e^{-2\sqrt{-1}y}))\\
    &&+O(\log r)\\
    &\stackrel{Proposition\ \ref{prD}}{=}&\frac{r}{4\pi}\im V(x,y,z,s,m,n,l)+O(\log r)
  \end{eqnarray*}
  Where the $O(\log r)$ here is uniform. By Lemma \ref{lemim<3.5} we have for $r$ sufficiently large, there exist a $\delta>0$ such that
  $$\frac{r}{4\pi}\im V(x,y,z,s,m,n,l)+O(\log r)<\frac{r}{4\pi}(3.5-\delta)$$
  We conclude the proof.

\end{proof}

\subsection{Poisson summation formula}\label{subsec4.2}

 We will use the Poisson summation formula to express $\mathrm{RT}_r(M)$ as a sum of integrals in this section.

We can choose a $C^{\infty}$-smooth bump function $\psi$  on $ \mathbb{R}^{3}$ such that
$$\begin{cases}
  \psi(x, y,z)=1, & \mbox{if }(x, y,z)\in\overline{D_{\frac{\epsilon}{2}}}\\
  0<\psi(x, y,z) <1, & \mbox{if }(x, y,z)\in D_{\frac{\epsilon}{2}}\backslash\overline{D_{\epsilon}}\\
  \psi(x, y,z)=0, & \mbox{if }(x, y,z)\in\mathbb{R}^3\backslash\overline{D_{\frac{\epsilon}{2}}}.
\end{cases}$$
and $\psi(x,y,z)=\psi(-x,y,z),\psi(x,y,\pi-z)=\psi(x,y,z)$, since the region $D_{\epsilon}$ is symmetric with respect to the plane $x=0$ and $z=\frac{\pi}{2}$.

Let$$f_r(s,x,y,z)= \psi(x,y,z)g_r(s,x,y,z)$$

 Then by the Theorem \ref{thm4}
 \begin{equation}\label{e9}
    RT_r({\mathcal{K}_{p'}(p,q)})=c_r\sum_{s=0}^{|q|-1}\sum_{(m,n,l')\in  (\mathbb{Z} +\frac{1}{2})^3}f_r(s,x,y,z)+o(e^{\frac{r}{4\pi}\cdot3.5})
 \end{equation}
we recall the Possion summation Formula (see \cite{stein}), for any function $f$ in Schwarz space we have
$$\sum_{\left(m_{1}, \ldots, m_{k}\right) \in \mathbb{Z}^{k}} f\left(m_{1}, \ldots, m_{k}\right)=\sum_{\left(n_{1}, \ldots, n_{k}\right) \in \mathbb{Z}^{k}} \check{f}\left(n_{1}, \ldots, n_{k}\right)$$
where
$$\check{f}\left(n_{1}, \ldots, n_{k}\right)=\int_{\mathbb{R}^{k}} f(x_1, \cdots, x_k) e^{\sum_{j=1}^{k}\limits 2 \pi \sqrt{-1} n_j x_j}\D x_{1} \cdots\D x_{k}$$

 As a consequence, we have
\begin{pr}\label{pr2}
 $$\mathrm{RT}_r({\mathcal{K}_{p'}(p,q)})=c_r\sum_{s=0}^{|q|-1}\sum_{(m,n,l) \in \mathbb{Z}^3}\check{f}_r(s,m,n,l)+o(e^{\frac{r}{4\pi}\cdot3.5})$$
where
\begin{eqnarray}
  \check{f}_r(s,m,n,l)&=&4(-1)^{m+n+l}(\frac{r}{2\pi})^3\int_{D_{\frac{\epsilon}{2}}}\psi(x,y,z)\sin(\frac{x}{q}-J(s)\pi)\sin (2z)\label{ef}\\
  &&\cdot e^{-\sqrt{-1}(x+y)+\frac{r}{4 \pi \sqrt{-1}}V_r(x,y,z,s,m,n,l)}\D x\D y\D z\nonumber
\end{eqnarray}
\end{pr}

\begin{proof}[Proof of Propositon \ref{pr2}]
    We abuse the notation $f_r(s,m,n,l')=f_r(s,x,y,z)$, which does not cause ambiguity. We have
  \begin{eqnarray*}
    &&\sum_{(m,n,l')\in (\mathbb{Z}+\frac{1}{2})^3}f_r(s,m,n,l')\\
   &=&\sum_{(m,n,l)\in \mathbb{Z}^3}f_r(s,m-\frac{1}{2},n-\frac{1}{2},l-\frac{1}{2})\\
   &=&\sum_{(m,n,l)\in \mathbb{Z}^3}\int_{\mathbb{R}^{3}} f_r(s,m'-\frac{1}{2},n'-\frac{1}{2},l'-\frac{1}{2}) e^{2 \pi \sqrt{-1}(mm'+nn'+ll')}\D m'\D n'\D l'\\
   &=&\sum_{(m,n,l)\in \mathbb{Z}^3}\int_{D_{\frac{\epsilon}{2}}} f_r(s,m',n',l') e^{2 \pi \sqrt{-1}(m(m'+\frac{1}{2})+n(n'+\frac{1}{2})+l(l'+\frac{1}{2}))}\D m'\D n'\D l'\\
   &=&(\frac{r}{2\pi})^3\sum_{(m,n,l)\in \mathbb{Z}^3}(-1)^{m+n+l}\int_{D_{\frac{\epsilon}{2}}} f_r(s,x',y',z') e^{2 \pi \sqrt{-1}(mm'+nn'+ll')}\D x'\D y'\D z'\\
   &=&(\frac{r}{2\pi})^3\sum_{(m,n,l)\in \mathbb{Z}^3}(-1)^{m+n+l}\int_{D_{\frac{\epsilon}{2}}} f_r(s,x',y',z') e^{\frac{r}{4\pi\sqrt{-1}}(-4\pi)(mx'+ny'+lz')}\D x'\D y'\D z'
  \end{eqnarray*}

  plug the previous formulas into (\ref{e9}), we get Propositon \ref{pr2}.
\end{proof}

\begin{pr}\label{pr2.1}
    \begin{eqnarray}
        \check{f}_r(s,m,n,l)&=&4\sqrt{\frac{2}{r}}e^{-\frac{3\pi\sqrt{-1}}{4}}(-1)^{m+n+l}(\frac{r}{2\pi})^3\int_{D_{\frac{\epsilon}{2}}}g(s,x,y,z)e^{\frac{r}{4 \pi \sqrt{-1}}V(x,y,z,s,m,n,l)}\D x\D y\D z\nonumber\\
        &&\cdot(1+O(\frac{1}{r}))\nonumber
    \end{eqnarray}
    $$g(s,x,y,z)=\frac{\psi(x,y,z)\sin(\frac{x}{q}-J(s)\pi)\sin (2z)}{\sqrt{(1-e^{2\sqrt{-1}(y+x)})(1-e^{2\sqrt{-1}(y-x)})(1-e^{-2\sqrt{-1}y})}}$$
    and the amplitude of $\sqrt{\cdots}$ is between $-\frac{\pi}{2}$ and $\frac{\pi}{2}$.
\end{pr}

\begin{proof}[Proof of Proposition \ref{pr2.1}]
  By Lemma \ref{convergeV}, we have
  \begin{eqnarray*}
    &&e^{-\sqrt{-1}(x+y)+\frac{r}{4 \pi \sqrt{-1}}V_r(x,y,z,s,m,n,l)}\\
    &=&e^{\frac{r}{4 \pi \sqrt{-1}}V(x,y,z,s,m,n,l)}e^{-\frac{1}{2}(\log(1-e^{2\sqrt{-1}(y+x)})+\log(1-e^{2\sqrt{-1}(y-x)})+\log(1-e^{-2\sqrt{-1}y})+\log r+\frac{3\pi\sqrt{-1}}{2}-\log2))}+o(e^{\frac{r}{4\pi}\cdot3.5})\\
    &=&\sqrt{\frac{2}{r}}e^{-\frac{3\pi\sqrt{-1}}{4}}\frac{1}{\sqrt{(1-e^{2\sqrt{-1}(y+x)})(1-e^{2\sqrt{-1}(y-x)})(1-e^{-2\sqrt{-1}y})}}e^{\frac{r}{4 \pi \sqrt{-1}}V(x,y,z,s,m,n,l)}+o(e^{\frac{r}{4\pi}\cdot3.5})
  \end{eqnarray*}
\end{proof}

\subsection{Symmetry of the Fourier coefficients}\label{subsec4.3}
We will discuss the symmetry of the Fourier coefficients in this section, which generalizes the equality of the leading Fourier terms and the \textbf{Big Cancellation} discovered by Chen-Zhu \cite{Chen-Zhu-1}.

\begin{pr}\label{fx=f-x}
    For $(s,m)$ and $(s',m')$ such that $k(s,m)+k(s',m')=0$, we have
    $$\check{f}_r(s,m,n,l)=\check{f}_r(s',m',n,l)$$
    Particularly, we have
    $$\check{f}_r(s^+,m^+,n,l)=\check{f}_r(s^-,m^-,n,l)$$
\end{pr}

\begin{proof}[Proof of Proposition \ref{fx=f-x}]
    \begin{eqnarray*}
      &&\check{f}_r(s',m',n,l)\\
      &=&4(-1)^{m'+n+l}(\frac{r}{2\pi})^3\int_{D_{\frac{\epsilon}{2}}}\psi(x,y,z)\sin(\frac{x}{q}-J(s')\pi)\sin (2z)\cdot e^{-\sqrt{-1}(x+y)+\frac{r}{4 \pi \sqrt{-1}}V_r(x,y,z,s',m',n,l)}\D x\D y\D z\\
      &=&4(-1)^{m'+n+l}(\frac{r}{2\pi})^3\int_{D_{\frac{\epsilon}{2}}}\psi(-x,y,z)\sin(\frac{-x}{q}-J(s')\pi)\sin (2z)\cdot e^{-\sqrt{-1}(-x+y)+\frac{r}{4 \pi \sqrt{-1}}V_r(-x,y,z,s',m',n,l)}\D x\D y\D z\\
      &=&4(-1)^{m'+n+l+1}(\frac{r}{2\pi})^3\int_{D_{\frac{\epsilon}{2}}}\psi(x,y,z)\sin(\frac{x}{q}-J(s)\pi)\sin (2z)\cdot e^{-\sqrt{-1}(-x+y)+\frac{r}{4 \pi \sqrt{-1}}V_r(-x,y,z,s',m',n,l)}\D x\D y\D z\\
      &\stackrel{(\ref{e20})}{=}&4(-1)^{m'+n+l+1}(\frac{r}{2\pi})^3\int_{D_{\frac{\epsilon}{2}}}\psi(x,y,z)\sin(\frac{x}{q}-J(s)\pi)\sin (2z)\\
      &&\cdot e^{-\sqrt{-1}(-x+y)+\frac{r}{4 \pi \sqrt{-1}}(V_r(x,y,z,s,m,n,l)+\frac{8\pi}{r}x-(K(s)-K(s'))\pi^2)}\D x\D y\D z\\
      &=&4(-1)^{m'+n+l+1}(\frac{r}{2\pi})^3\int_{D_{\frac{\epsilon}{2}}}\psi(x,y,z)\sin(\frac{x}{q}-J(s)\pi)\sin (2z)\\
      &&\cdot e^{-\sqrt{-1}(x+y)+\frac{r}{4 \pi \sqrt{-1}}(V_r(x,y,z,s,m,n,l)-(K(s)-K(s'))\pi^2)}\D x\D y\D z\\
      &=&(-1)^{1+m'-m+\frac{K(s)-K(s')}{4}r}\check{f}_r(s,m,n,l)\\
      &\stackrel{Lemma\ \ref{IJK}}{=}&(-1)^{1+m'-m+\frac{K(s)-K(s')}{4}}\check{f}_r(s,m,n,l)\\
      &\stackrel{Lemma\ \ref{IJK}}{=}&\check{f}_r(s,m,n,l)
    \end{eqnarray*}

    The third equality is because we have $\sin(\frac{-x}{q}-J(s')\pi)=-\sin(\frac{x}{q}-J(s)\pi)$ by Lemma \ref{IJK2}.

\end{proof}

Similarly by (\ref{ef}) and (\ref{e21}) we have the following \textbf{Big Cancellation} as showed in \cite{Chen-Zhu-1}.
\begin{pr}\label{cancell}
    \begin{equation}\label{e22}
      \check{f}_r(s,m,n,l)=(-1)^{p'+l}\check{f}_r(s,m,n,-2p'-2-l)
    \end{equation}
\end{pr}
For example, we have
$$\check{f}_r(s,m,n,-p'-1)=0$$
$$\check{f}_r(s,m,n,-p')=\check{f}_r(s,m,n,-p'-2)$$

\begin{rem}\label{rem1}
    Actually, $V(x_0,y_0,z_0,s^+,m^+,0,-p'-1)$ and $V(x_0,y_0,z_0,s^-,m^-,0,-p'-1)$\footnote{$(x_0,y_0,z_0)$ is the critical point of $V$} are equal to $\vol(\mathcal{W}((p+4q,-q),(0,0)))=\vol(\mathcal{W}((p,q),(0,0)))$, which is bigger than $\vol(\mathcal{K}_{p'}(p,q))=\vol(\mathcal{W}((p,q),(1,-p')))$. This is an insurmountable problem in the estimate part (see Proposition \ref{pr<Vol}). Fortunately, the \textbf{Big Cancellation} tell us $\check{f}_r(s,m,n,-p'-1)=0$ and we overcome this difficulty.
\end{rem}

We'll talk more about symmetry of the Fourier coefficients in Section \ref{sec7.3} and Section \ref{sec7.5}.

\section{Geometry of the critical points}\label{sec5}
The goal of this section is to understand the geometric meaning of the critical points and the critical values of the functions $V$ defined in the previous section. The main result is Proposition \ref{Vol} in Section \ref{sec5.4} which shows that the critical values of $V$ equal to the $\sqrt{-1}\cv(M)$ by Neumann-Zagier-Yoshida's Theory \cite{NZ}\cite{Yoshida}.

We construct a ideal triangulation of $\mathcal{W}$ in Section \ref{sec5.1}, which we want it's hyperbolicity equations corresponding to the critical point equations of $V$ directly. However, this part of the study was not as simple as we expected. We discover that the system of critical point equations of $V$ corresponding to the hyperbolicity equations of $\mathcal{W}((p+4q,-q),(1,p'-\frac{1}{2}))$ instead of $\mathcal{W}((p,q),(1,-p'))$ in Section \ref{secV}.

We find that this intriguing phenomenon is closely related to Hodgson-Meyerhoff-Weeks' conclusion \cite{HMW} (Figure \ref{sister manifolds}). To overcome this difficulty, we proposed a theory called the sister potential function in Section \ref{sec5.3}.


\subsection{Geometry of Whitehead link complement}\label{sec5.1}
First, let's discuss the hyperbolicity of $\mathcal{W}$.

By \cite{MP}, we have
\begin{lem}\label{nonhyperbolic1}
    $\mathcal{W}((p,q),(m,l))$ is a closed hyperbolic manifold unless one of the following holds:
    \begin{enumerate}[(1)]
      \item $\frac{p}{q}$ or $\frac{m}{l}\in\{0,1,2,3,4,\infty\}$
      \item $(\frac{p}{q},\frac{m}{l})$ or $(\frac{m}{l},\frac{p}{q})\in$
      $$\{(-4,-1),(-3,-1),(-2,-2),(-2,-1),(\frac{3}{2},5),(\frac{4}{3},5),(\frac{5}{2},\frac{7}{2})\}$$
    \end{enumerate}
\end{lem}

Particularly, we have

\begin{lem}\label{nonhyperbolic2}
    $\mathcal{K}_{p'}(p,q)$ is a closed hyperbolic manifold unless one of the following holds:
    \begin{enumerate}[(1)]
      \item $p'\in\{0,-1\}$
      \item $\frac{p}{q}\in\{0,1,2,3,4,\infty\}$
      \item $p'=1,\frac{p}{q}\in\{-2,-3,-4\}$
    \end{enumerate}
\end{lem}

\begin{figure}[H]
      \centering
        \begin{minipage}[b]{0.36\linewidth}
          \includegraphics[width=1\linewidth]{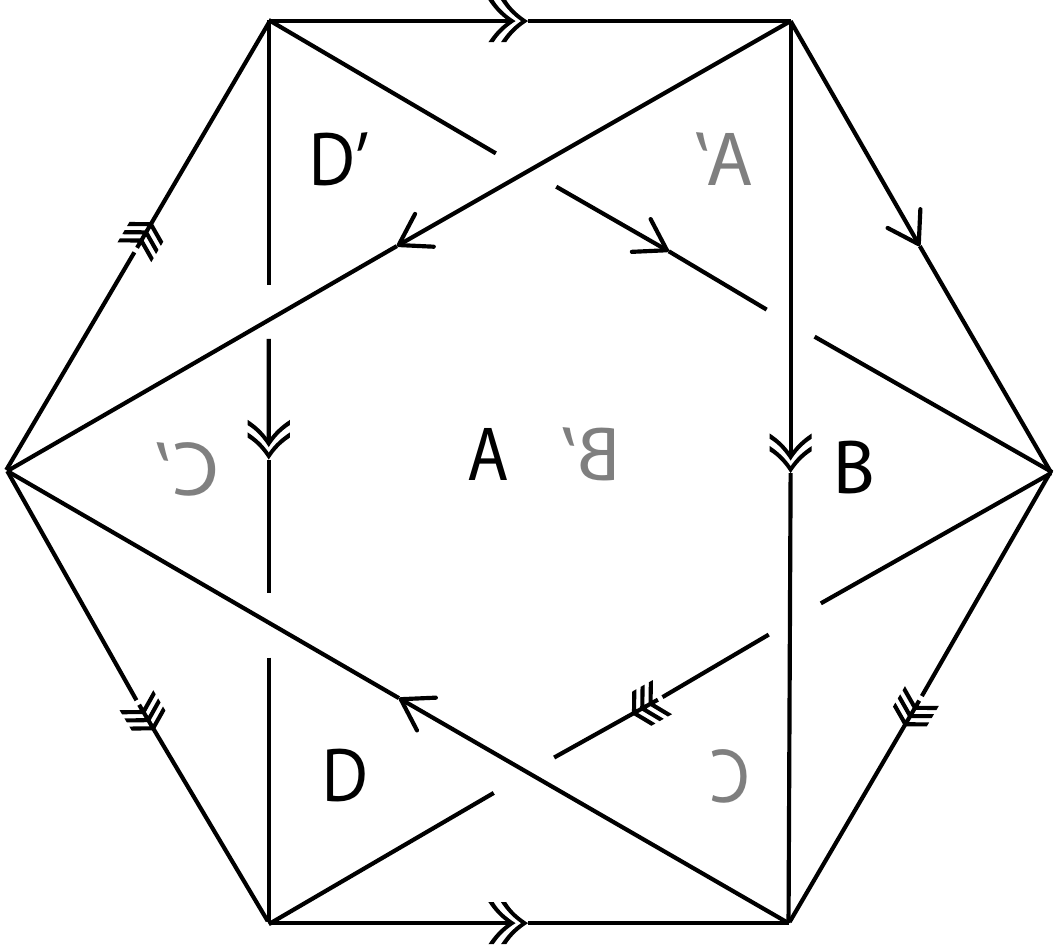}
            \caption{An ideal octahedron}\label{octahedron}
        \end{minipage}
        \begin{minipage}[b]{0.4\linewidth}
          \includegraphics[width=1\linewidth]{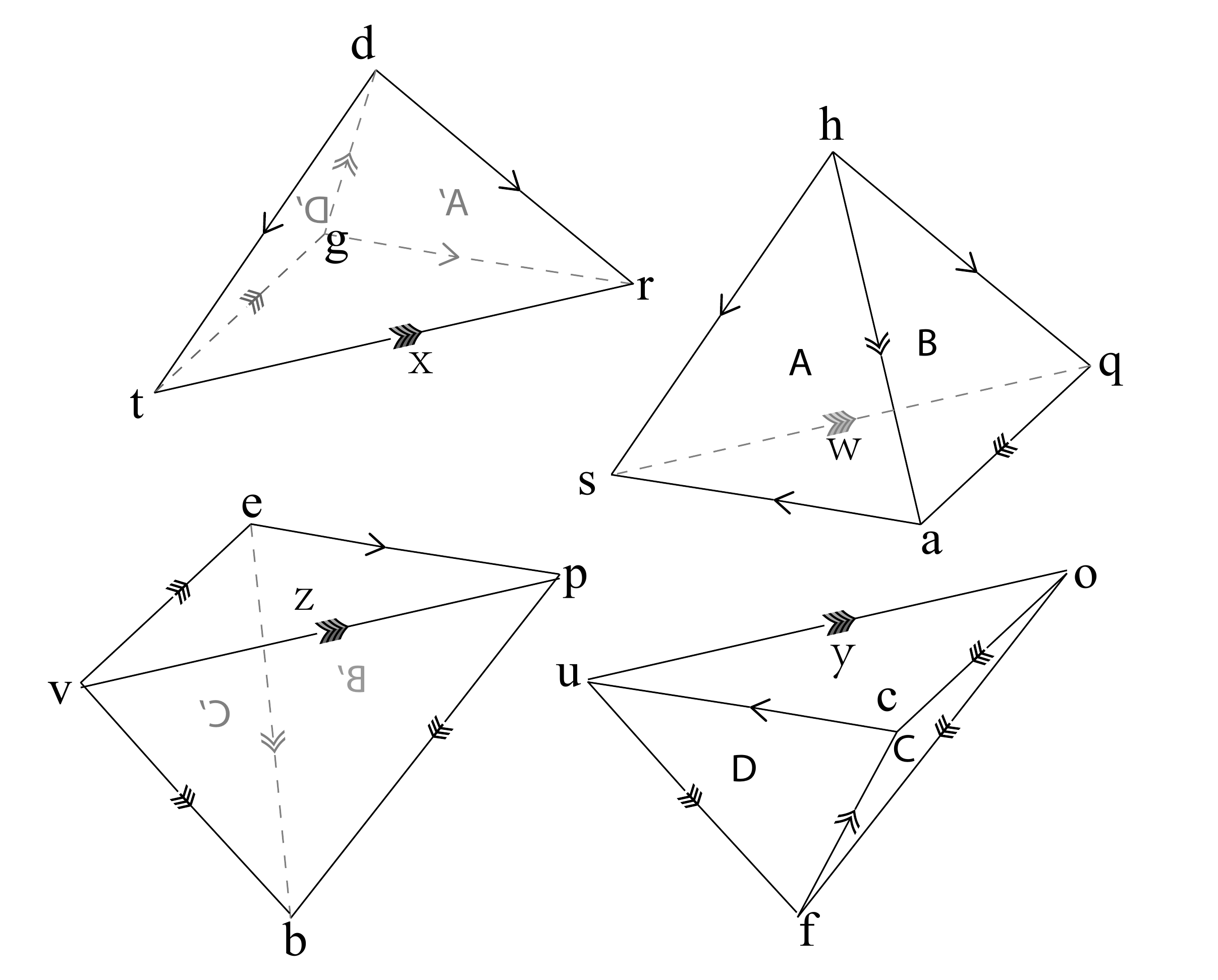}
          \caption{$4$ ideal tetrahedrons}\label{4tetrahedron}
        \end{minipage}
    \end{figure}

According to Thurston's notes \cite{Thurston} and \cite{NR}, the Whitehead link complement has an ideal triangulation(see Figure \ref{4tetrahedron}) of $4$ ideal tetrahedrons comes from cutting an ideal octahedron along the horizontal diagonal line in Figure \ref{octahedron}.

\begin{figure}[H]
      \centering
        \includegraphics[width=0.4\linewidth]{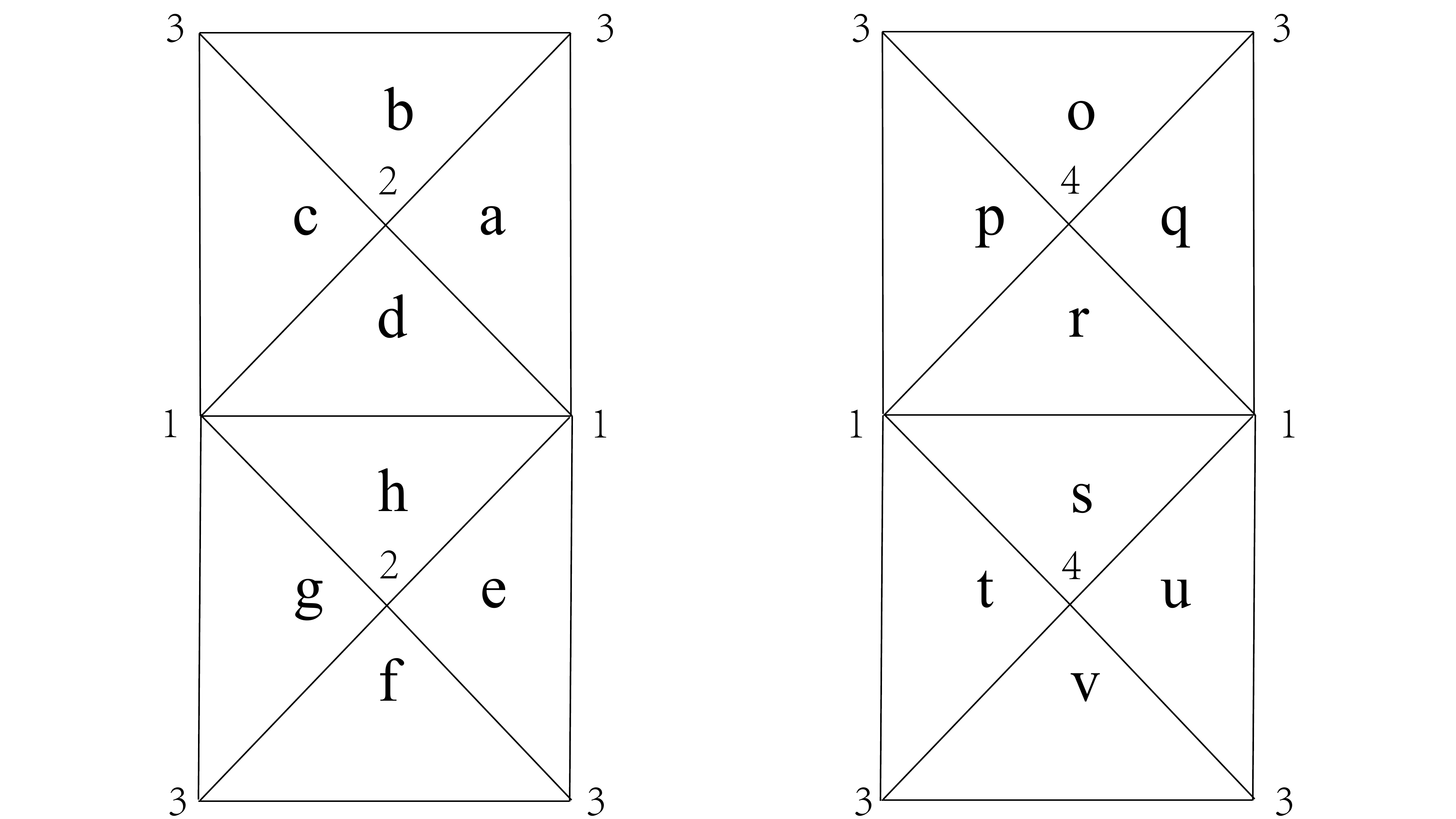}
        \caption{Fundamental domains of the two boundaries of the tubular neighborhood of the Whitehead link complement.}\label{cusp1}
    \end{figure}

Figure \ref{cusp1} is the fundamental domains of the two boundarys of the tubular neighborhood of the Whitehead link complement.

\begin{figure}[H]
      \centering
        \includegraphics[width=0.7\linewidth]{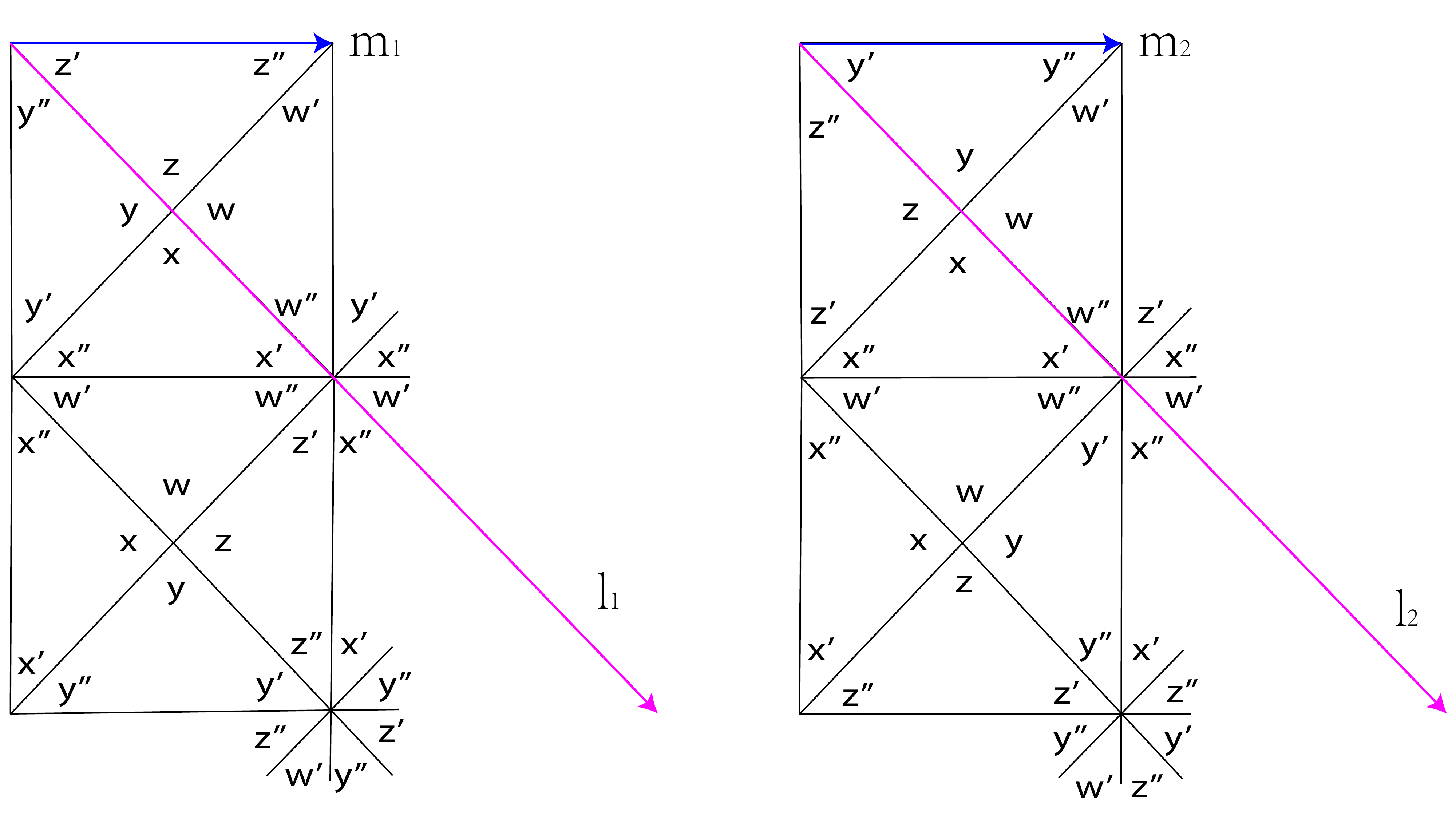}
        \caption{meridians and longitudes of $W$ in the fundamental domains}\label{cusp2}
    \end{figure}

Figure \ref{cusp2} shows the meridians and longitudes of $W$ in the fundamental domains.

The gluing equations of this triangulation are
\begin{eqnarray}
  \log w+\log x+\log y+\log z&=&2\pi \sqrt{-1} \label{e12}\\
  \log(1-w)+\log(1-x)-\log(1-y)-\log(1-z) &=&0 \label{e13}
\end{eqnarray}
where (\ref{e12}) represent the glue equations of the edges marked with $2$ arrows or $4$ arrows, (\ref{e13}) represent the glue equations of the edges marked with $1$ arrow or $3$ arrows.

The holonomy of meridians and longitudes of $W$ can be represented as
\begin{eqnarray}\label{e14}
  m_1&=&\log(w-1)+\log x+\log y-\log(y-1)-\pi\sqrt{-1}\\
  l_1 &=&2\log x+2\log y-2\pi\sqrt{-1}\nonumber\\
  m_2&=&\log(w-1)+\log x+\log z-\log(z-1)-\pi\sqrt{-1}\nonumber\\
  l_2 &=&2\log x+2\log z-2\pi\sqrt{-1}\nonumber
\end{eqnarray}

The complete hyperbolic structure of the Whitehead link complement corresponding to
$$m_1=l_1=m_2=l_2=0\Leftrightarrow w=x=y=z=\sqrt{-1}$$

Our next goal is to transform the ideal triangulation above to a new  triangulation corresponding to $V$. The new triangulation has five tetrahedrons, which corresponding to the five $\li(\cdots)$ in $V$.

\begin{rem}
    We will show some figures about the ideal triangulation of $\mathcal{W}$. We will put a letter on each face, such as $A,A',B,b,c,c',C,d$. $'$ means glue together. For example, we glue face $A$ and face $A'$ together. Please note that the letters are case sensitive.

    The letters on the back are displayed in reverse. Notice that the ``$b$" on the back looks like a "$d$" when written in reverse. Please don't be confused in this case.
\end{rem}

The new construction requires the following four steps.
\begin{enumerate}[(1)]
  \item The above  triangulation of Neumann is constructed by cutting the octahedron into four tetrahedrons(see Figure \ref{octahedron} and Figure \ref{4tetrahedron}). Our first step is to take a different approach to cut the octahedron into four tetrahedrons (see Figure \ref{4tetrahedron2}) along the diagonal line from bottom left to top right in Figure \ref{octahedron}.
          \begin{figure}[H]
            \centering
            \includegraphics[width=0.7\linewidth]{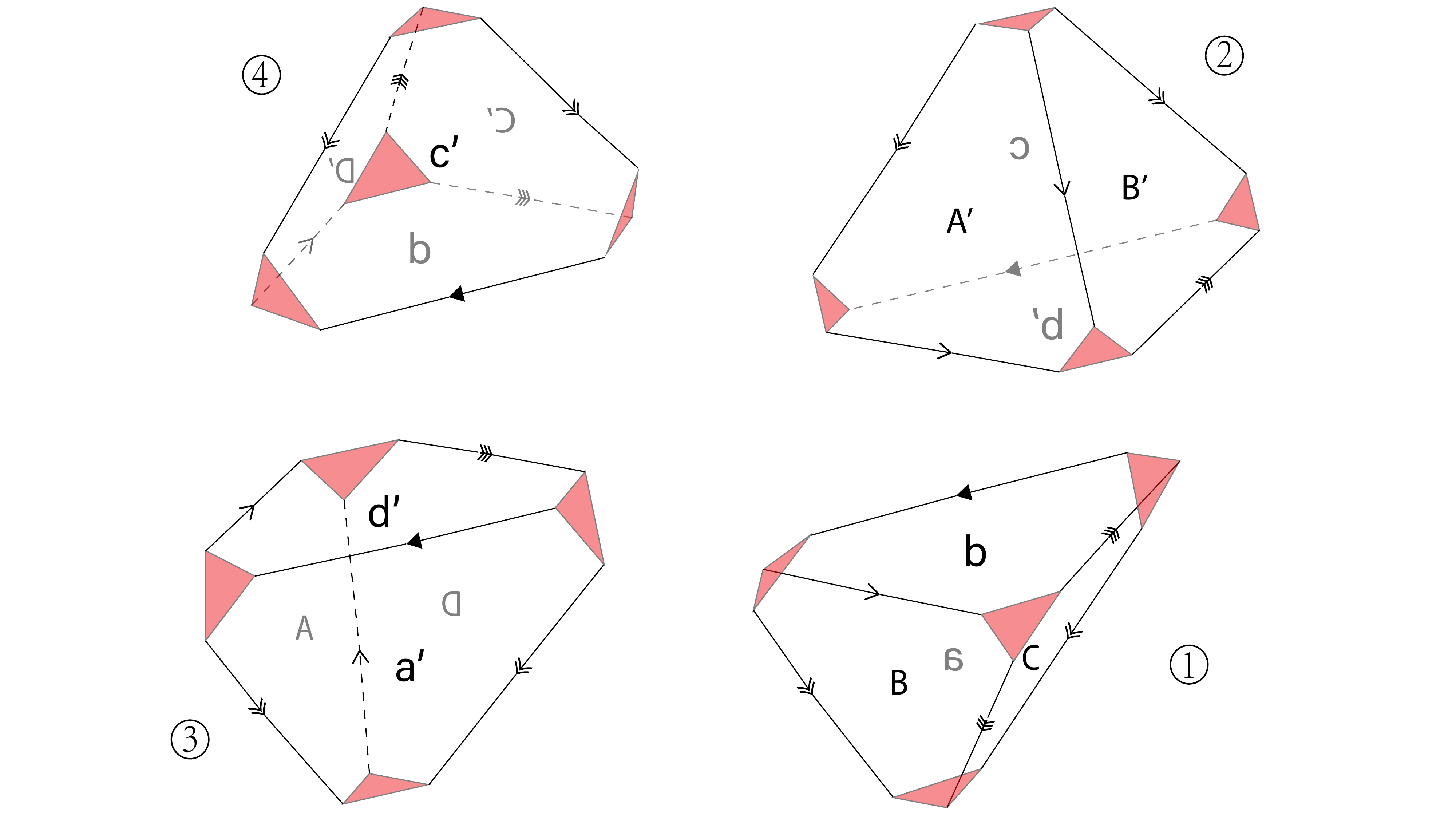}
            \caption{A different approach to cut  the octahedron into four tetrahedrons}\label{4tetrahedron2}
         \end{figure}

  \item We recombine these four tetrahedrons as Figure \ref{4tetrahedron3} and Figure \ref{platform}.
        \begin{figure}[H]
            \centering
            \includegraphics[width=0.7\linewidth]{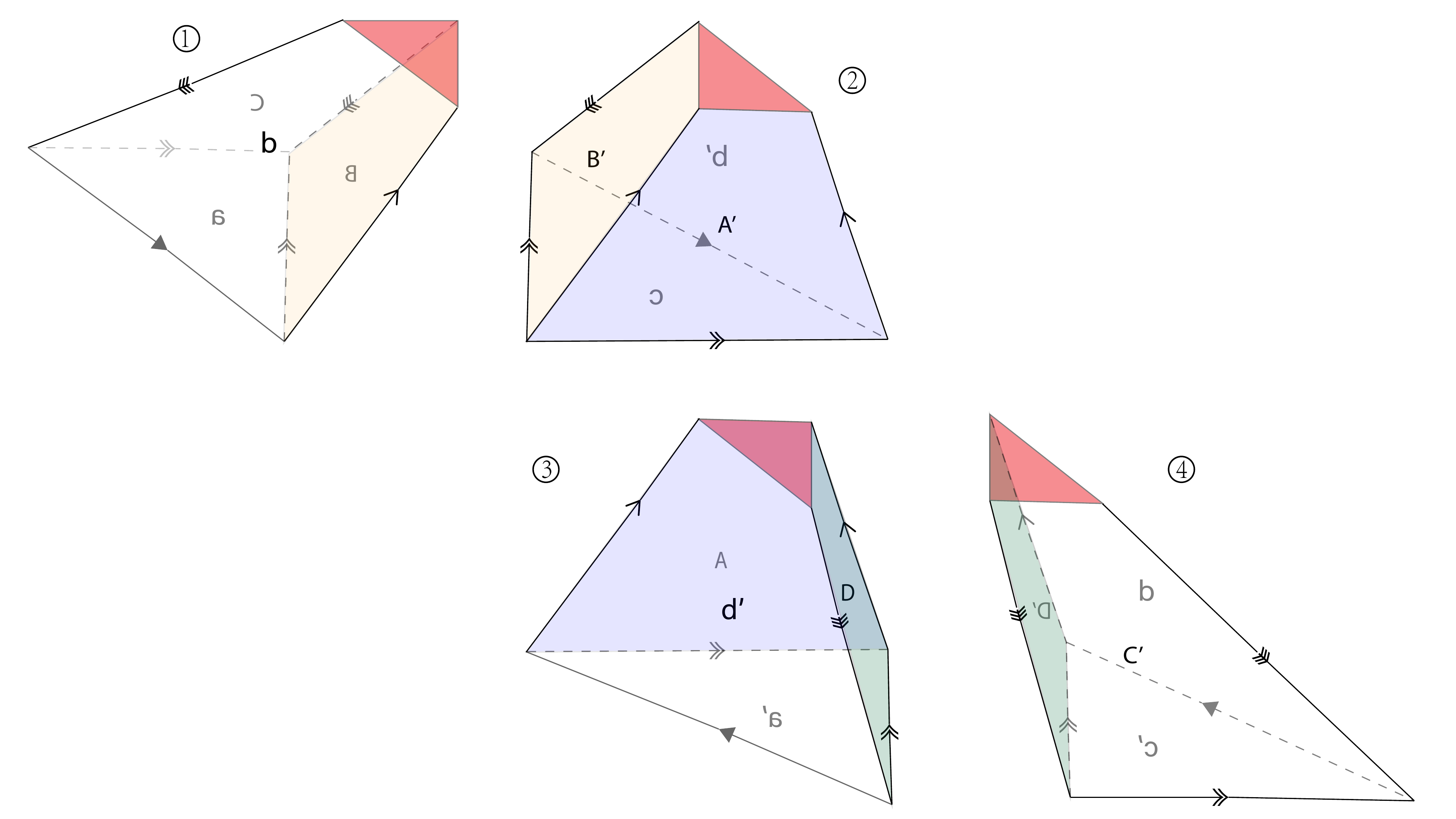}
            \caption{Recombine these four tetrahedrons}\label{4tetrahedron3}
         \end{figure}

         \begin{figure}[H]
            \centering
            \includegraphics[width=0.54\linewidth]{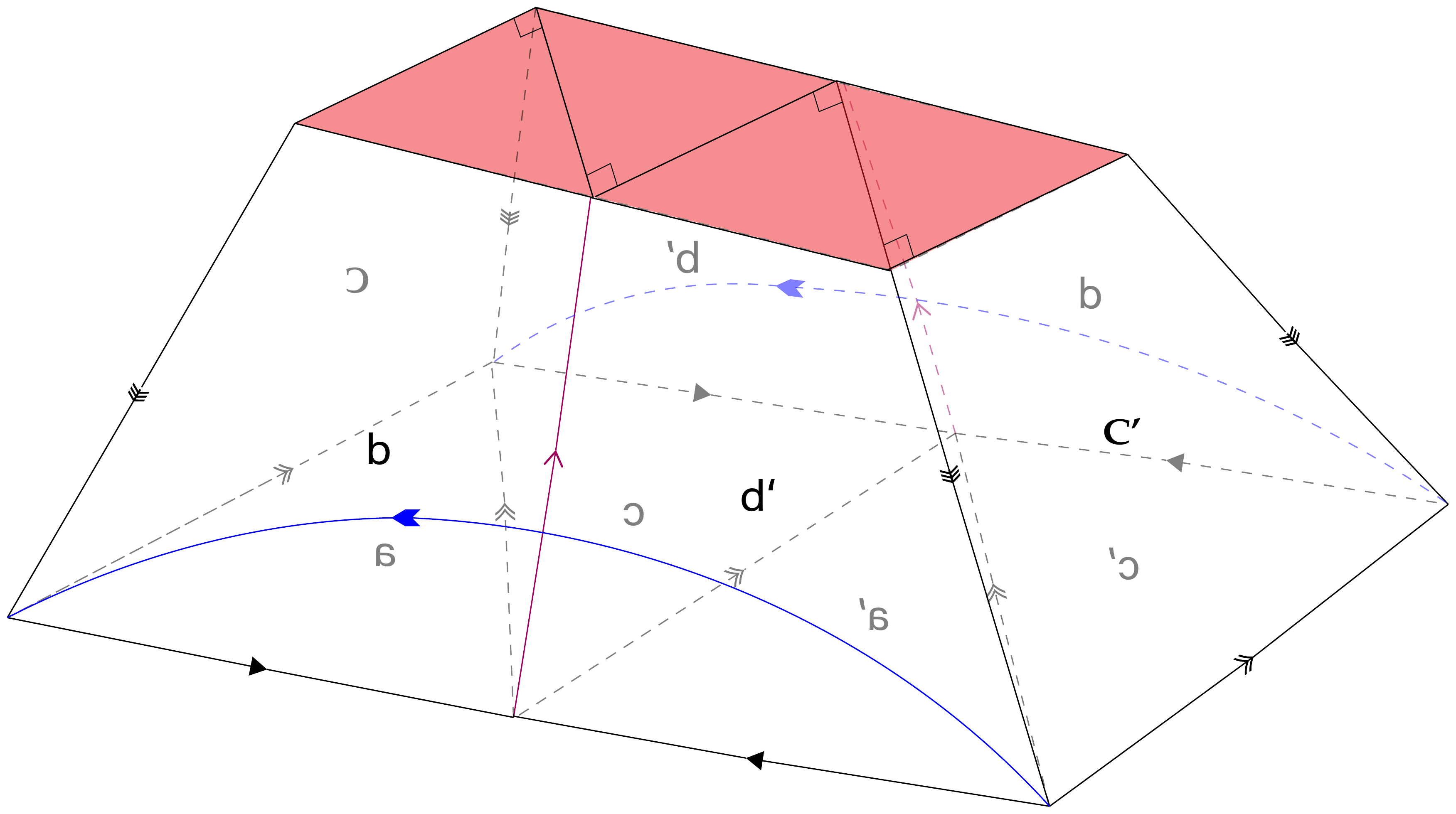}
            \caption{Glue the four tetrahedrons in a new way}\label{platform}
         \end{figure}

  \item We remove the purple edge with an arrow\footnote{This edge are $E\infty$ and $F\infty$ in Figure \ref{quadrangular pyramid}.} (see Figure \ref{platform}) and add the edge with blue label ($BC$ and $DA$ in Figure \ref{quadrangular pyramid}).
        \begin{figure}[H]
            \centering
            \includegraphics[width=0.6\linewidth]{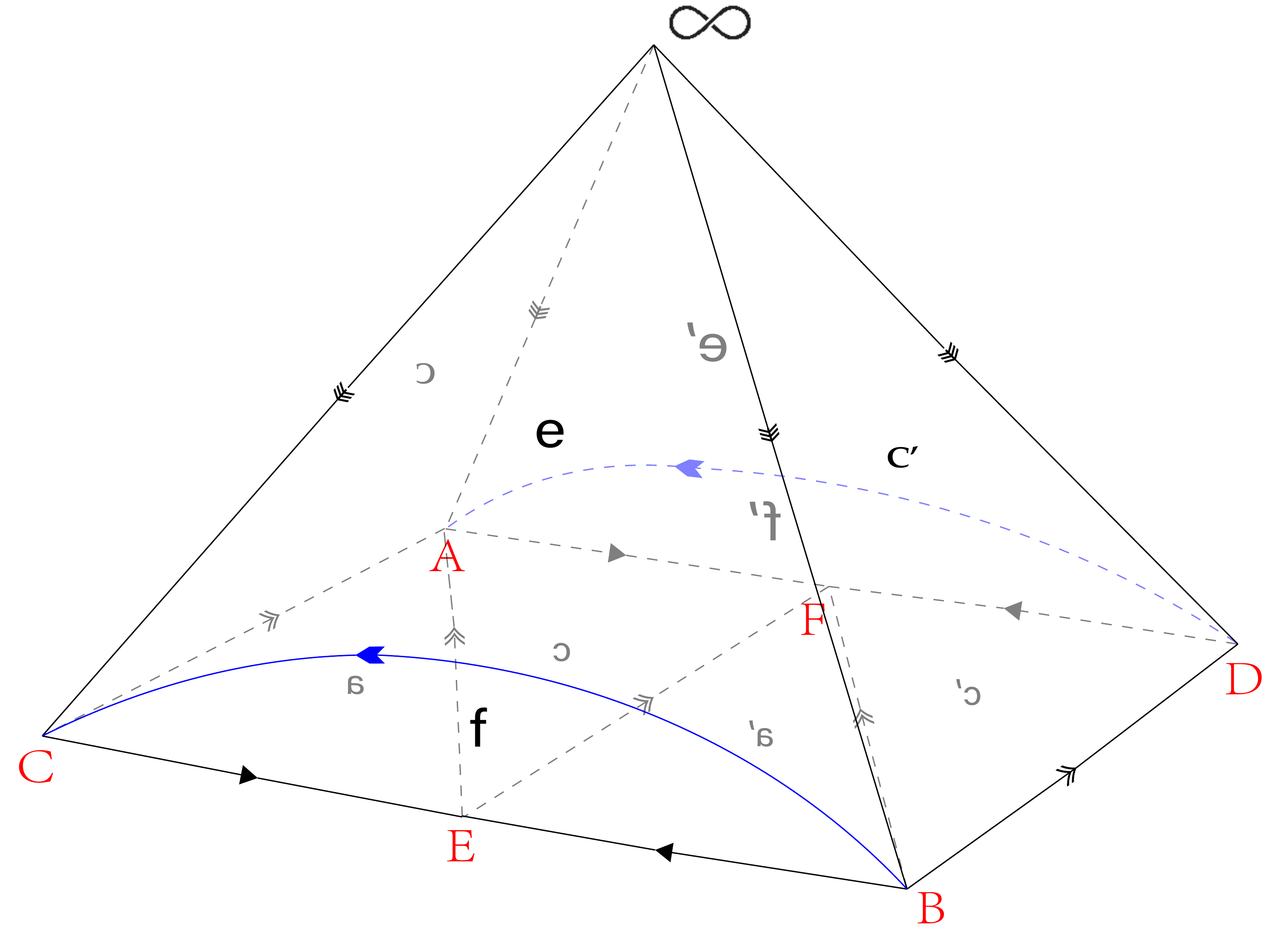}
            \caption{Add an edge marked with blue label}\label{quadrangular pyramid}
         \end{figure}

  \item We think of the quadrangular pyramid in Figure \ref{quadrangular pyramid} as the pentagonal bipyramid in Figure \ref{pentagonal bipyramid}.

      The link $L(\infty)$ (The upper surface of the "platform" in Figure \ref{platform}) corresponding to the fundamental domains of a boundary of the tubular neighborhood of the Whitehead link complement. $L(A)$, $L(B)$, $L(C)$, $L(D)$, $L(E)$, $L(F)$ are put together to form the fundamental domains of another boundary.

       We can divide this pentagonal bipyramid into five tetrahedrons $AB\infty C$, $ABCE$, $ABEF$, $ABFD$, $ABD\infty$ with shape $c_1$, $c_2$, $c_3$, $c_4$, $c_5$.
          \begin{figure}[H]
            \centering
            \includegraphics[width=0.5\linewidth]{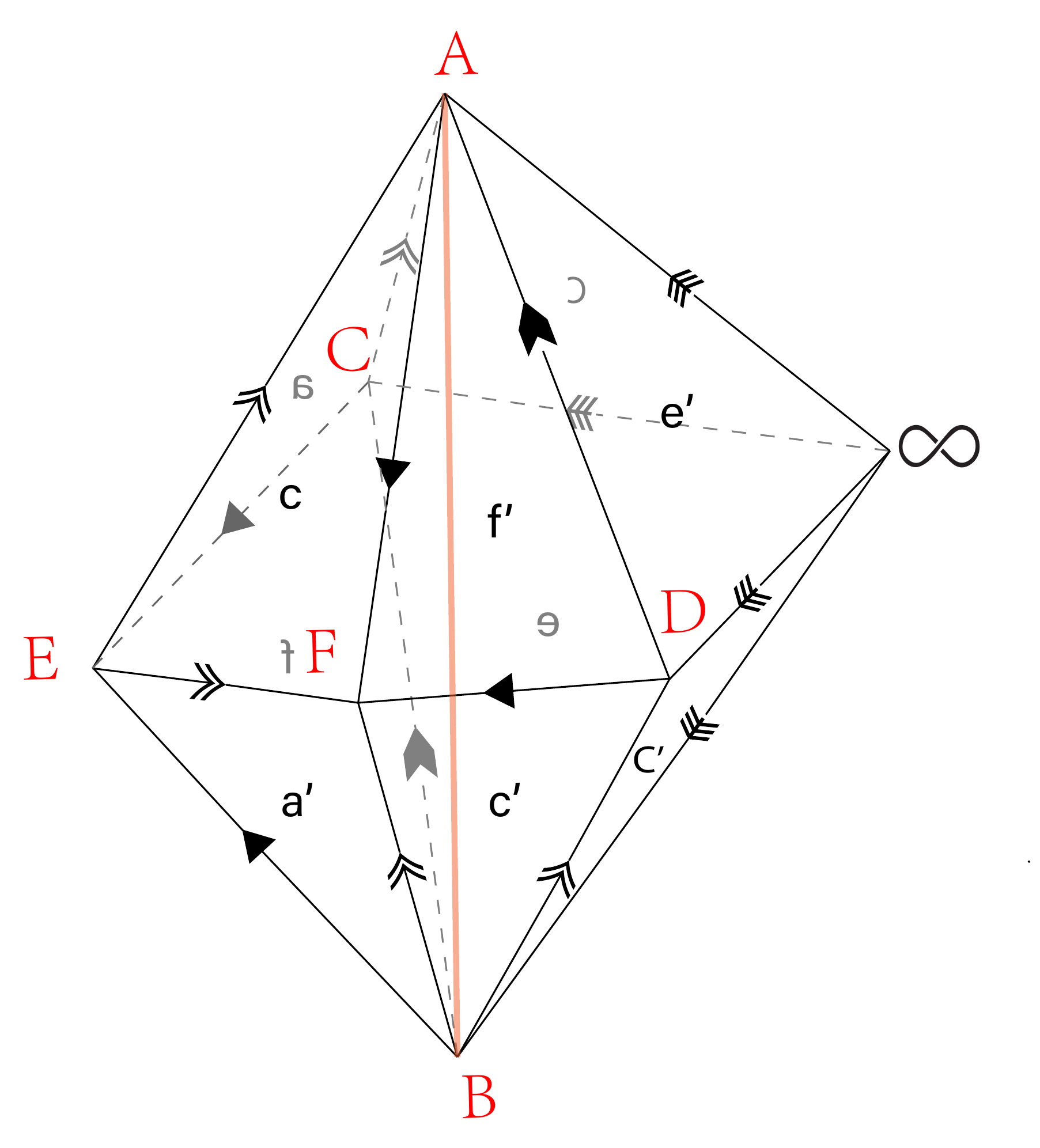}
            \caption{An ideal pentagonal bipyramid}\label{pentagonal bipyramid}
         \end{figure}
\end{enumerate}

Our next goal is to find the gluing equations (\ref{e15}) and the Dehn-surgery equations (\ref{e16}) for this triangulation. We can get the  Fundamental domains of the two boundaries from Figure \ref{pentagonal bipyramid}, and then get the holonomy of meridians and longitudes. This method involves too much stereogeometric imagination, we can use another simple method to achieve this purpose. In fact, our new triangulation comes from Thurston's triangulation by a Pachner $2$-$3$ move and a Pachner $4$-$4$ move, so we can get the relationship between $c_i$ and $x,y,z,w$ (\ref{ecx}) (\ref{exc}) by using the $5$-term formula (Proposition \ref{5terms1}) for three times. Then plug (\ref{ecx}) (\ref{exc}) into (\ref{e12}) (\ref{e13}) (\ref{e14}), we achieved our goal.

The calculation of five-term formulas is tedious and boring, we have the following trick to get the relationship between $c_i$ and $x,y,z,w$ directly without calculate three times. We can use $w,x,y,z$ to represent the shape of $L(\infty)$, so we can get the coordinate (in $\mathbb{C}$) of $A$, $B$, $C$, $D$, $E$, $F$ as the figure below.
    \begin{figure}[H]
        \centering
        \includegraphics[width=0.9\linewidth]{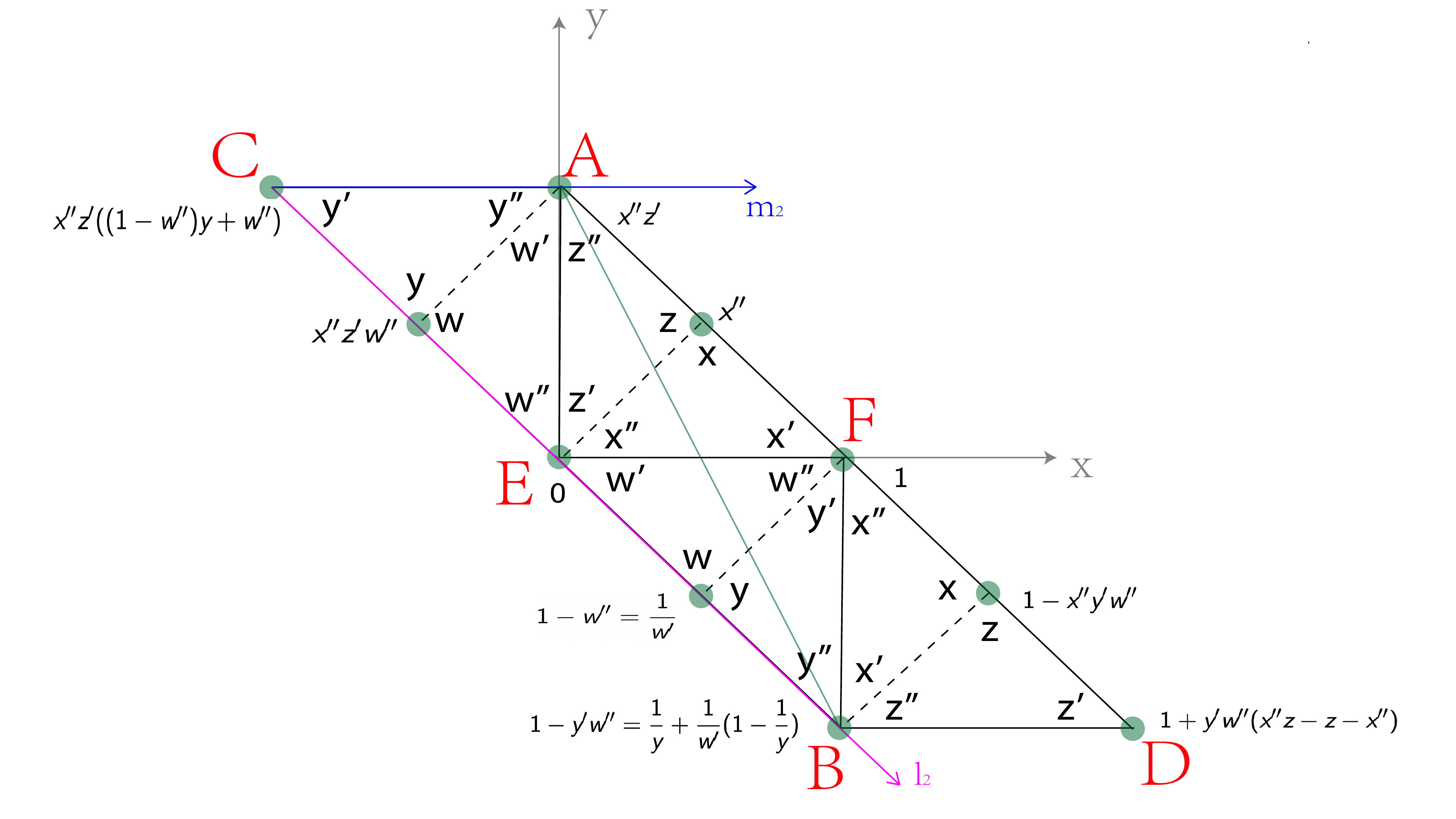}
        \caption{The shape of $L(\infty)$}\label{cusp3}
    \end{figure}

Suppose $E=0$ and $F=1$ in the complex plane, then we have
$$A=x''z',C=x''z'((1-w'')y+w''),B=1-y'w'',D=1+y'w''(x''z-z-x'')$$
as showed in Figure \ref{cusp3}. We have
\begin{eqnarray}
  c_1&=&[A,B;\infty,C]=\frac{wy(y-1)(z-1)}{(wy-1)(yz+xz-y-z)}\label{ecx}\\
  c_2&=&[A,B;C,E]=\frac{(w-1)(yz+xz-y-z)}{wz(x-1)(y-1)}\nonumber\\
  c_3&=&[A,B;E,F]=\frac{(wy-1)(xz-1)}{(y-1)(z-1)}\nonumber\\
  c_3&=&[A,B;F,D]=\frac{(x-1)(wy+yz-y-z)}{xy(z-1)(w-1)}\nonumber\\
  c_5&=&[A,B;D,\infty]=\frac{xz(y-1)(z-1)}{(xz-1)(wy+yz-y-z)}\nonumber
\end{eqnarray}
so we have
$$c_1c_2c_3c_4c_5=1$$
which corresponding to the edge equation of $AB$ in Figure \ref{pentagonal bipyramid} (five tetrahedrons glue together around $AB$). We can represent $w,x,y,z$ by $c_1,c_2,c_3,c_4,c_5$ (this step can be calculated with the help of Mathmatica).
\begin{eqnarray}
  x&=&\frac{1-c_1c_2}{c_1c_2(c_4-1)}\label{exc}\\
  y&=&-\frac{c_1(c_2-1)(c_3-1)}{(c_1-1)(\frac{1}{c_4c_5}-1)}\nonumber\\
  z&=&-\frac{c_1c_2(c_3-1)(c_4-1)}{(c_1c_2-1)(\frac{1}{c_5}-1)}\nonumber\\
  w&=&\frac{\frac{1}{c_4c_5}-1}{c_2-1}\nonumber
\end{eqnarray}
By plugging the above formula into (\ref{e12}) (\ref{e13}) (\ref{e14}) and simplifying it(use $c_1',c_2'',c_3'',c_4'',c_5'$ instead of $c_1,c_2,c_3,c_4,c_5$), we get the gluing equation of our new ideal triangulation is
\begin{eqnarray}\label{e15}
  \sum_{i=1}^5\log c_i&=&2\pi\sqrt{-1}\\
  2\log c_3''+\log c_1'+\log c_5'&=&2\pi\sqrt{-1}\nonumber\\
  2\log c_3''-\log c_2''-\log c_4''&=&0\nonumber
\end{eqnarray}
and holonomies of the meridians and longitudes of $W$ can be represented as
\begin{eqnarray}\label{e16}
  m_1&=&\log\frac{c_3''}{c_2''}\\
  l_1 &=&2\log\frac{c_3''}{c_2''}+2\log(1-\frac{{c_3''}^2}{c_2''})-2\log(1-c_2'')\nonumber\\
  m_2&=&2\log(-c_1'c_3'')+\log(1-\frac{1}{c_1'{c_3''}^2})-\log(1-c_1')\nonumber\\
  l_2&=&2\log(-c_1'c_3'')\nonumber
\end{eqnarray}

the Dehn-surgery equations of $\mathcal{K}_{p'}(p,q)=\mathcal{W}((p,q),(1,-p'))$ are
\begin{eqnarray}
  pm_1+ql_1&=&2\pi\sqrt{-1}\label{eD}\\
  m_2-p'l_2&=&2\pi\sqrt{-1}\nonumber
\end{eqnarray}

By eliminating $c_4,c_5$ from (\ref{e15}), we get the gluing equation (\ref{e15}) is equivalent to
\begin{equation}\label{e17}
  2\log c_3''-\log(1-\frac{{c_3''}^2}{c_2''})-\log(1-c_2'')-\log(1-\frac{1}{c_3''})+\log(1-c_1')+\log(1-\frac{1}{c_1'{c_3''}^2})=0
\end{equation}

The complete hyperbolic structure of the Whitehead link complement corresponding to
$$c_1'=c_5'=-1+2\sqrt{-1},c_2''=c_3''=c_4''=\frac{1+2\sqrt{-1}}{5}$$
or put it another way below.
$$c_1=c_5=\frac{1+\sqrt{-1}}{4},c_2=c_3=c_4=2\sqrt{-1}$$

Let's simplify (\ref{e16}) and (\ref{e17}) variable substitution. Take
\begin{eqnarray*}
  a&=&\frac{c_3''}{c_2''}\\
  b&=&c_3''\\
  c&=&c_1'c_3''
\end{eqnarray*}
 we have
 \begin{eqnarray}
    c_1'&=&\frac{c}{b}\label{eca}\\
    c_2''&=&\frac{b}{a}\nonumber\\
    c_3''&=&b\nonumber\\
    c_4''&=&ab\nonumber\\
    c_5'&=&\frac{1}{bc}\nonumber
 \end{eqnarray}
 or put it another way below.
 \begin{eqnarray*}
   c_1&=&\frac{b}{b-c}\\
   c_2&=&\frac{b-a}{b}\\
   c_3&=&\frac{b-1}{b}\\
   c_4&=&\frac{ab-1}{ab}\\
   c_5&=&\frac{bc}{bc-1}
 \end{eqnarray*}
  so $(a,b,c)$ determines $c_i$.

 The holonomies of the meridians and longitudes of $W$ (\ref{e16}) can be represented as
\begin{eqnarray}
  m_1&=&\log a\label{e17.01}\\
  l_1 &=&2\log a+2\log(1-ab)-2\log(1-\frac{b}{a})\nonumber\\
  m_2&=&2\log(-c)+\log(1-\frac{1}{bc})-\log(1-\frac{c}{b})\nonumber\\
  l_2&=&2\log(-c)\nonumber
\end{eqnarray}
and the gluing equation (\ref{e17}) is equivalent to
\begin{equation}\label{e17.1}
  2\log b-\log(1-ab)-\log(1-\frac{b}{a})-\log(1-\frac{1}{b})+\log(1-\frac{c}{b})+\log(1-\frac{1}{bc})=0
\end{equation}
 So we can use $(a,b,c)$ to represent the hyperbolic structure equations of $W((p,q),(1,-p'))$.

\subsection{Critical point equation of $V$}\label{secV}

Let's now compute the derivative of $V$.
\begin{eqnarray}
  V_x'(x,y,z,s,m,n,l)&=&\frac{\sqrt{-1}}{q}((p+4q)(2\sqrt{-1}x)-q(4\sqrt{-1}x+2\log(1-e^{2\sqrt{-1}(y+x)})\label{e23}\\
  &&-2\log(1-e^{2\sqrt{-1}(y-x)}))-qk(s,m)2\pi\sqrt{-1})\nonumber\\
  V_y'(x,y,z,s,m,0,l)&=&2\sqrt{-1}(4\sqrt{-1}y-\log(1-e^{2\sqrt{-1}(y+x)})-\log(1-e^{2\sqrt{-1}(y-x)})\label{e24}\\
  &&-\log(1-e^{-2\sqrt{-1}y})+\log(1-e^{2\sqrt{-1}(-y+z)})+\log(1-e^{2\sqrt{-1}(-y-z)}))\nonumber\\
  V_z'(x,y,z,s,m,n,l)&=&2\sqrt{-1}((4\sqrt{-1}z-2\pi\sqrt{-1}+\log(1-e^{2\sqrt{-1}(-y-z)})-\log(1-e^{2\sqrt{-1}(-y+z)}))\label{e25}\\
  &&+(p'-\frac{1}{2})(4\sqrt{-1}z-2\pi\sqrt{-1})+(p'+l+1)2\pi\sqrt{-1})\nonumber
\end{eqnarray}

By comparing (\ref{e23}) (\ref{e24}) (\ref{e25}) and (\ref{e16}) (\ref{e17}),  we found that the system of critical point equations of $V(x,y,z,s,m,0,l)$ is equivalent to the system of hyperbolic structure equations of $\mathcal{W}(\frac{(p+4q,-q)}{qk(s,m)},\frac{(1,p'-\frac{1}{2})}{p'+l+1})$. We can relate these two sets of equation by the following variable substitution. For example, if $(s,m)=(s^+,m^+),l=-p'-2$, we can take $a=e^{2\sqrt{-1}x},b=e^{2\sqrt{-1}y},c=e^{2\sqrt{-1}z}$. Then we have $V_y'[x,y,z,s^+,m^+,0,-p'-2]=0$ is equivalent to the gluing equation (\ref{e17}), and $V_x',V_z'$ corresponding to the Dehn-surgery equation of $\mathcal{W}((p+4q,-q),(1,p'-\frac{1}{2}))$.
\begin{eqnarray}
    V_x'(x,y,z,s^+,m^+,0,-p'-2)&=&\frac{\sqrt{-1}}{q}((p+4q)m_1-ql_1-2\pi\sqrt{-1})\label{eV}\\
    V_z'(x,y,z,s^+,m^+,0,-p'-2)&=&2\sqrt{-1}(m_2+(p'-\frac{1}{2})l_2-2\pi\sqrt{-1})\nonumber
\end{eqnarray}

\subsection{Sister manifolds and sister potential functions}\label{sec5.3}
Take $(s,m)=(s^+,m^+)$ or $(s^-,m^-)$ and $l=-p'$ or $-p'-2$, then we have $|qk(s,m)|=|p'+l+1|=1$. In this condition $\mathcal{W}(\frac{(p,q)}{qk(s,m)},\frac{(1,-p')}{p'+l+1})=\mathcal{W}((p,q),(1,-p'))$ is a normal manifold without strange cone-angle structure.

To prove the volume conjecture, the most straightforward idea is to prove the system of critical point equations of $V(x,y,z,s,m,0,l)$ is equivalent to the system of hyperbolic structure equations of $\mathcal{W}((p,q),(1,-p'))$, and then prove
$$V(x_0,y_0,z_0,s,m,0,l)=\cv(\mathcal{W}((p,q),(1,-p')))$$
where $(x_0,y_0,z_0)$ is the critical point of $V$.

However, the above calculation of the critical point equations of $V(x,y,z,s,m,0,l)$ lead us to prove the system of critical point equations of $V(x,y,z,s,m,0,l)$ is equivalent to the system of hyperbolic structure equations of $\mathcal{W}((p+4q,-q),(1,p'-\frac{1}{2}))$, and then we prove
$$V(x_0,y_0,z_0,s,m,0,l)=\cv(\mathcal{W}((p+4q,-q),(1,p'-\frac{1}{2})))$$
 Hodgson, Meyerhoff, and Weeks \cite{HMW} found a relation bewteen $\mathcal{W}((p,q),(m,l))$ and $\mathcal{W}((p+4q,-q),(m,-l-\frac{m}{2}))$ (see Figure \ref{sister manifolds}).
 if we can prove
 $$\cv(\mathcal{W}((p+4q,-q),(1,p'-\frac{1}{2})))=\cv(\mathcal{W}((p,q),(1,-p')))$$
by this result, we will get our goal.

Unfortunately, the hyperbolic structure of $\mathcal{W}((p+4q,-q),(1,p'-\frac{1}{2}))$ is difficult to describe.  $\mathcal{W}((p+4q,-q),(1,p'-\frac{1}{2}))$ have a cone-angle $4\pi$, which is not a hyperbolic manifold in the usual sense.

The purpose of this section is to find a sister potential function $W$ of $V$ to circumvent the above difficulties. We will define a function $W$ which is similar to $V$, such that $W(x'_0,y'_0,z'_0)=V(x_0,y_0,z_0)$, where $(x'_0,y'_0,z'_0)$ is the critical point of $W$. $W$ also has the geometric properties we expect, the system of critical point equations of $W$ is equivalent to the system of hyperbolic structure equations of $\mathcal{W}((p,q),(1,-p'))$, and in Section \ref{sec5.4} we will prove
$$W(x'_0,y'_0,z'_0)=\cv(\mathcal{W}((p,q),(1,-p')))$$

At first we study the gluing equation (\ref{e17.1}).

Taking the exponential, we have
\begin{eqnarray}
(\ref{e17.1})&\Rightarrow&\frac{b^2(1-\frac{1}{bc})(1-\frac{c}{b})}{(1-\frac{1}{b})(1-ab)(1-\frac{b}{a})}=1\nonumber\\
&\Leftrightarrow&b^2(c-a)(ac-1)=c(ab-1)(b-a)\nonumber\\
&\Leftrightarrow&\frac{1}{b^2}-(a+\frac{1}{a})\frac{1}{b}+(1+a+\frac{1}{a}-c-\frac{1}{c})=0 \label{e28}
\end{eqnarray}
this is a quadratic equation of $b$ (or $\frac{1}{b}$), there are two solutions $b_1(a,c),b_2(a,c)$ of this equation such that
\begin{eqnarray}
    \frac{1}{b_1}+\frac{1}{b_2}&=&a+\frac{1}{a}\label{e26}\\
    \frac{1}{b_1b_2}&=&1+a+\frac{1}{a}-c-\frac{1}{c}\label{e27}
\end{eqnarray}

By simple calculation we have
\begin{eqnarray}
    (1-\frac{b_1}{a})(1-\frac{b_2}{a})&=&(1-ab_1)(1-ab_2)\label{b1b2}\\
    (-c)^3(1-\frac{1}{b_1c})(1-\frac{1}{b_2c})&=&(1-\frac{c}{b_1})(1-\frac{c}{b_2})\nonumber
\end{eqnarray}

We have the following theorem corresponding to the result of  Hodgson, Meyerhoff, and Weeks \cite{HMW} (Figure \ref{sister manifolds}).
\begin{thm}\label{thmsister}
    Suppose $b_1,b_2$ are the solutions of  (\ref{e28}). Then we have
    \begin{eqnarray*}
      &&(a,b_1,c)\text{ is the solution of the hyperbolic structure equations (\ref{eD}) (\ref{e17.1}) of }\mathcal{W}((p,q),(1,-p'))\\
      &\Leftrightarrow&(a,b_2,c)\text{ is a solution of the hyperbolic structure equations of }\mathcal{W}((p+4q,-q),(1,p'-\frac{1}{2}))
    \end{eqnarray*}

\end{thm}

\begin{proof}[Proof of Theorem \ref{thmsister}]
    The gluing equations (\ref{e17.1}) of  $\mathcal{W}((p,q),(1,-p'))$ and $\mathcal{W}((p+4q,-q),(1,p'-\frac{1}{2}))$ are both equivalent to (\ref{e28}).

    Denote
    \begin{eqnarray}
        m_1'&=&\log a\\
         l_1'&=&2\log a+2\log(1-ab_2)-2\log(1-\frac{b_2}{a})\nonumber\\
        m_2'&=&2\log(-c)+\log(1-\frac{1}{b_2c})-\log(1-\frac{c}{b_2})\nonumber\\
        l_2'&=&2\log(-c)\nonumber
    \end{eqnarray}
    then we have the Dehn-surgery equations of $\mathcal{W}((p+4q,-q),(1,p'-\frac{1}{2}))$ is equivalent to
    \begin{eqnarray*}
      (p+4q)m_1'-ql_1'&=&2\pi\sqrt{-1}\\
      m_2'+(p'-\frac{1}{2})l_2'&=&-2\pi\sqrt{-1}
    \end{eqnarray*}

    we have
    \begin{eqnarray*}
      &&(p+4q)m_1'-ql_1'-(pm_1+ql_1)\\
      &=&q(4\log a-l_1-l_1')\\
      &=&2q(\log(1-\frac{b_1}{a})+\log(1-\frac{b_2}{a})-\log(1-ab_1)-\log(1-ab_2))\\
      &=&0\\
      &&m_2'+(p'-\frac{1}{2})l_2'+(m_2-p'l_2)\\
      &=&m_2+m_2'-\log(-c)\\
      &=&3\log(-c)+\log(1-\frac{1}{b_1c})+\log(1-\frac{1}{b_2c})-\log(1-\frac{c}{b_1})-\log(1-\frac{c}{b_2})\\
      &=&0
    \end{eqnarray*}
    Through the above two relations, we have completed the proof.
\end{proof}

\begin{de}[Sister potential function $W$]
    \begin{eqnarray*}
      W(x,y,z,s,m,n,l)=\li(e^{2\sqrt{-1}(y+x)})+\li(e^{2\sqrt{-1}(y-x)})+\li(e^{2\sqrt{-1}(-y+z)})+\li(e^{2\sqrt{-1}(-y-z)})\\
      -\li(e^{-2\sqrt{-1}y})+\frac{p+2q}{q}x^2-4y^2+4(p'-1)z^2-2\pi k(s,m)x-4\pi ny+4\pi l z-\frac{\pi^2}{2}-K(s)\pi^2\nonumber
    \end{eqnarray*}
\end{de}

Similiarly as $V$, we have the following symmetries of $W$.

\begin{pr}\label{prW}
    For $(s,m)$ and $(s',m')$ such that $k(s,m)+k(s',m')=0$, we have
    \begin{enumerate}[(1)]
        \item
            \begin{equation}\label{ec2}
              \overline{W(x,y,z,s,m,n,l)}=W(-\overline{x},-\overline{y},-\overline{z},s',m',-n,-l)-(K(s)-K(s'))\pi^2
            \end{equation}
      \item
        \begin{equation}\label{e29}
          W(x,y,z,s,m,n,l)=W(-x,y,z,s',m',n,l)-(K(s)-K(s'))\pi^2
        \end{equation}
      \item
        \begin{equation}
          W(x,y,z,s,m,n,l)=W(x,y,\pi-z,s,m,n,-2p'+2-l)-4(-p'-l+1)\pi^2
        \end{equation}
    \end{enumerate}
\end{pr}

Let's discuss the critical point equations of $W$.

\begin{eqnarray}
  W_x'(x,y,z,s,m,n,l)&=&\frac{\sqrt{-1}}{q}(-p(2\sqrt{-1}x)-q(4\sqrt{-1}x+2\log(1-e^{2\sqrt{-1}(y+x)})\label{eWx}\\
  &&-2\log(1-e^{2\sqrt{-1}(y-x)}))+qk(s,m)2\pi\sqrt{-1})\nonumber\\
  W_y'(x,y,z,s,m,n,l)&=&V_y'(x,y,z,s,m,l,n)\label{eWy}\\
  W_z'(x,y,z,s,m,n,l)&=&2\sqrt{-1}((4\sqrt{-1}z-2\pi\sqrt{-1}+\log(1-e^{2\sqrt{-1}(-y-z)})-\log(1-e^{2\sqrt{-1}(-y+z)}))\label{eWz}\\
  &&-p'(4\sqrt{-1}z-2\pi\sqrt{-1})+(-p'-l+1)2\pi\sqrt{-1})\nonumber
\end{eqnarray}

The following Lemma shows that $W$ corresponding to the geometry of $\mathcal{W}((p,q),(1,-p'))$ directly.

\begin{lem}\label{=}
    The system of critical point equations of $W$ is equivalent to the system of hyperbolic structure equations of $\mathcal{W}((p,q),(1,-p'))$ (\ref{e17.01}) (\ref{e17.1}) as followings:
    \begin{enumerate}[(1)]
      \item For $W(x,y,z,s^+,m^+,0,-p')$, let $a=e^{2\sqrt{-1}x},b=e^{2\sqrt{-1}y},c=e^{-2\sqrt{-1}z}$;
      \item For $W(x,y,z,s^+,m^+,0,-p'+2)$, let $a=e^{2\sqrt{-1}x},b=e^{2\sqrt{-1}y},c=e^{2\sqrt{-1}z}$;
      \item For $W(x,y,z,s^-,m^-,0,-p')$, let $a=e^{-2\sqrt{-1}x},b=e^{2\sqrt{-1}y},c=e^{-2\sqrt{-1}z}$;
      \item For $W(x,y,z,s^-,m^-,0,-p'+2)$, let $a=e^{-2\sqrt{-1}x},b=e^{2\sqrt{-1}y},c=e^{2\sqrt{-1}z}$.
    \end{enumerate}
\end{lem}

     Similarly as (\ref{e28}) (\ref{e26}) (\ref{e27}), suppose $x,y_1,z$ and $x,y_2,z$ are the solutions of $V'_y(x,y,z,s,m,0,l)=0$(we assume that $a=e^{\pm2\sqrt{-1}x}$,$b=e^{2\sqrt{-1}y}$,$c=e^{\pm2\sqrt{-1}z}$).
     Then $y_1(x,z), y_2(x,z)$ are analytic function of $x,z$.

     Similarly as (\ref{b1b2}), we have

      $\log(1-e^{2\sqrt{-1}(y_1-x)})+\log(1-e^{2\sqrt{-1}(y_2-x)})=\log(1-e^{2\sqrt{-1}(y_1+x)})+\log(1-e^{2\sqrt{-1}(y_2+x)})$

      $6\sqrt{-1}(z-\frac{\pi}{2})+\log(1-e^{2\sqrt{-1}(-y_1-z)})+\log(1-e^{2\sqrt{-1}(-y_2-z)})=\log(1-e^{2\sqrt{-1}(-y_1+z)})+\log(1-e^{2\sqrt{-1}(-y_2+z)})$

We have discussed the symmetries about $x$ and $z$ of potential functions in Proposition \ref{pr3} and Proposition \ref{prW}. The following lemma is more difficult to detect, which is kind of symmetries about $y$ of potential functions.

\begin{lem}\label{V+W}
        For $x,z,y_1,y_2$ as above, we have
        $$V(x,y_1,z,s,m,0,l)+W(x,y_2,z,s,m,0,l+2)=0$$
\end{lem}

\begin{proof}[Proof of Lemma \ref{V+W}]
    \begin{eqnarray*}
      &&V'_x(x,y_1,z,s,m,0,l)+W'_x(x,y_2,z,s,m,0,l+2)\\
      &=&2\sqrt{-1}(\log(1-e^{2\sqrt{-1}(y_1-x)})+\log(1-e^{2\sqrt{-1}(y_2-x)})-\log(1-e^{2\sqrt{-1}(y_1+x)})-\log(1-e^{2\sqrt{-1}(y_2+x)}))\\
      &=&0\\
      &&V'_z(x,y_1,z,s,m,0,l)+W'_z(x,y_2,z,s,m,0,l+2)\\
      &=&2\sqrt{-1}(6\sqrt{-1}(z-\frac{\pi}{2})+\log(1-e^{2\sqrt{-1}(-y_1-z)})+\log(1-e^{2\sqrt{-1}(-y_2-z)})-\log(1-e^{2\sqrt{-1}(-y_1+z)})\\
      &&-\log(1-e^{2\sqrt{-1}(-y_2+z)}))\\
      &=&0
    \end{eqnarray*}
    so we have
    $$V(x,y_1,z,s,m,0,l)+W(x,y_2,z,s,m,0,l+2)=C(s,m,l)$$
    which is independent of $x,z$. Let $x=0,z=\frac{\pi}{2}$, we have $$y_{1,2}(0,\frac{\pi}{2})=\pm\frac{\arctan2}{2}+\frac{\log5}{4}\sqrt{-1}$$
    \begin{eqnarray*}
      C(s,m,l)&=&V(0,y_1,\frac{\pi}{2},s,m,0,l)+W(0,y_2,\frac{\pi}{2},s,m,0,l+2)\\
      &=&2\li(e^{2\sqrt{-1}y_1})+2\li(-e^{-2\sqrt{-1}y_1})-\li(e^{-2\sqrt{-1}y_1})-4y_1^2+2\li(e^{2\sqrt{-1}y_2})+2\li(-e^{-2\sqrt{-1}y_2})\\
      &&-\li(e^{-2\sqrt{-1}y_2})-4y_2^2+\frac{\pi^2}{2}\\
      &\stackrel{\footnote{See Appendix Lemma \ref{f(y0)}}}{=}&0
    \end{eqnarray*}
\end{proof}

%

\subsection{Potential functions and Volume}\label{sec5.4}
The goal of this section is that the critical values of $V$ equal to the $\sqrt{-1}\cv(M)$. To prove this we use the Neumann-Zagier-Yoshida Theory\cite{NZ}\cite{Yoshida} on $W$, and then use the relationship between $V$ and $W$ established in Section \ref{sec5.3}.

    By Thurston's notes \cite{Thurston}, for $(p,q,p')$ doesn't mention in Lemma \ref{nonhyperbolic2}, there is a unique solution $(a,b,c)$(or $(c_1,c_2,c_3,c_4,c_5)$) of hyperbolic structure equations (\ref{eD}) (\ref{e17.1})(or (\ref{eD}) (\ref{e17})) of $\mathcal{W}((p,q),(1,-p'))$ with $\im c_i>0$. $b$ is a solution of (\ref{e28}). WLOG, we denote $b=b_1$ and $b_2$ is the another solution. We have $\im b=\im c_3''>0$.

    Take
    \begin{eqnarray}
      x_0&=&\frac{\log a}{2\sqrt{-1}}\label{ex0}\\
      y_1&=&\frac{\log b_1}{2\sqrt{-1}}\nonumber\\
      y_2&=&\frac{\log b_2}{2\sqrt{-1}}\nonumber\\
      z_0&=&\frac{\log(-c)}{2\sqrt{-1}}+\frac{\pi}{2}\nonumber
    \end{eqnarray}
    we have $(x_0,y_1,z_0)$ is a critical point of $W(x,y,z,s^+,m^+,0,-p'+2)$, $-\frac{\pi}{2}<\re(x_0)<\frac{\pi}{2}$,$0<\re(y_1)<\frac{\pi}{2}, 0<\re(z_0)<\pi$ because of $\im c_i>0$, and
    \begin{eqnarray*}
      m_1&=&2\sqrt{-1}x_0\\
      l_2&=&4\sqrt{-1}z_0-2\pi\sqrt{-1}\\
      pm_1+ql_1&=&2\pi\sqrt{-1}\\
      m_2-p'l_2&=&2\pi\sqrt{-1}
    \end{eqnarray*}

    The complete hyperbolic structure of the Whitehead link complement corresponding to
    $$m_1=l_1=m_2=l_2=0\Leftrightarrow x_0=0,y_1=\frac{\log\frac{1+2\sqrt{-1}}{5}}{2\sqrt{-1}}=\frac{\arctan2}{2}+\frac{\log5}{4}\sqrt{-1},z_0=\frac{\pi}{2}$$

    \begin{pr}\label{Wvol}
        $$W(x_0,y_1,z_0,s^+,m^+,0,-p'+2)\equiv\sqrt{-1}(\vol-\sqrt{-1}\cs)(\mathcal{W}((p,q),(1,-p')))\mod\pi^2$$
    \end{pr}

    In order to prove the above proposition, we need the following theorem of Yoshida to describe the characteristic of complex volume.
    \begin{thm}[Yoshida \cite{Yoshida} 1985]\label{thmYoshida}
        Suppose $M$ is a oriented complete hyperbolic $3$-manifold of finite volume with $h$ cusps. Let $U$ be the deformation space of the hyperbolic structure on $M$ and $\vec{u}\in U$ be the parameter\footnote{$\vec{u}=(u_1,\cdots,u_h)$, where $u_i=m_i$ or $l_i$, and $m_i,l_i$ are holonomies of the meridian or longitude of the $i^{th}$ cusp. $\vec{u_0}=\vec{0}\in U$ is the point representing the original complete hyperbolic structure on $M$.} of $U$. Then there is a analytic function $f(\vec{u})$ of $\vec{u}$ on a neighborhood $V$ of $\vec{0}$, such that if $M_{\vec{u}}$ is the closed hyperbolic manifold obtained by doing hyperbolic Dehn filling $h$ geodesic loops $\gamma_1,\cdots,\gamma_h$ to the $h$ cusps of $M$, then $\vec{u}\in V$ and
        $$f(\vec{u})\equiv\vol(M_{\vec{u}})-\sqrt{-1}\cs(M_{\vec{u}})\footnote{The Chern-Simons invariant in our paper is $-2\pi^2$ times the Chern-Simons invariant in \cite{Yoshida}}+\frac{\pi}{2}\sum_{i=1}^{h}\mathrm{H}(\gamma_i)\mod\pi^2\sqrt{-1}$$
        where $\mathrm{H}(\gamma_i)$ is the holonomy of $\gamma_i$.
    \end{thm}

    Consider $\mathcal{W}((p,q),(1,-p'))$ with core curves $\gamma_1$ and $\gamma_2$. Choose $\vec{u}=(m_1,l_2)$. By Theorem \ref{thmYoshida}, we have an analytic function $f(m_1,l_2)$, such that
    $$f(m_1,l_2)\equiv\sqrt{-1}(\vol-\sqrt{-1}\cs)(\mathcal{W}((p,q),(1,-p')))+\frac{\pi\sqrt{-1}}{2}(\mathrm{H}(\gamma_1)+\mathrm{H}(\gamma_2))\mod\pi^2$$

    We have $\mathrm{H}(\gamma_1)\equiv-\tilde{q}m_1+\tilde{p}l_1\equiv\frac{-m_1+2\pi\sqrt{-1}\tilde{p}}{q}\mod2\pi\sqrt{-1}$, $\mathrm{H}(\gamma_2)\equiv l_2\mod2\pi\sqrt{-1}$ and
    \begin{eqnarray*}
      &&pm_1+ql_1=2\pi\sqrt{-1}\\
      &\Leftrightarrow&\frac{2\sqrt{-1}(px_0-\pi)}{q}=-l_1\\
      &\Leftrightarrow&\frac{px_0^2-\pi x_0}{q}=\frac{-l_1}{2\sqrt{-1}}\frac{m_1}{2\sqrt{-1}}\\
      &\Rightarrow&\frac{px_0^2-\pi x_0}{q}=\frac{m_1l_1}{4}\\
      &&m_2-p'l_2=2\pi\sqrt{-1}\\
      &\Leftrightarrow&p'(4z_0-2\pi)=\frac{m_2}{\sqrt{-1}}-2\pi
    \end{eqnarray*}

    Now we can simplify $W$ as below. We're going to eliminate $p,q,p'$.

    \begin{eqnarray*}
      &&W(x_0,y_1,z_0,s^+,m^+,0,-p'+2)+\frac{\pi\sqrt{-1}}{2}(\mathrm{H}(\gamma_1)+\mathrm{H}(\gamma_2))\\
      &\equiv&\li(e^{2\sqrt{-1}(y_1+x_0)})+\li(e^{2\sqrt{-1}(y_1-x_0)})+\li(e^{2\sqrt{-1}(-y_1+z_0)})+\li(e^{2\sqrt{-1}(-y_1-z_0)})-\li(e^{-2\sqrt{-1}y_1})\\
      &&+\frac{p+2q}{q}x_0^2-4y_1^2+4(p'-1)z_0^2-\frac{2\pi x_0}{q}-4\pi(p'-2)z_0-\frac{\pi^2}{2}-K(s^+)\pi^2+\frac{\pi x_0}{q}-\frac{\tilde{p}\pi^2}{q}+\frac{\pi\sqrt{-1}l_2}{2}\\
      &\equiv&\li(e^{2\sqrt{-1}(y_1+x_0)})+\li(e^{2\sqrt{-1}(y_1-x_0)})+\li(e^{2\sqrt{-1}(-y_1+z_0)})+\li(e^{2\sqrt{-1}(-y_1-z_0)})-\li(e^{-2\sqrt{-1}y_1})\\
      &&-4y_1^2+2x_0^2+\frac{px_0^2-\pi x_0}{q}+p'(z_0(4z_0-2\pi)-\pi(2z_0-\pi)-\pi^2)-4z_0^2+8\pi z_0-\frac{\pi^2}{2}+\frac{\pi\sqrt{-1}l_2}{2}\\
      &\equiv&\li(e^{2\sqrt{-1}(y_1+x_0)})+\li(e^{2\sqrt{-1}(y_1-x_0)})+\li(e^{2\sqrt{-1}(-y_1+z_0)})+\li(e^{2\sqrt{-1}(-y_1-z_0)})-\li(e^{-2\sqrt{-1}y_1})\\
      &&-4y_1^2-\frac{m_1^2}{2}+\frac{m_1l_1}{4}+(-\frac{\sqrt{-1}l_2}{4}+\frac{\pi}{2})(\frac{m_2}{\sqrt{-1}}-2\pi)-\pi\frac{\frac{m_2}{\sqrt{-1}}-2\pi}{2}-4(-\frac{\sqrt{-1}l_2}{4}+\frac{\pi}{2})^2+8\pi (-\frac{\sqrt{-1}l_2}{4}+\frac{\pi}{2})\\
      &&-\frac{\pi^2}{2}+\frac{\pi\sqrt{-1}l_2}{2}\\
      &\equiv&\li(e^{2\sqrt{-1}(y_1+x_0)})+\li(e^{2\sqrt{-1}(y_1-x_0)})+\li(e^{2\sqrt{-1}(-y_1+z_0)})+\li(e^{2\sqrt{-1}(-y_1-z_0)})-\li(e^{-2\sqrt{-1}y_1})\\
      &&-4y_1^2-\frac{m_1^2}{2}+\frac{m_1l_1}{4}+\frac{l_2^2}{4}-\frac{m_2l_2}{4}-\frac{\pi^2}{2}\mod\pi^2
    \end{eqnarray*}

    Define
    \begin{eqnarray*}
        \tilde{f}(m_1,l_2)&=&\li(e^{2\sqrt{-1}(y_1+x)})+\li(e^{2\sqrt{-1}(y_1-x)})+\li(e^{2\sqrt{-1}(-y_1+z)})+\li(e^{2\sqrt{-1}(-y_1-z)})-\li(e^{-2\sqrt{-1}y_1})\\
        &&-4y_1^2-\frac{m_1^2}{2}+\frac{m_1l_1}{4}+\frac{l_2^2}{4}-\frac{m_2l_2}{4}-\frac{\pi^2}{2}
    \end{eqnarray*}
    $x=-\frac{\sqrt{-1}m_1}{2},z=-\frac{\sqrt{-1}l_2}{4}+\frac{\pi}{2},y_1(m_1,l_2)=y_1(x,z),l_1,m_2$ are analytic function of $m_1,l_2$, so $\tilde{f}(m_1,l_2)$ is an  analytic function of $m_1,l_2$. To prove Proposition \ref{Wvol}, it suffices to prove $f(m_1,l_2)\equiv\tilde{f}(m_1,l_2)\mod\pi^2$.

    \begin{pr}\label{im!=0}
        $$\im x_0\ne0, \im z_0\ne0$$
    \end{pr}

    \begin{proof}[Proof of Proposition \ref{im!=0}]
      \begin{eqnarray*}
        \mathrm{Length}(\gamma_1)&=&\re \mathrm{H}(\gamma_1)\\
        &=&\re(\frac{-m_1}{q})\\
        &=&\frac{2\im x_0}{q}>0\\
        \mathrm{Length}(\gamma_2)&=&\re \mathrm{H}(\gamma_2)\\
        &=&\re(l_2)\\
        &=&-4\im z_0>0
      \end{eqnarray*}
    \end{proof}

    Proposition \ref{im!=0} will be used in Section \ref{subsecxz}.

    \begin{pr}\label{W=vol}
        $$\im W(x_0,y_1,z_0,s^+,m^+,0,-p'+2)=\vol(\mathcal{W}((p,q),(1,-p')))$$
    \end{pr}

    \begin{proof}[Proof of Proposition \ref{W=vol}]
        \begin{eqnarray*}
          &&\im W(x_0,y_1,z_0,s^+,m^+,0,-p'+2)\\
          &\stackrel{Lemma\ \ref{lem2}}{=}&\dd(e^{2\sqrt{-1}(y_1+x_0)})+\dd(e^{2\sqrt{-1}(y_1-x_0)})+\dd(e^{2\sqrt{-1}(-y_1+z_0)})+\dd(e^{2\sqrt{-1}(-y_1-z_0)})-\dd(e^{-2\sqrt{-1}y_1})\\
          &\stackrel{(\ref{eca})\ (\ref{ex0})}{=}&\dd(c_4'')+\dd(c_2'')+\dd(c_1')+\dd(c_5')-\dd(\frac{1}{c_3''})\\
          &\stackrel{Proposition\ \ref{prD}}{=}&\dd(c_1')+\dd(c_2'')+\dd(c_3'')+\dd(c_4'')+\dd(c_5')\\
          &\stackrel{Proposition\ \ref{prvol}}{=}&\vol(\mathcal{W}((p,q),(1,-p')))
        \end{eqnarray*}
    \end{proof}

    \begin{cor}\label{imf}
      $$\im f(m_1,l_2)=\im \tilde{f}(m_1,l_2)$$
    \end{cor}

    \begin{lem}\label{f(0,0)}
        $\tilde{f}(0,0)\equiv\sqrt{-1}(\vol(\mathcal{W})-\sqrt{-1}\cs(\mathcal{W}))\mod\pi^2$
    \end{lem}

    \begin{proof}[Proof of Lemma \ref{f(0,0)}]
     \begin{eqnarray*}
        \tilde{f}(0,0)&=&(\li(e^{2\sqrt{-1}(y_1+x)})+\li(e^{2\sqrt{-1}(y_1-x)})+\li(e^{2\sqrt{-1}(-y_1+z)})+\li(e^{2\sqrt{-1}(-y_1-z)})-\li(e^{-2\sqrt{-1}y_1})\\
        &&-4y_1^2-\frac{\pi^2}{2})|_{(x,y_1,z)=(0,\frac{\arctan2}{2}+\frac{\log5}{4}\sqrt{-1},\frac{\pi}{2})}\\
        &=&(2\li(e^{2\sqrt{-1}y_1})+\li(-e^{-2\sqrt{-1}y_1})-\li(e^{-2\sqrt{-1}y_1})-4y_1^2-\frac{\pi^2}{2})|_{y_1=\frac{\arctan2}{2}+\frac{\log5}{4}\sqrt{-1}}\\
        &\equiv&\frac{\pi^2}{4}+4\dd(\sqrt{-1})\sqrt{-1}\\
        &\equiv&\sqrt{-1}(\vol(\mathcal{W})-\sqrt{-1}\cs(\mathcal{W}))\mod\pi^2
     \end{eqnarray*}

     The penultimate "$\equiv$" is because of Appendix Lemma \ref{f(y0)}. The last "$\equiv$" is from \cite{Mutation} Proof of Proposition 2.6 in page 109. The Chern-Simons invariant in our paper is $-2\pi^2$ times the Chern-Simons invariant in \cite{Mutation}.
    \end{proof}

    \begin{cor}\label{f(0,0)=}
        $$f(0,0)\equiv\tilde{f}(0,0)\mod\pi^2$$
    \end{cor}

    \begin{proof}[Proof of Proposition \ref{Wvol}]
        Consider $f(m_1,l_2)$ and $\tilde{f}(m_1,l_2)$. They are both analytic function. By Corollary \ref{imf} we know they have the same imaginary part. By Corollary \ref{f(0,0)=} we know they have the same initial value. So       $$f(m_1,l_2)\equiv\tilde{f}(m_1,l_2)\mod\pi^2$$
    \end{proof}

    Now we use the relationship between $V$ and $W$ established in Section \ref{sec5.3} to prove the critical values of $V$ equal to the $\sqrt{-1}\cv(M)$. We have
    \begin{eqnarray*}
      &&W(x_0,y_1,z_0,s^+,m^+,0,-p'+2)\\
      &\stackrel{\text{Lemma }\ref{V+W}}{=}&-V(x_0,y_2,z_0,s^+,m^+,0,-p')\\
      &\stackrel{(\ref{ec1})}{\equiv}&-\overline{V(-\overline{x_0},-\overline{y_2},-\overline{z_0},s^-,m^-,0,p'-1)}\\
      &\stackrel{(\ref{eVz})}{\equiv}&-\overline{V(-\overline{x_0},-\overline{y_2},\pi-\overline{z_0},s^-,m^-,0,-p'-2)}\mod\pi^2
    \end{eqnarray*}
    where $x_0,y_1,z_0$ and $x_0,y_2,z_0$ are the solutions of $W'_y(x,y,z,s^+,m^+,0,-p'+2)=0$.

    So we have
    \begin{pr}\label{Vol}
        \begin{eqnarray*}
          &&V(-\overline{x_0},-\overline{y_2},\pi-\overline{z_0},s^-,m^-,0,-p'-2)\\
          &\equiv&V(\overline{x_0},-\overline{y_2},\pi-\overline{z_0},s^+,m^+,0,-p'-2)\\
          &\equiv&V(-\overline{x_0},-\overline{y_2},\overline{z_0},s^-,m^-,0,-p')\\
          &\equiv&V(\overline{x_0},-\overline{y_2},\overline{z_0},s^+,m^+,0,-p')\\
          &\equiv&\sqrt{-1}(\vol+\sqrt{-1}\cs)(\mathcal{W}((p,q),(1,-p')))\mod\pi^2
        \end{eqnarray*}
    \end{pr}

\section{Asymptotics of the Reshetikhin-Turaev invariants}\label{sec6}
    The goal of this section is to prove Theorem \ref{thm1} by estimating each of the Fourier cofficients $\check{f_r}(s,m,n,l)$ in Proposition \ref{pr2}.

    By Proposition \ref{pr2.1}, we just need to consider
    $$\int_{D_{\frac{\epsilon}{2}}}g(s,x,y,z)e^{\frac{r}{4 \pi \sqrt{-1}}V(x,y,z,s,m,n,l)}\D x\D y\D z$$
    We can deform $D_{\frac{\epsilon}{2}}$ in $\mathbb{C}^3$ without changing the value of the integral. We will deform $D_{\frac{\epsilon}{2}}$ in Section \ref{subsecy} and Section \ref{subsecxz} in order for $\im V(x,y,z,s,m,n,l)<\im \vol(M)$.

     In Section \ref{subsecy} we estimate the terms $\check{f_r}(s,m,n,l)$ for $n\ne0$.

     In Section \ref{subsecxz} we estimate the terms $\check{f_r}(s,m,0,l)$ for $(s,m)\ne(s^+,m^+),(s^-,m^-)$ or $l>-p'$ or $l<-p'-2$.

     In Section \ref{subsecleading} we estimate the four leading terms $\check{f_r}(s^+,m^+,0,-p')$, $\check{f_r}(s^+,m^+,0,-p'-2)$, $\check{f_r}(s^-,m^-,0,-p')$ and $\check{f_r}(s^-,m^-,0,-p'-2)$.

      Finally in Section \ref{subsecpf} we show that the sum of leading terms has the desired asymptotic behavior and the sum of all the other terms are neglectable, which completes the proof. We will also prove Corollary \ref{thmtv} in the end of Section \ref{subsecpf}.

\subsection{Estimates of the Fourier coefficients $n\ne0$}\label{subsecy}
In this section we will deform the region along the imaginary part of $y$.


Define
$$S_{\text{top}}(Y)=\{(x,y,z)\in\mathbb{C}^3|\re(x,y,z)\in D_{\epsilon}, \im(x,y,z)=(0,Y,0)\}$$
$$S_{\text{side}}(Y)=\{(x,y+tY\sqrt{-1},z)|(x,y,z)\in\partial D_{\epsilon},t\in[0,1]\}$$
and $S(Y)=S_{\text{top}}(Y)\cup S_{\text{side}}(Y)$, then we have
\begin{eqnarray*}
  &&\int_{D_{\frac{\epsilon}{2}}}g(s,x,y,z)e^{\frac{r}{4 \pi \sqrt{-1}}V(x,y,z,s,m,n,l)}\D x\D y\D z\\
  &=&\int_{S(Y)\cup(D_{\frac{\epsilon}{2}}\backslash D_{\epsilon})}g(s,x,y,z)e^{\frac{r}{4 \pi \sqrt{-1}}V(x,y,z,s,m,n,l)}\D x\D y\D z
\end{eqnarray*}

By Lemma \ref{lemim<3.5}, we have
\begin{pr}\label{prD1}
    For $(x,y,z)\in D_{\frac{\epsilon}{2}}\backslash D_{\epsilon}$,
    $$\im V(x,y,z,s,m,n,l)<3.5$$
\end{pr}

So the integral in $D_{\frac{\epsilon}{2}}\backslash D_{\epsilon}$ neglectable.

Suppose $\im(x,z)=(0,0)$. we have
\begin{eqnarray*}
  &&\im V(x,y,z,s,m,n,l)\\
  &=&\im(\li(e^{2\sqrt{-1}(y+x)})+\li(e^{2\sqrt{-1}(y-x)})+\li(e^{2\sqrt{-1}(-y+z)})+\li(e^{2\sqrt{-1}(-y-z)})-\li(e^{-2\sqrt{-1}y}))\\
  &&-4(2\re y+n\pi)\im y
\end{eqnarray*}

    Define
    \begin{eqnarray*}
      \tilde{F}(x,y,z,n)&=&F(x+\re y,\im y)+F(-x+\re y,\im y)+F(z-\re y,-\im y)+F(-z-\re y,-\im y)\\
      &&-F(-\re y,-\im y)-4(2\re y+n\pi)\im y
    \end{eqnarray*}
    where $F(\cdot,\cdot)$ is defined as (\ref{eF}). Corollary \ref{cor1.1} shows that $F$ is close related to $\im V$. We have
    $$\tilde{F}(x,y,z,n)=\begin{cases}
                 -2(2\pi n+2\re y-\pi)\im y, & \mbox{if }\im y\ge0;\\
                 -4\pi(n+1)\im y, & \mbox{if }\im y\le0.
               \end{cases}$$
    So we have
    \begin{enumerate}
      \item if $n\ge1$ and $\im y>0$,
        $$\tilde{F}(x,y,z,n)\le-2\pi(2n-1)\im y\le-2\pi\im y$$
      \item  if $n\le-2$ and $\im y<0$,
        $$\tilde{F}(x,y,z,n)<4\pi\im y$$
    \end{enumerate}
    By Corollary \ref{cor1.1}, we have
    $$|\im V(x,y,z,s,m,n,l)-\tilde{F}(x,y,z,n)|\le\frac{5}{2}$$
    for $|\im y|>\frac{\log5}{4}$. So we have the following lemma.

\begin{lem}\label{lemimy}
    \begin{enumerate}[(1)]
      \item if $n\ge1$, and $\im y=\frac{\log5}{4}$, then
        $$\im V(x,y,z,s,m,n,l)<0$$
      \item if $n\le-1$, and $\im y=-\frac{\log5}{4}$, then
        $$\im V(x,y,z,s,m,n,l)<\frac{5}{2}$$
    \end{enumerate}
\end{lem}


The following part of this section shows the convexity of $\im V$ about $\im y$.

For $x\in \mathbb{C}\backslash\pi\mathbb{Z}$, define
$$h(x)=\frac{1}{e^{2\sqrt{-1}x}-1}$$
then we have
$$\overline{h(x)}=h(-\overline{x})$$
$$h(x)+h(-x)=-1$$
$$\frac{\D^2\li(e^{2\sqrt{-1}x})}{\D x^2}=-4h(-x)$$
\begin{equation}\label{eh}
  \im h(x)=-\im h(-x)=-\frac{e^{-2\im x}\sin2\re x}{|e^{2\sqrt{-1}x}-1|^2}
\end{equation}

\begin{lem}\label{lemyconvex}
    Suppose $\im(x,z)=(0,0)$, $\re(x,y,z)\in D_0$, $|\im y|\le\frac{\log5}{4}$. Then $\im f(x,y,z,s,m,n,l)$ is a strictly concave up function of $\im y$, that is
    $$\frac{\D^2\im f}{(\D\im y)^2}>0$$
\end{lem}

\begin{proof}[Proof of Lemma \ref{lemyconvex}]
    \begin{eqnarray*}
      \frac{\D^2\im f}{(\D\im y)^2}&=&-\im f_{yy}''\\
      &=&4\im(h(y+z)+h(y-z)-h(y+x)-h(y-x)-h(y))
    \end{eqnarray*}

    For $(x,\re y,z)\in D_0$, we have $0<\re y<\frac{\pi}{2}$, so $\sin2\re y>0$. Together with (\ref{eh}), we  have
    \begin{equation}\label{eyconvex1}
      \im h(y)<0
    \end{equation}

    \begin{eqnarray*}
      &&\im(h(y+z)+h(y-z))\\
      &=&-(\frac{e^{-2\im y}\sin2(\re y+z)}{|e^{2\sqrt{-1}(y+z)}-1|^2}+\frac{e^{-2\im y}\sin2(\re y-z)}{|e^{2\sqrt{-1}(y-z)}-1|^2})\\
      &=&\frac{4e^{-4\im y}\sin2\re y(\cos2\re y-\cosh(2\im y)\cos2z)}{|e^{2\sqrt{-1}(y+z)}-1|^2|e^{2\sqrt{-1}(y-z)}-1|^2}
    \end{eqnarray*}

    If $z\in [\frac{\pi}{4},\frac{3\pi}{4}]$, then $\cos 2z\le0$, we have
    \begin{eqnarray*}
      &&\cos2\re y-\cosh(2\im y)\cos2z\\
      &\ge&\cos2\re y-\cos2z\\
      &=&2\sin(z+\re y)\sin(z-\re y)>0
    \end{eqnarray*}

    if $z\in (\frac{\pi}{8},\frac{\pi}{4})\cup(\frac{3\pi}{4},\frac{7\pi}{8})$, then $\cos2z>0$, $\cosh(2\im y)\le\cosh(\frac{\log5}{2})=\frac{3}{\sqrt{5}}$, we have
    \begin{eqnarray*}
      &&\cos2\re y-\cosh(2\im y)\cos2z\\
      &\ge&\cos2\re y-\frac{3}{\sqrt{5}}\cos2z\\
      &>&\cos2\re y-\sqrt{2}\cos2z
    \end{eqnarray*}

    By the symmetry of $z$ and $\pi-z$, we can suppose $z\in (\frac{\pi}{8},\frac{\pi}{4})$, then because of $\re y+\frac{\pi}{4}<2z$ and $\re y>0$, we have $\re y\in(0,\frac{\pi}{4})$ and
    \begin{eqnarray}
      &&\cos2\re y-\sqrt{2}\cos2z\nonumber\\
      &>&\cos2\re y-\sqrt{2}\cos(\re y+\frac{\pi}{4})\nonumber\\
      &\stackrel{\theta=\re y+\frac{\pi}{4}}{=}&\cos(2\theta-\frac{\pi}{2})-\sqrt{2}\cos\theta\nonumber\\
      &=&2\cos\theta(\sin\theta-\frac{\sqrt{2}}{2})\label{etheta1}
    \end{eqnarray}
    because of $\theta\in(\frac{\pi}{4},\frac{\pi}{2})$, we have (\ref{etheta1})$>0$, so
    \begin{equation}\label{eyconvex2}
      \im(h(y+z)+h(y-z))>0
    \end{equation}

    \begin{eqnarray*}
      &&\im(h(y+x)+h(y-x))\\
      &=&\frac{4e^{-4\im y}\sin2\re y(\cos2\re y-\cosh(2\im y)\cos2x)}{|e^{2\sqrt{-1}(y+x)}-1|^2|e^{2\sqrt{-1}(y-x)}-1|^2}
    \end{eqnarray*}
    we have $\cos2x>0$,
    \begin{eqnarray*}
      &&\cos2\re y-\cosh(2\im y)\cos2x\\
      &\le&\cos2\re y-\cos2x\\
      &=&-2\sin(\re y+x)\sin(\re y-x)<0
    \end{eqnarray*}
      so we have
      \begin{equation}\label{eyconvex3}
        \im(h(y+x)+h(y-x))<0
      \end{equation}

      Combine (\ref{eyconvex1}), (\ref{eyconvex2}), (\ref{eyconvex3}), we conclude the proof.
\end{proof}

    By Lemma \ref{lemyconvex}, we have
    \begin{cor}\label{corside1}
    \begin{enumerate}[(1)]
      \item if $n\ge1$, for $(x,y,z)\in S_{\text{side}}(\frac{\log 5}{4})$
    $$\im V(x,y,z,s,m,n,l)\le\max\{\im V(x,\re y,z,s,m,n,l),\im V(x,\re y+\frac{\log 5}{4}\sqrt{-1},z,s,m,n,l)\}<3.5$$
      \item if $n\le-1$, for $(x,y,z)\in S_{\text{side}}(-\frac{\log 5}{4})$,
    $$\im V(x,y,z,s,m,n,l)\le\max\{\im V(x,\re y,z,s,m,n,l),\im V(x,\re y-\frac{\log 5}{4}\sqrt{-1},z,s,m,n,l)\}<3.5$$
    \end{enumerate}
    \end{cor}

    Let $D_{(m,n,l)}=\begin{cases}
                    S(\frac{\log5}{4}), & \mbox{if } n\ge1; \\
                    S(-\frac{\log5}{4}), & \mbox{if } n\le-1.
                  \end{cases}$

    Combine  Lemma \ref{lemimy} with Corollary \ref{corside1}, we have the following Proposition.
    \begin{pr}\label{n!=0}
      For $(x,y,z)\in D_{(m,n,l)}$,
      $$\im V(x,y,z,s,m,n,l)<3.5$$
    \end{pr}


%

\subsection{Estimates of the Fourier coefficients $n=0$}\label{subsecxz}

    Define $l^+=-p',l^-=-p'-2$.

\begin{rem}
    For $l=-p'-1$, by Proposition \ref{cancell} and Remark \ref{rem1}, we have $\check{f_r}=0$.
\end{rem}

    Similiar as before, we have
\begin{eqnarray*}
  &&\int_{D_{\frac{\epsilon}{2}}}g(s,x,y,z)e^{\frac{r}{4 \pi \sqrt{-1}}V(x,y,z,s,m,0,l)}\D x\D y\D z\\
  &=&\int_{S(\frac{\log5}{4})\cup(D_{\frac{\epsilon}{2}}\backslash D_{\epsilon})}g(s,x,y,z)e^{\frac{r}{4 \pi \sqrt{-1}}V(x,y,z,s,m,0,l)}\D x\D y\D z
\end{eqnarray*}
and the integral in $D_{\frac{\epsilon}{2}}\backslash D_{\epsilon}$ neglectable.

Define $D_H=\{(x,y,z)\in\mathbb{R}^3|y-x\in(0,\frac{\pi}{2}),y+x\in(0,\frac{\pi}{2}),z+y\in(\frac{\pi}{2},\pi),z-y\in(0,\frac{\pi}{2})\}$, we have $D_H$ is octahedron and $D_H\in D_0$.

\begin{figure}[H]
  \centering
  \includegraphics[width=0.5\textwidth]{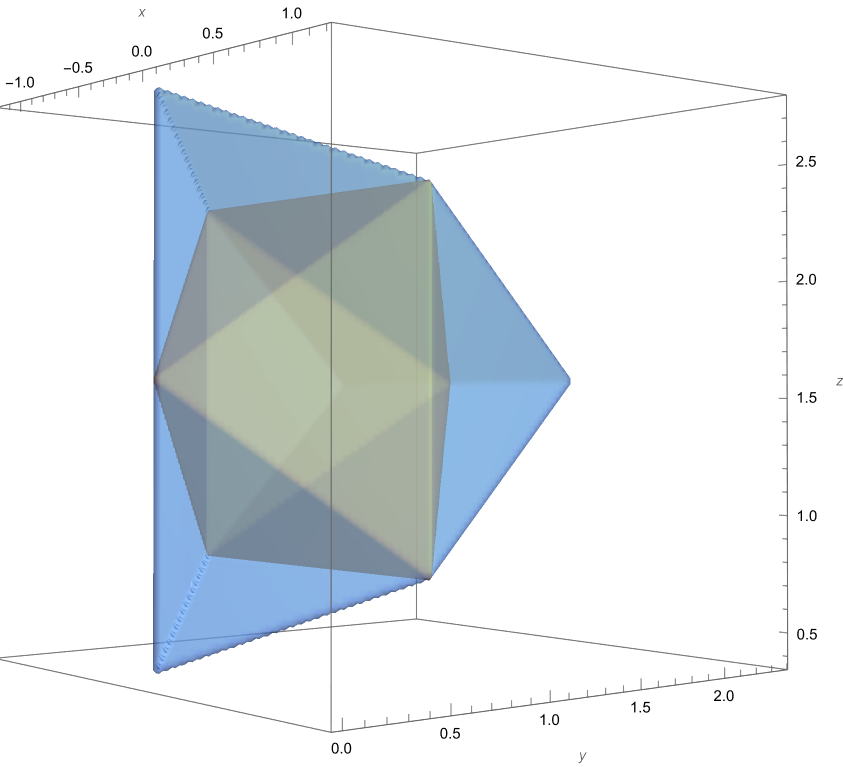}
  \caption{The blue region is $D_0$ and the orange region is $D_H$.}\label{DHD0}
\end{figure}

\begin{lem}\label{H<0}
        If $\re(x,y,z)\in D_H$, then $\im V(x,y,z,s,m,n,l)$ is strictly concave down in $\re x,\re y,\re z$, and is strictly concave up in $\im x,\im y,\im z$.
\end{lem}

\begin{proof}[Proof of Lemma \ref{H<0}]
    We need to prove $\hess_{\re x,\re y,\re z}(\im V)$ is negative definite and $\hess_{\im x,\im y,\im z}(\im V)$ is positive definite. However, we have
    $$\hess_{\re x,\re y,\re z}(\im V)=\im \hess_{x,y,z}(V)=-\hess_{\im x,\im y,\im z}(\im V)$$
    so it suffices to prove $\im \hess_{x,y,z}(V)$ is negative definite.

    \begin{eqnarray*}
      &&\im \hess_{x,y,z}(V)\\
      &=&-4\im\left(\begin{matrix}
                -h(y+x)-h(y-x) & h(y-x)-h(y+x) & 0\\
                h(y-x)-h(y+x) & -h(y+x)-h(y-x)-h(y)+h(y+z)+h(y-z) & h(y-z)-h(y+z) \\
                0 & h(y-z)-h(y+z) & h(y+z)+h(y-z)
              \end{matrix}\right)\\
      &=&-4(\left(\begin{matrix}
                     1 & -1 & 0 \\
                     1 & 1 & 0 \\
                     0 & 0 & 0
                   \end{matrix}\right)\im\left(\begin{matrix}
                     -h(y+x) & 0 & 0 \\
                     0 & -h(y-x) & 0 \\
                     0 & 0 & 0
                   \end{matrix}\right)\left(\begin{matrix}
                     1 & 1 & 0 \\
                     -1 & 1 & 0 \\
                     0 & 0 & 0
                   \end{matrix}\right)+\im\left(\begin{matrix}
                     0 & 0 & 0 \\
                     0 & -h(y) & 0 \\
                     0 & 0 & 0
                   \end{matrix}\right)\\
      &&+\left(\begin{matrix}
                     0 & 0 & 0 \\
                     0 & 1 & -1 \\
                     0 & 1 & 1
                   \end{matrix}\right)\im\left(\begin{matrix}
                     0 & 0 & 0 \\
                     0 & h(y-z) & 0 \\
                     0 & 0 & h(y+z)
                   \end{matrix}\right)\left(\begin{matrix}
                     0 & 0 & 0 \\
                     0 & 1 & 1 \\
                     0 & -1 & 1
                   \end{matrix}\right))
    \end{eqnarray*}
    By (\ref{eh}) and the definition of $D_H$, we have if $\re(x,y,z)\in D_H$, then
    $$h(y+x)<0,h(y-x)<0,h(y)<0,h(y-z)>0,h(y+z)>0$$
    so $\im \hess_{x,y,z}(V)$ is negative definite, we conclude the proof.
\end{proof}

Define $D_{H,\epsilon}=D_H\cap D_\epsilon$.

\begin{lem}\label{imV>3.5}
    $$\{(x,y,z)\in D_0|\im V(x,y+\frac{\log5}{4}\sqrt{-1},z,s,m,0,l)>3.5\}\subset D_{H,\epsilon}$$
\end{lem}

\begin{figure}[H]
      \centering
        \begin{minipage}[t]{0.45\linewidth}
            \null
            \includegraphics[width=1\linewidth]{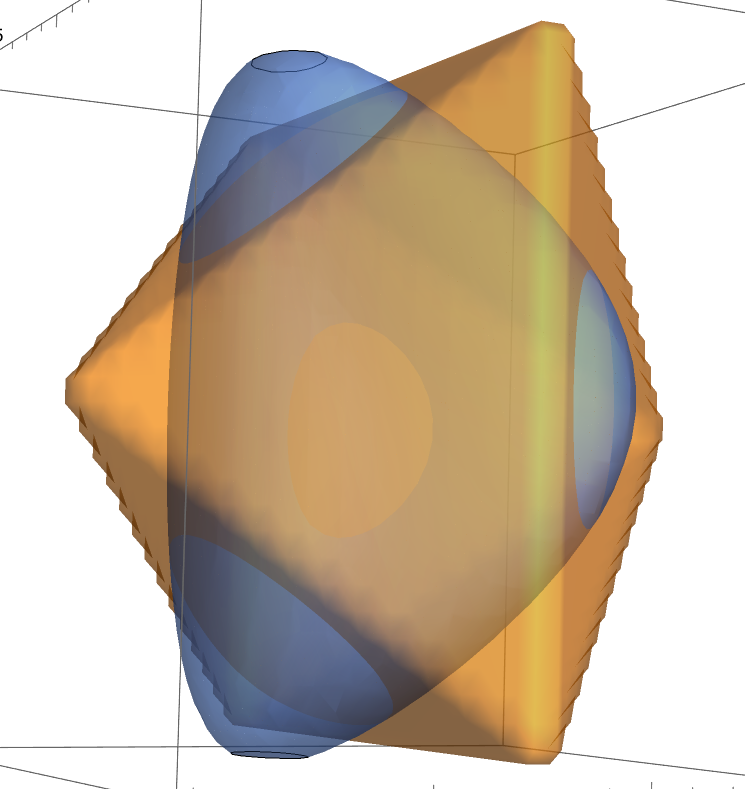}
            \caption{The region $\im V(x,y,z)>3.5$}\label{big}
        \end{minipage}
        \hfill
        \begin{minipage}[t]{0.4\linewidth}
            \null
            \includegraphics[width=1\linewidth]{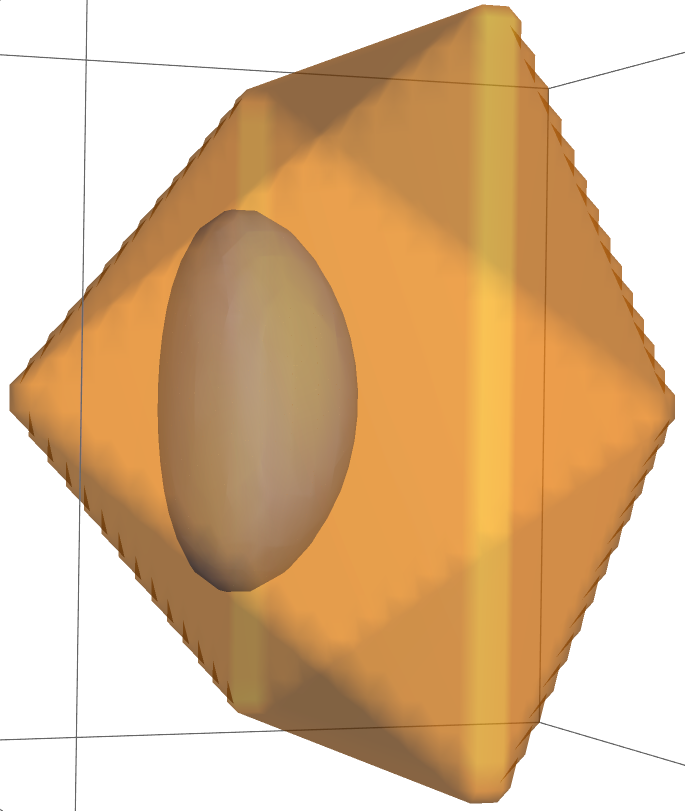}
            \caption{The region $\im V(x,y+\frac{\log5}{4}\sqrt{-1},z)>3.5$}\label{small}
        \end{minipage}
        \hfill
    \end{figure}

    In the following figures $D_H$ is the orange regions, The blue regions correspond to $\im V>3.5$. The reason of why we move $\im y$ with $\frac{\log5}{4}$ is that if we don't move, then the $4$-corners of the blue region(see Figure \ref{big}) is out of the orange region $D_H$, which is a troublesome situation. We can overcome this difficulty by moving $\im y$ with $\frac{\log5}{4}$. The choice of $\frac{\log5}{4}$ is because of $\im y_0=\frac{\log5}{4}$, where $(x_0,y_0,z_0)$ corresponding to the solution of hyperbolic structure equations of $\mathcal{W}$. $\im V$ is tame at the critical point of $V$. For $p,q,p'$ sufficiently large, the critical point of $V$ is going to approach $(x_0,y_0,z_0)$.
    \begin{rem}
       Our $x_0,y_0,z_0$ notation here is just for convenience. And then we will use $x_0,z_0$ in a different sense.
    \end{rem}

    The proof of Lemma \ref{imV>3.5} is in Appendix.

    \begin{cor}\label{cortop2}
        If $(x,y,z)\in \partial S_{\mathrm{top}}(\frac{\log5}{4})$, then
        $$\im V(x,y,z,s,m,0,l)<3.5$$
    \end{cor}
     Combine Corollary \ref{cortop2}, Lemma \ref{lemim<3.5} and Lemma \ref{lemyconvex}, we have for $(x,y,z)\in S_{\text{side}}(\frac{\log 5}{4})$
    $$\im V(x,y,z,s,m,0,l)\le\max\{\im V(x,\re y,z,s,m,0,l),\im V(x,\re y+\frac{\log 5}{4}\sqrt{-1},z,s,m,0,l)\}<3.5$$


    Now we just need to estimate the integral in $S_{\text{top}}(\frac{\log5}{4})$.

    Define
    $$S_{\text{top}}^{\epsilon_1 x,\epsilon_2 z}=\{(x,y,z)|\re(x,y,z)\in D_{H,\epsilon},\im(x,y,z)=\im(\epsilon_1 \overline{x_0},-\overline{y_2},\epsilon_2\overline{z_0})\}$$
    $$S_{\text{side}}^{\epsilon_1 x,\epsilon_2 z}=\{(x,y+\frac{\log5}{4}\sqrt{-1},z)+t\sqrt{-1}\im(\epsilon_1 \overline{x_0},-\overline{y_2}-\frac{\log5}{4}\sqrt{-1},\epsilon_2\overline{z_0})|(x,y,z)\in \partial D_{H,\epsilon},t\in[0,1]\}$$
    $$S_{\text{bottom}}=\{(x,y+\frac{\log5}{4}\sqrt{-1},z)|(x,y,z)\in D_\epsilon\backslash D_{H,\epsilon}\}$$
    where $(x_0,y_2,z_0)$ is mentioned in Propositon \ref{Vol}. We can deform $S_{\text{top}}(\frac{\log5}{4})$ to $S^{\epsilon_1 x,\epsilon_2 z}=S_{\text{top}}^{\epsilon_1 x,\epsilon_2 z}\cup S_{\text{side}}^{\epsilon_1 x,\epsilon_2 z}\cup S_{\text{bottom}}$.

    By Lemma \ref{imV>3.5}, we have the following corollary.
    \begin{cor}\label{corbot}
        If $(x,y,z)\in S_{\text{bottom}}$, then
        $$\im V(x,y,z,s,m,0,l)<3.5$$
    \end{cor}



    Define $x_0^{\epsilon_1}=\epsilon_1\overline{x_0}$, $z_0^{\epsilon_2}=\begin{cases}
            \overline{z_0}, & \mbox{if }\epsilon_2=+\\
            \pi-\overline{z_0}, & \mbox{if }\epsilon_2=-
        \end{cases}$.

    \begin{lem}\label{leminD}
        For $\epsilon$ sufficiently small and $|p|,|q|,|p'|$ sufficiently large,
        $$\re(x_0^{\epsilon_1},-\overline{y_2},z_0^{\epsilon_2})\in D_{H,\epsilon}$$
    \end{lem}

    Lemma \ref{leminD} is because of
    $$\lim_{(p.q.p')\to \infty}\re(x_0^{\epsilon_1},-\overline{y_2},z_0^{\epsilon_2})=(0,\frac{\arctan2}{2},\frac{\pi}{2})$$

    Combine Proposition \ref{Vol}, Lemma \ref{H<0} and Lemma \ref{leminD}, we have
    \begin{pr}\label{prinDH}
        For $\epsilon$ sufficiently small and $|p|+|q|,|p'|$ sufficiently large, if $(x,y,z)\in S_{\text{top}}^{\epsilon_1 x,\epsilon_2 z}$, then
        $$\im(V(x,y,z,s^{\epsilon_1},m^{\epsilon_1},0,l^{\epsilon_2}))$$
        takes the unique maximum value $\vol(\mathcal{K}_{p'}(p,q))$ at the unique critical point $(x_0^{\epsilon_1},-\overline{y_2},z_0^{\epsilon_2})$.
    \end{pr}

    From now on we assume $\epsilon$ is sufficiently small and $|p|+|q|,|p'|$ are sufficiently large.

    By Lemma \ref{IJK}, we have $k(s,m)\in\frac{2\mathbb{Z}+1}{q}$, $(s,m)=(s^+,m^+)$ is the only pair such that $k(s,m)=\frac{1}{q}$, and $(s,m)=(s^-,m^-)$ is the only pair such that $k(s,m)=-\frac{1}{q}$. So for each pair $(s,m)$ other than $(s^\pm,m^\pm)$, we have $\mathrm{sgn}k(s,m)=\mathrm{sgn}(k(s,m)-k(s^\pm,m^\pm))$.

    For $l>l_+$ or $l<l_-$, we have $\mathrm{sgn}(l+p')=\mathrm{sgn}(l-l^\pm)$.

    \begin{pr}\label{pr<Vol}
        Suppose $l\ne-p'-1$ and $(s,m,l)$ other than $(s^{\epsilon_1},m^{\epsilon_1},l^{\epsilon_2})$. For $(x,y,z)\in S_{\text{top}}^{\epsilon_1 x,\epsilon_2z}$, we have
        $$\im V(x,y,z,s,m,0,l)\le\vol(\mathcal{K}_{p'}(p,q))-\epsilon_0$$
        where $\epsilon_0=\min\{|\frac{4\pi\im x_0}{q}|,|4\pi\im z_0|\}$.
    \end{pr}

    \begin{rem}
        By Proposition \ref{im!=0}, we have $\epsilon_0>0$.
    \end{rem}

    \begin{rem}
        We can give a geometric explanation for Proposition \ref{pr<Vol}. Similarly as Section \ref{sec5},we have $\im V(x,y,z,s,m,0,l)=\vol(\mathcal{W}(\frac{(p,q)}{qk(s,m)},\frac{(1,-p')}{p'+l+1}))$. Under the condition of Proposition \ref{pr<Vol}, we have $|qk(s,m)|\ge1, |p'+l+1|\ge1$ and one of it is strictly greater than $1$. So the filling coefficients $\frac{(p,q)}{qk(s,m)},\frac{(1,-p')}{p'+l+1}$ is strictly "less" than $(p,q),(1,-p')$. By the hyperbolic Dehn filling theory we have $\vol_{\vec{u}}$ is an "increasing" function of $\vec{u}$, so we believe $\vol(\mathcal{W}(\frac{(p,q)}{qk(s,m)},\frac{(1,-p')}{p'+l+1}))<\vol(\mathcal{W}((p,q),(1,-p'))$.
    \end{rem}

    \begin{proof}[Proof of Proposition \ref{pr<Vol}]
        Take
        $$\epsilon_1=\epsilon_1(k)=\mathrm{sgn}(k(s,m)\im x_0),\epsilon_2=\epsilon_2(l)= -\mathrm{sgn}((l+p')\im z_0)$$
        We have
\begin{eqnarray*}
  &&\im V(x,y,z,s,m,0,l)\\
  &=&\im V(x,y,z,s^{\epsilon_1},m^{\epsilon_1},0,l^{\epsilon_2})+2\pi(k(s,m)-k(s^{\epsilon_1},m^{\epsilon_1}))\im(\epsilon_1 \overline{x_0})\\
  &&-4\pi(l-l^{\epsilon_2})\im(\epsilon_2 \overline{z_0})\\
  &\le&\im V(x,y,z,s^{\epsilon_1},m^{\epsilon_1},0,l^{\epsilon_2})-\epsilon_0\\
  &\le& \im V(x_0^{\epsilon_1},-\overline{y_2},z_0^{\epsilon_2},s^{\epsilon_1},m^{\epsilon_1},0,l^{\epsilon_2})-\epsilon_0\\
  &=&\vol(\mathcal{K}_{p'}(p,q))-\epsilon_0
\end{eqnarray*}
    \end{proof}


    Combine Proposition \ref{pr<Vol}, Lemma \ref{imV>3.5} and Lemma \ref{lemyconvex}, we have for $(x,y,z)\in S_{\text{side}}^{\epsilon_1 x,\epsilon_2 z}$,
    \begin{eqnarray*}
      &&\im V(x,y,z,s,m,0,l)\\
      &\le&\max\{\im V(\re x+\sqrt{-1}\im(\epsilon_1 \overline{x_0}),\im \re y+\sqrt{-1}\im(-\overline{y_2}),\re z+\sqrt{-1}\im(\epsilon_1 \overline{z_0}),s,m,0,l),\\
      &&\im V(\re x,\re y+\frac{\log 5}{4}\sqrt{-1},\re z,s,m,0,l)\}\\
      &\le&\max\{\vol(\mathcal{K}_{p'}(p,q))-\epsilon_0,3.5\}
    \end{eqnarray*}
    Because of
    $$\lim_{(p.q.p')\to \infty}\vol(\mathcal{K}_{p'}(p,q))=\vol(\mathcal{W})=3.66\cdots>3.5$$
    we can choose $\delta=\min\{\vol(\mathcal{K}_{p'}(p,q))-3.5,\epsilon_0\}$, then we have

    \begin{cor}\label{corside2}
        for $(x,y,z)\in S_{\text{side}}^{\epsilon_1 x,\epsilon_2 z}$,
        $$\im V(x,y,z,s,m,0,l)\le\vol(\mathcal{K}_{p'}(p,q))-\delta$$
    \end{cor}


    Define $D_{(m,0,l)}=S_{\text{side}}(Y)\cup S^{\epsilon_1(k(s,m)) x,\epsilon_2(l) z}$.

    Combine Corollary \ref{corside1}, Corollary \ref{corbot}, Proposition \ref{pr<Vol}, Corollary \ref{corside2}, we conclude the following proposition.

\begin{pr}\label{n=0}
    Suppose $l\ne-p'-1$, $(s,m,l)$ is other than $(s^{\epsilon_1},m^{\epsilon_1},l^{\epsilon_2})$ where $\epsilon_1,\epsilon_2\in\{+,-\}$. There exists $\delta>0$ which is only depend on the geometry of $\mathcal{K}_{p'}(p,q))$, such that for $(x,y,z)\in D_{(m,0,l)}$
    $$\im V(x,y,z,s,m,n,l)<\vol(\mathcal{K}_{p'}(p,q))-\delta$$
\end{pr}

\subsection{Estimates of the leading Fourier coefficients}\label{subsecleading}
    In this section we will estimates the last four terms $\check{f_r}(s^{\epsilon_1},m^{\epsilon_1},0,l^{\epsilon_2})$, which is the leading terms. Suppose $p,q,p'$ is sufficient large.

    By Proposition \ref{pr2.1}, Proposition \ref{prD1}, Corollary \ref{corside1} and Corollary \ref{corbot}, we have
    \begin{eqnarray}
        \check{f}_r(s^{\epsilon_1},m^{\epsilon_1},0,l^{\epsilon_2})&=&(-1)^{m^{\epsilon_1}+l^{\epsilon_2}}\frac{e^{-\frac{3\pi\sqrt{-1}}{4}}r^{\frac{5}{2}}}{\sqrt{2}\pi^3}\int_{S_{\mathrm{top}}^{\epsilon_1 x,\epsilon_2z}\cup S_{\mathrm{side}}^{\epsilon_1 x,\epsilon_2z}}g(s^{\epsilon_1},x,y,z)e^{\frac{r}{4 \pi \sqrt{-1}}V(x,y,z,s,m,n,l)}\D x\D y\D z\nonumber\\
        &&\cdot(1+O(\frac{1}{r}))+o(e^{\frac{r}{4\pi}\cdot3.5})\nonumber
    \end{eqnarray}
    where
    $$g(s,x,y,z)=\frac{\sin(\frac{x}{q}-J(s)\pi)\sin (2z)}{\sqrt{(1-e^{2\sqrt{-1}(y+x)})(1-e^{2\sqrt{-1}(y-x)})(1-e^{-2\sqrt{-1}y})}}$$

    By Lemma \ref{H<0}, Lemma \ref{leminD} and Propositon \ref{prinDH}, there exists $\delta>0$ such that if $(x,y,z)\in\partial S_{\mathrm{top}}^{\epsilon_1 x,\epsilon_2z}$, we have
    $$\im V(x,y,z,s^{\epsilon_1},m^{\epsilon_1},0,l^{\epsilon_2})<\vol(\mathcal{K}_{p'}(p,q))-\delta$$

    Combine Lemma \ref{imV>3.5} and Lemma \ref{lemyconvex}, we have for $(x,y,z)\in S_{\text{side}}^{\epsilon_1 x,\epsilon_2 z}$,
    \begin{eqnarray*}
      &&\im V(x,y,z,s^{\epsilon_1},m^{\epsilon_1},0,l^{\epsilon_2})\\
      &\le&\max\{\im V(\re x+\sqrt{-1}\im(\epsilon_1 \overline{x_0}),\im \re y+\sqrt{-1}\im(-\overline{y_2}),\re z+\sqrt{-1}\im(\epsilon_1 \overline{z_0}),s^{\epsilon_1},m^{\epsilon_1},0,l^{\epsilon_2}),\\
      &&\im V(\re x,\re y+\frac{\log 5}{4}\sqrt{-1},\re z,s^{\epsilon_1},m^{\epsilon_1},0,l^{\epsilon_2})\}\\
      &\le&\max\{\vol(\mathcal{K}_{p'}(p,q))-\delta,3.5\}
    \end{eqnarray*}
    we choose $\delta<\vol(\mathcal{K}_{p'}(p,q))-3.5$, then we have for $(x,y,z)\in S_{\text{side}}^{\epsilon_1 x,\epsilon_2 z}$,
    $$\im V(x,y,z,s^{\epsilon_1},m^{\epsilon_1},0,l^{\epsilon_2})<\vol(\mathcal{K}_{p'}(p,q))-\delta$$

    \begin{cor}\label{corside3}
        For $|p|,|q|,|p'|$ sufficiently large, there exists $\delta>0$ such that
        $$|\int_{S_{\mathrm{side}}^{\epsilon_1 x,\epsilon_2z}}g(s,x,y,z)e^{\frac{r}{4 \pi \sqrt{-1}}V(x,y,z,s^{\epsilon_1},m^{\epsilon_1},0,l^{\epsilon_2})}\D x\D y\D z|<O(e^{\frac{r}{4\pi}(\vol(\mathcal{K}_{p'}(p,q))-\delta)})$$
    \end{cor}

    To estimate the integral in $S_{\mathrm{top}}^{\epsilon_1 x,\epsilon_2z}$, we need the saddle point method which was developed by Ohtsuki \cite{52}, Ka Ho Wong and Tian Yang \cite{Wong-Yang41}.

\begin{pr}[Saddle Point Method]\cite[Proposition 6.1]{Wong-Yang41}\label{saddle}
Let $D$ be a region in $\mathbb C^n$ and let $f(z_1,\dots, z_n)$ and $g(z_1,\dots, z_n)$ be holomorphic functions on $D$ independent of $r$. Let $f_r(z_1,\dots,z_n)$ be a holomorphic function of the form
$$ f_r(z_1,\dots, z_n) = f(z_1,\dots, z_n) + \frac{\upsilon_r(z_1,\dots,z_n)}{r^2}.$$
Let $S$ be an embedded real $n$-dimensional closed disk in $D$ and let $(c_1,\dots, c_n)$ be a point on $S.$ If
\begin{enumerate}[(1)]
\item $(c_1, \dots, c_n)$ is a critical point of $f$ in $D,$
\item $\mathrm{Re}(f)(c_1,\dots,c_n) > \mathrm{Re}(f)(z_1,\dots,z_n)$ for all $(z_1,\dots,z_n) \in S\setminus \{(c_1,\dots,c_n)\},$
\item the domain $\{(z_1,\dots,z_n)\in D\ |\ \mathrm{Re}{f}(z_1,\dots,z_n)<\mathrm{Re}{f}(c_1,\dots,c_n)\}$ deformation retracts to $S\setminus\{(c_1,\dots,c_n)\},$
\item the Hessian matrix $\mathrm{Hess}(f)(c_1,\dots,c_n)$ of $f$ at $(c_1,\dots,c_n)$ is non-singular,
\item $g(c_1,\dots,c_n) \neq 0,$ and
\item $|\upsilon_r(z_1,\dots,z_n)|$ is bounded from above by a constant independent of $r$ in $D,$
\end{enumerate}
then
\begin{equation*}
\begin{split}
 \int_S g(z_1,\dots, z_n) &e^{rf_r(z_1,\dots,z_n)} dz_1\dots dz_n\\
 &= \Big(\frac{2\pi}{r}\Big)^{\frac{n}{2}}\frac{g(c_1,\dots,c_n)}{\sqrt{\det(-\mathrm{Hess}(f)(c_1,\dots,c_n))}} e^{rf(c_1,\dots,c_n)} \Big( 1 + O \Big( \frac{1}{r} \Big) \Big).
 \end{split}
 \end{equation*}
\end{pr}

Let's verify the conditions (1),$\cdots$,(6) of Proposition \ref{saddle}.
\begin{enumerate}
  \item Take the $D_S=\{(x,y,z)\in\mathbb{C}^3|\re(x,y,z)\in D_{H,\epsilon},\im(x,y,z)\in[-L,L]^3\}$\footnote{$L$ is sufficiently large.} as the region $D$ in Proposition \ref{saddle}. (1) is because of Lemma \ref{leminD};
  \item (2) is because of Proposition \ref{prinDH};
  \item For each $(X,Y,Z)\in D_{H,\epsilon}$, let
  $$P_{(X,Y,Z)}=\{(x,y,z)\in D_S|\re(x,y,z)=(X,Y,Z)\}$$
  $$D_{(X,Y,Z)}^{\epsilon_1x,\epsilon_2z}=\{(x,y,z)\in P_{(X,Y,Z)}|\im V(x,y,z,s^{\epsilon_1},m^{\epsilon_1},0,l^{\epsilon_2})<\im V(x_0^{\epsilon_1},-\overline{y_2},z_0^{\epsilon_2},s^{\epsilon_1},m^{\epsilon_1},0,l^{\epsilon_2})\}$$
  Then we claim that $D_{\re(x_0^{\epsilon_1},-\overline{y_2},z_0^{\epsilon_2})}^{\epsilon_1x,\epsilon_2z}=\emptyset$ and $D_{(X,Y,Z)}^{\epsilon_1x,\epsilon_2z}$ is homeomorphic to a disk for $(X,Y,Z)\ne(x_0^{\epsilon_1},-\overline{y_2},z_0^{\epsilon_2})$, from which we conclude that the domains
  $$D_{(X,Y,Z)}^{\epsilon_1x,\epsilon_2z}=\{(x,y,z)\in D_S|\im V(x,y,z,s^{\epsilon_1},m^{\epsilon_1},0,l^{\epsilon_2})<\im V(x_0^{\epsilon_1},-\overline{y_2},z_0^{\epsilon_2},s^{\epsilon_1},m^{\epsilon_1},0,l^{\epsilon_2})\}$$
  respectively deformation retract to $S_{\mathrm{top}}^{\epsilon_1 x,\epsilon_2z}\backslash(x_0^{\epsilon_1},-\overline{y_2},z_0^{\epsilon_2})$ by shrinking each $D_{(X,Y,Z)}^{\epsilon_1x,\epsilon_2z}$ respectively to $\{(X,Y,Z)+\sqrt{-1}\im(x_0^{\epsilon_1},-\overline{y_2},z_0^{\epsilon_2})\}$, verifying (3) of Proposition \ref{saddle}.

  By Lemma \ref{H<0}, For $(x,y,z)\in P_{(x_0^{\epsilon_1},-\overline{y_2},z_0^{\epsilon_2})}$, $\im V$ respectively achieve the absolute minimum at $(x_0^{\epsilon_1},-\overline{y_2},z_0^{\epsilon_2})$, hence $D_{\re(x_0^{\epsilon_1},-\overline{y_2},z_0^{\epsilon_2})}^{\epsilon_1x,\epsilon_2z}=\emptyset$.

  For $(X,Y,Z)\ne\re(x_0^{\epsilon_1},-\overline{y_2},z_0^{\epsilon_2})$, we have
  \begin{eqnarray*}
    &&\min_{(x,y,z)\in P_{(X,Y,Z)}}\im V(x,y,z,s^{\epsilon_1},m^{\epsilon_1},0,l^{\epsilon_2})\\
    &\le&\im V(X+\sqrt{-1}\im(\epsilon_1 \overline{x_0}),Y+\sqrt{-1}\im(-\overline{y_2}),Z+\sqrt{-1}\im(\epsilon_1 \overline{z_0}),s^{\epsilon_1},m^{\epsilon_1},0,l^{\epsilon_2})\\
    &<&\im V(x_0^{\epsilon_1},-\overline{y_2},z_0^{\epsilon_2},s^{\epsilon_1},m^{\epsilon_1},0,l^{\epsilon_2})
  \end{eqnarray*}
  then by Lemma \ref{H<0} $\im V$ is concave up on $P_{(X,Y,Z)}$, so we have $D_{(X,Y,Z)}^{\epsilon_1x,\epsilon_2z}$ are nonempty convex subset of $P_{(X,Y,Z)}$, which is homeomorphic to a disk.

  \item (4) is because of Lemma \ref{H<0};
  \item
        \begin{eqnarray*}
          &&g(s^{\epsilon_1},x_0^{\epsilon_1},-\overline{y_2},z_0^{\epsilon_2})\\
          &=&\frac{\sin(\frac{x}{q}-J(s)\pi)\sin (2z)}{\sqrt{(1-e^{2\sqrt{-1}(y+x)})(1-e^{2\sqrt{-1}(y-x)})(1-e^{-2\sqrt{-1}y})}}|_{(s,x,y,z)=(s^{\epsilon_1},x_0^{\epsilon_1},-\overline{y_2},z_0^{\epsilon_2})}\\
          &\stackrel{(\ref{e24})}{=}&\frac{\sin(\frac{x}{q}-J(s)\pi)\sin (2z)}{\sqrt{(1-e^{2\sqrt{-1}(y+z)})(1-e^{2\sqrt{-1}(y-z)})}}|_{(s,x,y,z)=(s^{\epsilon_1},x_0^{\epsilon_1},-\overline{y_2},z_0^{\epsilon_2})}
        \end{eqnarray*}

        Because of $\re(x_0^{\epsilon_1},-\overline{y_2},z_0^{\epsilon_2})\in D_H$, the denominator $\sqrt{(\cdots)(\cdots)}\ne0$. By Proposition \ref{im!=0}, we have the numerator $\sin(\cdots)\sin(\cdots)\ne0$. So (5) holds.

  \item In our case $v_r=0$, (6) holds.
\end{enumerate}

    By Proposition \ref{saddle}, we have
    \begin{eqnarray*}
      &&\int_{S_{\mathrm{top}}^{\epsilon_1 x,\epsilon_2z}}g(s,x,y,z)e^{\frac{r}{4 \pi \sqrt{-1}}V(x,y,z,s^{\epsilon_1},m^{\epsilon_1},0,l^{\epsilon_2})}\D x\D y\D z\\
      &=&(\frac{2\pi}{r})^{\frac{3}{2}}t^{\epsilon_1,\epsilon_2}(M)e^{\frac{r}{4 \pi \sqrt{-1}}V(x_0^{\epsilon_1},-\overline{y_2},z_0^{\epsilon_2},s^{\epsilon_1},m^{\epsilon_1},0,l^{\epsilon_2})}(1+O(\frac{1}{r}))\\
      &=&(\frac{2\pi}{r})^{\frac{3}{2}}t^{\epsilon_1,\epsilon_2}(M)e^{\frac{r}{4 \pi }\cv(\mathcal{W}((p,q),(1,-p')))}(1+O(\frac{1}{r}))
    \end{eqnarray*}
    where
    \begin{equation}\label{et}
      t^{\epsilon_1,\epsilon_2}(M)=\frac{g(s^{\epsilon_1},x_0^{\epsilon_1},-\overline{y_2},z_0^{\epsilon_2})}{\sqrt{\mathrm{det}(-\hess(V)(x_0^{\epsilon_1},-\overline{y_2},z_0^{\epsilon_2}))}}
    \end{equation}


\begin{pr}\label{prleading}
    For $p,q,p'$ sufficient large,
    $$\check{f}_r(s^{\epsilon_1},m^{\epsilon_1},0,l^{\epsilon_2})=\frac{2re^{-\frac{3\pi\sqrt{-1}}{4}}}{\pi^\frac{3}{2}}t_0(M)e^{\frac{r}{4 \pi }\cv(\mathcal{W}((p,q),(1,-p')))}(1+O(\frac{1}{r}))$$
\end{pr}

    where $t_0(M)=(-1)^{m^{\epsilon_1}+l^{\epsilon_2}}t^{\epsilon_1,\epsilon_2}(M)$ has nothing to do with the choice of $\epsilon_1,\epsilon_2$ because of the symmetry of the Fourier terms (see Propositon \ref{fx=f-x} and Proposition \ref{cancell})\footnote{This can also be directly verified as in \cite{Wong-Yang41}.}.

\subsection{Proof of Theorem \ref{thm1} and Theorem \ref{thmtv}}\label{subsecpf}
Let express Theorem \ref{thm1} precisely.

\begin{thm}\label{thm1.0}
    Let $M=\mathcal{K}_{p'}(p,q)$ where $(p,q,p')\in\mathbb{Z}^3$, for $|p|,|q|,|r|$ sufficiently large,
    $$\rt_r(M)=C(r)t(M)e^{\frac{r}{4\pi}\cv(M)}(1+O(\frac{1}{r}))$$
    where
    $$C(r)=e^{\pi\sqrt{-1}(\frac{1}{r}(3\sigma(L)-\sum_{i=1}^k\limits a_i-\sum_{i=2}^k\limits \frac{1}{C_{i-1}C_i}-p'-\frac{5}{2})+\frac{3k}{4}+\sum_{i=1}^k\limits a_i+\frac{1}{2}+\frac{r}{4}(\sigma(L)+3a_k+2))}$$
    $$t(M)=\frac{8(-1)^{m^+-p'}\sin(\frac{\overline{x_0}}{q}-J(s^+)\pi)\sin (2\overline{z_0})}{\sqrt{q\pi^3(1-e^{2\sqrt{-1}(-\overline{y_2}+\overline{z_0})})(1-e^{2\sqrt{-1}(-\overline{y_2}-\overline{z_0})})\mathrm{det}(-\hess(V)(\overline{x_0},-\overline{y_2},\overline{z_0}))}}$$
    $L$ is defined as Figure \ref{figL}, $\sigma(L)$ is the signature of linking matrix of the link $L$. $a_i,C_i,s^+,m^+,J(s)$ are defined in Section \ref{secfrac}. $V(x,y,z,s,m,n,l)$ is the potential function defined in Definition \ref{defV}, $(x_0,y_2,z_0)$ is the critical point defined in (\ref{ex0}).
\end{thm}

\begin{proof}[Proof of Theorem \ref{thm1} (or Theorem \ref{thm1.0})]
    By Proposition \ref{pr2} and Proposition \ref{prleading}, it suffices to prove
    \begin{equation}\label{eother}
      \sum_{(s,m,n,l)\ne(s^{\epsilon_1},m^{\epsilon_1},0,l^{\epsilon_2})}|\check{f}_r(s,m,n,l)|<O(e^{\frac{r}{4\pi}(\vol(M)-\delta}))\nonumber
    \end{equation}
    for some $\delta>0$. Define
    $$\mathcal{S}=\{(s,m,n,l)\in\{0,\cdots,|q|-1\}\times\mathbb{Z}^3|(s,m,n,l)\ne(s^{\epsilon_1},m^{\epsilon_1},0,l^{\epsilon_2}),l\ne -p'-1\text{, and }(m,n,l)\ne(0,0,0)\}$$

    By Proposition \ref{cancell}, it suffices to prove
    $$\sum_{(s,m,n,l)\in\mathcal{S}}|\check{f}_r(s,m,n,l)|<o(e^{\frac{r}{4\pi}(\vol(M)-\delta}))$$

    Define
    $$h_r(s,x,y,z)=g(s,x,y,z)e^{\frac{v_r(x,y,z)}{4\pi\sqrt{-1}r}}$$
    where $v_r$ is defined in Lemma \ref{convergeV}. Then $h_r(s,x,y,z)$ is a smooth function for $\re(x,y,z)\in D_{\frac{\epsilon}{2}}$ and an analytic function for $\re(x,y,z)\in D_{\epsilon}$. Moreover, $h_r(s,x,y,z)=g(s,x,y,z)(1+O(\frac{1}{r}))$.

    By Lemma \ref{convergeV} and Proposition \ref{pr2.1}, we have
    $$|\check{f}_r(s,m,n,l)|=\frac{r^{\frac{5}{2}}}{\sqrt{2}\pi^3}|\int_{D_{\frac{\epsilon}{2}}}h(s,x,y,z)e^{\frac{r}{4 \pi \sqrt{-1}}V(x,y,z,s,m,n,l)}\D x\D y\D z|$$

    we have
    \begin{eqnarray*}
      &&r^6(m^6+n^6+l^6)\int_{D_{\frac{\epsilon}{2}}}h(s,x,y,z)e^{\frac{r}{4 \pi \sqrt{-1}}V(x,y,z,s,m,n,l)}\D x\D y\D z\\
      &=&-\int_{D_{\frac{\epsilon}{2}}}h(s,x,y,z)e^{\frac{r}{4 \pi \sqrt{-1}}V(x,y,z,s,0,0,0)}((\frac{\partial^6}{\partial x^6}+\frac{\partial^6}{\partial y^6}+\frac{\partial^6}{\partial z^6})e^{\frac{r}{4 \pi \sqrt{-1}}(-4\pi mx-4\pi ny-4\pi lz)})\D x\D y\D z\\
      &=&-\int_{D_{\frac{\epsilon}{2}}}((\frac{\partial^6}{\partial x^6}+\frac{\partial^6}{\partial y^6}+\frac{\partial^6}{\partial z^6})(h(s,x,y,z)e^{\frac{r}{4 \pi \sqrt{-1}}V(x,y,z,s,0,0,0)}))e^{\frac{r}{4 \pi \sqrt{-1}}(-4\pi mx-4\pi ny-4\pi lz)}\D x\D y\D z\\
      &=&-\int_{D_{\frac{\epsilon}{2}}}\tilde{h}_r(s,x,y,z)e^{\frac{r}{4 \pi \sqrt{-1}}V(x,y,z,s,m,n,l)}\D x\D y\D z
    \end{eqnarray*}
    where
    $$\tilde{h}_r(s,x,y,z)=\frac{(\frac{\partial^6}{\partial x^6}+\frac{\partial^6}{\partial y^6}+\frac{\partial^6}{\partial z^6})(h(s,x,y,z)e^{\frac{r}{4 \pi \sqrt{-1}}V(x,y,z,s,0,0,0)})}{e^{\frac{r}{4 \pi \sqrt{-1}}V(x,y,z,s,0,0,0)}}$$
    is a smooth function independent of $(m,n,l)$, and has the form
    $$\tilde{h}_r(s,x,y,z)=\tilde{h}(s,x,y,z)r^6+O(r^5)$$
    for a smooth function $\tilde{h}(s,x,y,z)$ independent of $r$.

    On the compact set
    $$S=(D_{\frac{\epsilon}{2}}\backslash D_{\epsilon})\bigcup_{(s,m,n,l)\in\mathcal{S}}\limits D_{(m,n,l)}=(D_{\frac{\epsilon}{2}}\backslash D_{\epsilon})\cup S(\frac{\log 5}{4})\cup S(-\frac{\log 5}{4})\cup S_{\text{top}}^{\epsilon_1 x,\epsilon_2 z}\cup S_{\text{side}}^{\epsilon_1 x,\epsilon_2 z}$$
    $|\frac{\tilde{h}_r(s,x,y,z)}{r^6}|$ is bounded from above by some $C>0$ independent of $(s,m,n,l)$ and $r$. Then we have

    \begin{eqnarray*}
      &&\sum_{(s,m,n,l)\in\mathcal{S}}|\check{f}_r(s,m,n,l)|\\
      &=&\frac{r^{\frac{5}{2}}}{\sqrt{2}\pi^3}\sum_{(s,m,n,l)\in\mathcal{S}}|\int_{D_{\frac{\epsilon}{2}}}h(s,x,y,z)e^{\frac{r}{4 \pi \sqrt{-1}}V(x,y,z,s,m,n,l)}\D x\D y\D z|\\
      &=&\frac{r^{\frac{5}{2}}}{\sqrt{2}\pi^3}\sum_{(s,m,n,l)\in\mathcal{S}}\frac{1}{m^6+n^6+l^6}|\int_{D_{\frac{\epsilon}{2}}}\frac{\tilde{h}_r(s,x,y,z)}{r^6}e^{\frac{r}{4 \pi \sqrt{-1}}V(x,y,z,s,m,n,l)}\D x\D y\D z|\\
      &=&\frac{r^{\frac{5}{2}}}{\sqrt{2}\pi^3}\sum_{(s,m,n,l)\in\mathcal{S}}\frac{1}{m^6+n^6+l^6}|(\int_{D_{\epsilon}}\frac{\tilde{h}_r(s,x,y,z)}{r^6}e^{\frac{r}{4 \pi \sqrt{-1}}V(x,y,z,s,m,n,l)}\D x\D y\D z\\
      &&+\int_{D_{\frac{\epsilon}{2}}\backslash D_{\epsilon}}\frac{\tilde{h}_r(s,x,y,z)}{r^6}e^{\frac{r}{4 \pi \sqrt{-1}}V(x,y,z,s,m,n,l)}\D x\D y\D z)|
      \end{eqnarray*}
      \begin{eqnarray*}
      &\le&\frac{r^{\frac{5}{2}}}{\sqrt{2}\pi^3}\sum_{(s,m,n,l)\in\mathcal{S}}\frac{1}{m^6+n^6+l^6}(|\int_{D_{\epsilon}}\frac{\tilde{h}_r(s,x,y,z)}{r^6}e^{\frac{r}{4 \pi \sqrt{-1}}V(x,y,z,s,m,n,l)}\D x\D y\D z|\\
      &&+\int_{D_{\frac{\epsilon}{2}}\backslash D_{\epsilon}}|\frac{\tilde{h}_r(s,x,y,z)}{r^6}||e^{\frac{r}{4 \pi \sqrt{-1}}V(x,y,z,s,m,n,l)}|\D x\D y\D z)\\
      &\stackrel{\text{Proposition }\ref{prD1}}{<}&\frac{r^{\frac{5}{2}}}{\sqrt{2}\pi^3}\sum_{(s,m,n,l)\in\mathcal{S}}\frac{1}{m^6+n^6+l^6}(|\int_{D_{(m,n,l)}}\frac{\tilde{h}_r(s,x,y,z)}{r^6}e^{\frac{r}{4 \pi \sqrt{-1}}V(x,y,z,s,m,n,l)}\D x\D y\D z|\\
      &&+C\vol(D_{\frac{\epsilon}{2}}\backslash D_{\epsilon})e^{\frac{r}{4\pi}\cdot3.5})\\
      &\le&\frac{r^{\frac{5}{2}}}{\sqrt{2}\pi^3}\sum_{(s,m,n,l)\in\mathcal{S}}\frac{1}{m^6+n^6+l^6}(\int_{D_{(m,n,l)}}|\frac{\tilde{h}_r(s,x,y,z)}{r^6}||e^{\frac{r}{4 \pi \sqrt{-1}}V(x,y,z,s,m,n,l)}|\D x\D y\D z\\
      &&+C\vol(D_{\frac{\epsilon}{2}}\backslash D_{\epsilon})e^{\frac{r}{4\pi}\cdot3.5})\\
      &<&\frac{r^{\frac{5}{2}}}{\sqrt{2}\pi^3}\sum_{(s,m,n,l)\in\mathcal{S}}\frac{1}{m^6+n^6+l^6}(C\vol(D_{(m,n,l)})e^{\frac{r}{4\pi}(\vol(M)-\delta)}+C\vol(D_{\frac{\epsilon}{2}}\backslash D_{\epsilon})e^{\frac{r}{4\pi}\cdot3.5})\\
      &<&\frac{r^{\frac{5}{2}}}{\sqrt{2}\pi^3}C(\sum_{(s,m,n,l)\in\mathcal{S}}\frac{1}{m^6+n^6+l^6})\vol(S)e^{\frac{r}{4\pi}(\vol(M)-\delta)}\\
      &=&o(e^{\frac{r}{4\pi}(\vol(M)-\delta_0)})
    \end{eqnarray*}
    where the inequality in  antepenultimate line is because of Proposition \ref{n!=0} and Proposition \ref{n=0}, $\delta$ is mentioned in Corollary \ref{corside2}, and $0<\delta_0<\delta$. We conclude the proof.
\end{proof}

Because of $\rt_r(M)$ and $\cv(M)$ are invariant of $M$, so $C(r)t(M)$ is an invariant of $M$ and $r$. By definition, we know $C(r)$ is determined by $r$ and $M$ with it's $(p,q)$ Dehn-Surgey structure,$|C(r)|=1$, $t(M)$ is independent of $r$. So we have the following corollary.

\begin{cor}\label{cortor}
    $|t(M)|$ is an invariant of $M$ and $t(M)$ is an invariant of $M$ with it's Dehn-Surgey structure.
\end{cor}

\begin{proof}[Proof of Theorem \ref{thmtv}]
    Ka Ho Wong and Tian Yang proved(see \cite{Wong-Yang41} Section 6 Proof of Theorem 1.2)
    $$\tv_r(M)=2^{b_2(M)-b_0(M)+2}|\rt_r(M)|^2$$
    where $b_0(M)=1$, $b_2(M)$ is the zeroth and second $\mathbb{Z}_2$-Betti numbers of $M$.  So Theorem \ref{thmtv} is a corollary of Theorem \ref{thm1}.
\end{proof}

\section{Further discussion}\label{sec7}

\subsection{Numerical calculation of $p,q$ and $p'$}
    By Lemma \ref{leminD}, Appendix Proof of Lemma \ref{imV>3.5} and the analysis at the beginning of Section \ref{sec5.4}, we have
    \begin{thm}
        If $(p,q,p')$ satisfy the following conditions,
        \begin{enumerate}[(1)]
          \item $\mathcal{K}_{p'}(p,q)$ is hyperbolic;
          \item $\re(x_0^{\epsilon_1},-\overline{y_2},z_0^{\epsilon_2})\in D_{H}$
          \item $\vol(\mathcal{K}_{p'}(p,q))>\max_{0<y<\frac{\pi}{4}}\limits\im v(0,y+\frac{\log5}{4}\sqrt{-1},\frac{\pi}{2}-y)=3.4595\cdots$
        \end{enumerate}
        then Theorem \ref{thm1} holds for $M=\mathcal{K}_{p'}(p,q)$. Where $v(x,y,z)$ is defined in Appendix Proof of Lemma \ref{imV>3.5}.
    \end{thm}

    Condition (1) is completely settled by Lemma \ref{nonhyperbolic2}. Through numerical calculation, we found that Condition (2) is easily satisfied, and for
    \begin{equation}\label{con18}
      |p|+|q|\ge18,|p'|\ge10
    \end{equation}
    Condition (3) are satisfied. So we believe Theorem \ref{thm1} holds for \ref{con18}.

\subsection{Volume Conjecture for twist knots at $e^{\frac{2\pi\sqrt{-1}}{N+\frac{1}{2}}}$}

    Chen-Zhu showed that

    \begin{thm}[\cite{Chen-Zhu-1}]
        $$\lim_{N\to\infty}\frac{2\pi}{N}\log J_N(\mathcal{K}_p;e^{\frac{2\pi\sqrt{-1}}{N+\frac{1}{2}}})\equiv \cv(S^3\backslash\mathcal{K}_p)\mod\pi^2\sqrt{-1}$$
    \end{thm}

    Our method  might provide a new proof of the above theorem.

\subsection{Conjecture about symmetry of the Fourier coefficients}\label{sec7.3}

    In Section \ref{subsec4.3} we discuss the symmetry of the Fourier coefficients for $M=\mathcal{K}_{p'}(p,q)$. Actually, the proof of Proposition \ref{fx=f-x} applies to the general case. We have the following theorem.

    \begin{thm}\label{thmsym}
        For a knot $K\subset S^3$ with the Reshetikhin-Turaev invariants in the following form(see Proposition \ref{pr1} and the formula (3.9) of Proposition 3.4 in \cite{Wong-Yang41}),
        $$\rt_r(K(p,q))=\kappa_r\sum_{s=0}^{|q|-1}\sum_{m=1-\frac{r}{2}}^{\frac{r}{2}-1}\sum_{n=|m|}^{\frac{r}{2}-1}\sin(\frac{x}{q}-J(s)\pi)\left(e^{-\sqrt{-1}x+\frac{r}{4\pi\sqrt{-1}}(-\frac{p}{q}x^2+\frac{2\pi I(s)}{q}x+4xy+K(s)\pi^2)}\right)$$
        $$\cdot\frac{(t)_{r-m-n-1}}{(t)_{n-m}}f_{K}(\frac{r}{2}-n-1)$$
        assume we can take the discrete Fourier transform for the above formula, and we get the Fourier term $\check{f}_r(s,m,n,\cdots)$. Then for $(s,m)$ and $(s',m')$ such that $k(s,m)+k(s',m')=0$, we have
        $$\check{f}_r(s,m,n,\cdots)=\check{f}_r(s',m',n,\cdots)$$
        where $k(s,m)$ is defined in (\ref{ek}).
    \end{thm}

    A more general statement of the symmetry of the Fourier coefficients is in Conjecture \ref{conprogram} (2) (7).

\subsection{Formula of $\rt(\mathcal{W}((p,q),(1,-p')))$}\label{sec7.4}

    We use the formula of $\rt(\mathcal{W}((p,q),(1,-p')))$ comes from Masbaum, Wong and Yang as in Section \ref{sec3}. It's turns out that this formula corresponding to the geometry of $\mathcal{W}((p+4q,-q),(1,p'-\frac{1}{2}))$ instead of $\rt(\mathcal{W}((p,q),(1,-p')))$ as in Section \ref{secV}. So we have the following question.

    \begin{qu}
      Is there exits an ideal triangulation $T$ of $\mathcal{W}$ and an expression of $\rt(\mathcal{W}((p,q),(1,-p')))$, such that this expression corresponding to the geometry of $\mathcal{W}((p,q),(1,-p'))$ with respect to $T$ directly?
    \end{qu}

\subsection{Idea of a proof of the Volume Conjecture}\label{sec7.5}

    Further, we have the following conjecture, which is our program to prove volume conjectures.
    \begin{con}\label{conprogram}
        For odd $r$ and a hyperbolic link $L\subset S^3$ with $l$ components $L_1,\cdots,L_l$, and sufficient large $|p_1|+|q_1|,\cdots,|p_l|+|q_l|$, there is a ideal triangulation $T$ of $L$, such that
        \begin{enumerate}[(1)]
          \item  The Reshetikhin-Turaev invariants of $M=L((p_1,q_1),\cdots,(p_l,q_l))$ can be expressed in the following form
              $$\rt_r(M)=\sum_{(x_1,\cdots,x_l, y_1,\cdots,y_k)\in \mathbb{Z}^{l+k}\cap \Omega}c_r(x_1,\cdots,x_l, y_1,\cdots,y_k)e^{\frac{r}{4\pi\sqrt{-1}}V_r(x_1,\cdots,x_l, y_1,\cdots,y_k)}$$
              This is the generalization of (\ref{e8}).

          \item There is a subregion $\Omega_0\subset \Omega\subset\mathbb{R}^{l+k}$ such that $\Omega_0$ is symmetry for the first $l$ coordinates and
              \begin{eqnarray*}
                &&\rt_r(M)=\sum_{(x_1,\cdots,x_l, y_1,\cdots,y_k)\in \mathbb{Z}^{l+k}\cap \Omega_0}a_r(x_1,\cdots,x_l, y_1,\cdots,y_k)e^{\frac{r}{4\pi\sqrt{-1}}V_r(x_1,\cdots,x_l, y_1,\cdots,y_k)}+o(e^{\frac{r}{4\pi\sqrt{-1}}(\vol(M)-\delta)})
              \end{eqnarray*}
              for some small $\delta>0$. This is the generalization of (\ref{e9}).

          \item We might take the discrete Fourier transform for the above formula, then we have
            $$\rt_r(M)=\sum_{(m_1,\cdots,m_l,n_1,\cdots,n_k)\in \mathbb{Z}^{l+k}}\check{f}(m_1,\cdots,m_l,n_1,\cdots,n_k)+o(e^{\frac{r}{4\pi\sqrt{-1}}(\vol(M)-\delta)})$$
            where
            \begin{eqnarray*}
              &&\check{f}(m_1,\cdots,m_l,n_1,\cdots,n_k)=\int_{\Omega_0}b_r(x_1,\cdots,x_l, y_1,\cdots,y_k,m_1,\cdots,m_l,n_1,\cdots,n_k)\\
              &&\cdot e^{\frac{r}{4\pi\sqrt{-1}}V_r(x_1,\cdots,x_l, y_1,\cdots,y_k,m_1,\cdots,m_l,n_1,\cdots,n_k)}\D\vec{\mathbf{x}}\D\vec{\mathbf{y}}
            \end{eqnarray*}
            $$V_r(x_1,\cdots,x_l, y_1,\cdots,y_k,m_1,\cdots,m_l,n_1,\cdots,n_k)=V_r(x_1,\cdots,x_l, y_1,\cdots,y_k)-4\pi(\sum_{i=1}^l m_ix_i+\sum_{i=1}^k n_iy_i)$$

            This is the generalization of Proposition \ref{pr2}.

          \item For $(x_1,\cdots,x_l, y_1,\cdots,y_k)\in \Omega_0$, we have the following uniform convergence of potential function.
            $$\lim_{r\to+\infty}V_r(x_1,\cdots,x_l, y_1,\cdots,y_k,m_1,\cdots,m_l,n_1,\cdots,n_k)=V(x_1,\cdots,x_l, y_1,\cdots,y_k,m_1,\cdots,m_l,n_1,\cdots,n_k)$$

            This is the generalization of Lemma \ref{convergeV}.

          \item There exist a sister potential function $W(x_1,\cdots,x_l, y_1,\cdots,y_k,m_1,\cdots,m_l,n_1,\cdots,n_k)$, such that $W-V$ is a $2$-degree polynomial in $x_1,\cdots,x_l$, and
              $$W(x_1,\cdots,x_l, y_1,\cdots,y_k,m_1,\cdots,m_l,n_1,\cdots,n_k)=\pm V(x_1,\cdots,x_l, y'_1,\cdots,y'_k,m'_1,\cdots,m'_l,n_1,\cdots,n_k)$$
              where $(y_1,\cdots, y_k)$ and $(y'_1,\cdots, y'_k)$ both are crtical points of the critical equation $\frac{\partial V}{\partial y_i}=\frac{\partial W}{\partial y_i}=0$ for $i=1,\cdots,k$, and $m_i-m'_i$ are some constants $\in\mathbb{Z}$.

              This is the generalization of Lemma \ref{V+W}.

          \item  There exists a hyperbolic manifold $M'$, and $d_1,\cdots,d_l$ which are $1$-degree polynomial in $m_1,\cdots,m_l$ whose coefficients $\in 2\pi\mathbb{Z}$, $\exists!(n_{1,0},\cdots,n_{k,0})\in\mathbb{Z}^k$, such that $\vol(M')=\vol(M)$, and
            \begin{enumerate}
              \item the critical equation $\frac{\partial W}{\partial x_i}=0$ corresponding to the Dehn-surgery equation of the components $L_i$ with respect to $T$:
                  $$p_i\mathbf{m}_i+q_i\mathbf{l}_i=d_i(m_1,\cdots,m_l)\sqrt{-1}$$
                  where $\mathbf{m}_i$ and $\mathbf{l}_i$ are the holonomy of the meridian and the longitude of $L_i$.
              \item the critical equation $\frac{\partial W}{\partial y_i}(x_1,\cdots,x_l, y_1,\cdots,y_k,m_1,\cdots,m_l,n_{1,0},\cdots,n_{k,0})=0$ corresponding to the edge gluing equations of $L$ with respect to $T$.
            \end{enumerate}

            This is the generalization of Lemma \ref{=}.

          \item If $(m_1,\cdots,m_l)$ and $(m'_1,\cdots,m'_l)$ satisfy $d_j(m_1,\cdots,m_l)=d_j(m'_1,\cdots,m'_l)$ for $j\ne i$ and $d_i(m_1,\cdots,m_l)=-d_i(m'_1,\cdots,m'_l)$, then we have
              $$\check{f}(m_1,\cdots,m_l,n_1,\cdots,n_k)=(-1)^{1+\frac{d_i(m_1,\cdots,m_l)}{2\pi}}\check{f}(m'_1,\cdots,m'_l,n_1,\cdots,n_k)$$

              This is the generalization of Theorem \ref{thmsym} and the \textbf{Big Cancellation}.

          \item $\exists!(m_{1,0},\cdots,m_{l,0})\in\mathbb{Z}^m$ such that $d_i(m_{1,0},\cdots,m_{l,0})=2\pi$ for $i=1,\cdots,l$, and $\exists!$ the critical point $(x_{1,0},\cdots,x_{l,0}, y_{1,0},\cdots,y_{k,0})\in\Omega_0$ of the critical equations
              \begin{eqnarray*}
                \frac{\partial W}{\partial x_i}(x_1,\cdots,x_l, y_1,\cdots,y_k,m_{1,0},\cdots,m_{l,0},n_{1,0},\cdots,n_{k,0})&=&0\\
                \frac{\partial W}{\partial y_i}(x_1,\cdots,x_l, y_1,\cdots,y_k,m_{1,0},\cdots,m_{l,0},n_{1,0},\cdots,n_{k,0})&=&0
              \end{eqnarray*}

              This might be a corollary of (6).

          \item $$W(x_{1,0},\cdots,x_{l,0}, y_{1,0},\cdots,y_{k,0},m_{1,0},\cdots,m_{l,0},n_{1,0},\cdots,n_{k,0})=\sqrt{-1}(\vol(M')\pm\sqrt{-1}\cs(M'))$$

              This is the generalization of Proposition \ref{W=vol}. This is a corollary of (6) because of Lemma \ref{lem2}.

              Then we have
              $$V(\vec{\mathbf{x}},\vec{\mathbf{y}},\vec{\mathbf{m}},n_{1,0},\cdots,n_{k,0})=\sqrt{-1}\cv(M)$$
              for $2^l$ certain tuples $(\vec{\mathbf{x}},\vec{\mathbf{y}},\vec{\mathbf{m}})$ which is determined by $(x_{1,0},\cdots,x_{l,0}, y_{1,0},\cdots,y_{k,0},m_{1,0},\cdots,m_{l,0})$.

              This is the generalization of Proposition \ref{Vol}.

          \item We might prove
          $$\im V(x_1,\cdots,x_l, y_1,\cdots,y_k,m_1,\cdots,m_l,n_1,\cdots,n_k)<\vol(M)-\delta$$
           for some small $\delta$ and $(n_1,\cdots,n_k)\ne(n_{1,0},\cdots,n_{k,0})$, by deform $\Omega_0$ alone the imaginary part of $y_i$ as in Section \ref{subsecy}.

          \item We might prove
          $$\im V(x_1,\cdots,x_l, y_1,\cdots,y_k,m_1,\cdots,m_l,n_{1,0},\cdots,n_{k,0})<\vol(M)-\delta$$
          for some small $\delta$ and the terms which don't corresponding to the $2^l$ terms mentioned in (9) and don't corresponding to the cancelled terms mentioned in (7), as we prove Proposition \ref{pr<Vol}.

          \item We might use saddle point method to prove that for $2^l$ certain tuples $\vec{\mathbf{m}}$ mentioned in (9),
              $$\check{f}(\vec{\mathbf{m}}, n_{1,0},\cdots,n_{k,0})=c(M)e^{\frac{r}{4\pi}\cv(M)}(1+O(\frac{1}{r}))$$
              as in Section \ref{subsecleading}.

          \item We might prove the sum of all terms which is not mentioned in (12) are neglectable, as in Section \ref{subsecpf}.

              Finally, we have
              $$\rt(M)=2^lc(M)e^{\frac{r}{4\pi}\cv(M)}(1+O(\frac{1}{r}))$$
        \end{enumerate}
    \end{con}

    In our proof, $x$ corresponding to two integer parameters $s,m$. We can use $qk(s,m)\in\mathbb{Z}$ to reduce the number of parameters to one, then our prove applies to the above conjecture. We can restate Theorem \ref{thm1} as below.
    \begin{thm}
        For $M=\mathcal{W}((p,q),(1,-p'))$, Conjecture \ref{conprogram} is true.
    \end{thm}

    Conjecture \ref{conprogram} may also apply to other quantum invariants, such as the colored Jones polynomial of hyperbolic knots.

\section{Appendix}\label{sec8}

We put the elementary proof in the appendix. Mathematica can draw nice graphs to help us understand the proof, and it can simplify complex elementary expressions very well.

\begin{proof}[Proof of Lemma \ref{lemim<3.5}]
    It suffices to prove that for $(x,y,z)\in \overline{D\backslash D_0}$,
    $$\im V(x,y,z,s,m,n,l)<3.5$$

    Define $$v(x,y,z)=\li(e^{2\sqrt{-1}(y+x)})+\li(e^{2\sqrt{-1}(y-x)})+\li(e^{2\sqrt{-1}(-y+z)})+\li(e^{2\sqrt{-1}(-y-z)})\\
      -\li(e^{-2\sqrt{-1}y})$$
    then for $(x,y,z)\in \overline{D\backslash D_0}$, $\im v(x,y,z)=\im V(x,y,z,s,m,n,l)$.

    \begin{enumerate}[(1)]
      \item First let's study the imaginary part of the function
      $$v_1(x,y)=\li(e^{2\sqrt{-1}(y+x)})+\li(e^{2\sqrt{-1}(y-x)})$$
       on $D_{xy}=\{(x,y)\in\mathbb{R}^2|y\in(|x|,\pi)\}$. We can divide $D_{xy}$ into four smaller regions (see Figure \ref{Dxy}).
       \begin{figure}[H]
         \centering
         \includegraphics[width=0.5\textwidth]{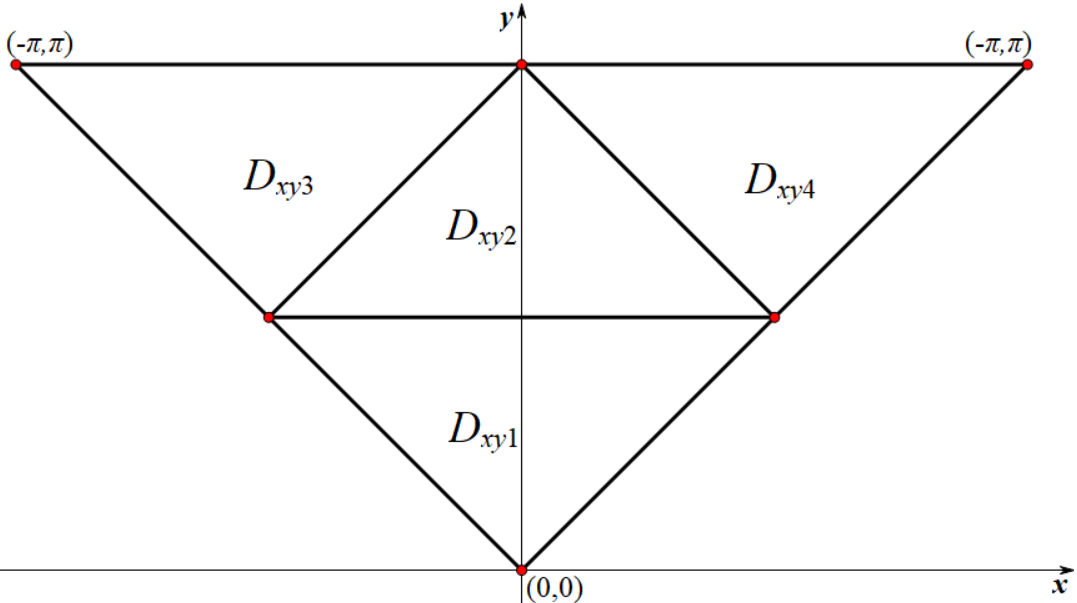}
         \label{Dxy}
       \end{figure}
       $$D_{xy}=D_{xy1}\cup D_{xy2} \cup D_{xy3}\cup D_{xy4}$$
       where
       \begin{eqnarray*}
         D_{xy1}&=&\{(x,y)\in\mathbb{R}^2|y\in[|x|,\frac{\pi}{2}]\}\\
         D_{xy2}&=&\{(x,y)\in\mathbb{R}^2|y\in[\frac{\pi}{2},\pi-|x|]\}\\
         D_{xy3}&=&\{(x,y)\in\mathbb{R}^2|y\in[|x+\frac{\pi}{2}|+\frac{\pi}{2},\pi)\}\\
         D_{xy4}&=&\{(x,y)\in\mathbb{R}^2|y\in[|x-\frac{\pi}{2}|+\frac{\pi}{2},\pi)\}
       \end{eqnarray*}

       \begin{figure}[H]
         \centering
         \includegraphics[width=0.8\textwidth]{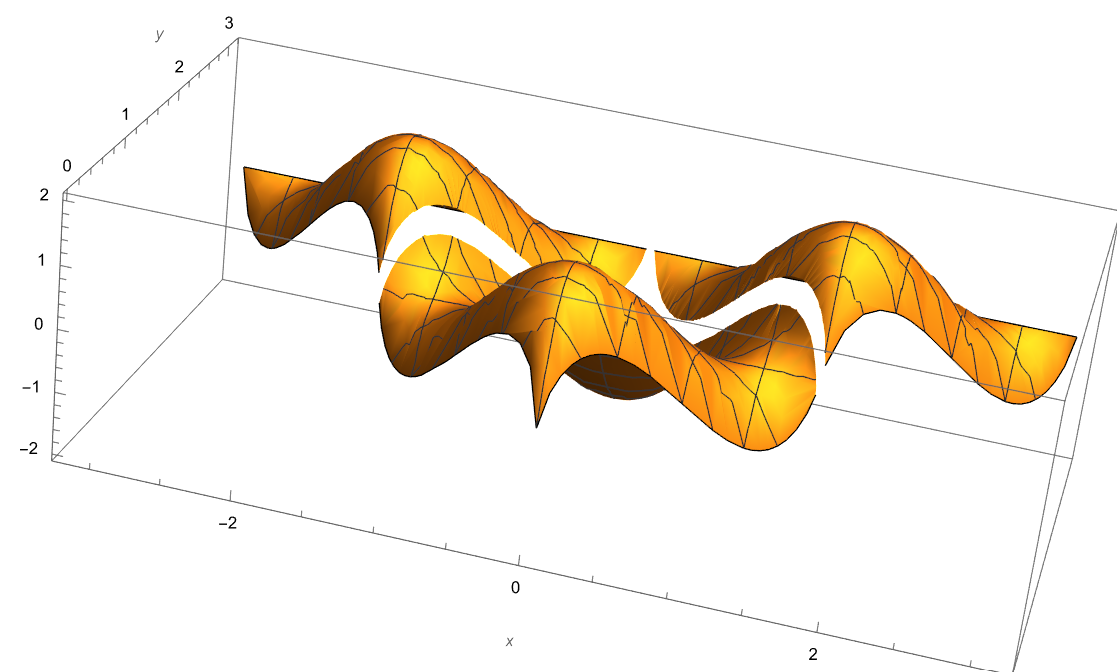}
         \caption{The graph of $\im v_1$ on $D_{xy}$. Due to the lack of precision in Mathematica, white ribbon gaps appear in the image. The graph is actually a continuous surface.}\label{v1}
       \end{figure}

       \begin{cl}
        \begin{enumerate}
          \item In $D_{xy1}$, $D_{xy3}$ and $D_{xy4}$, $\im v_1(x,y)$ is strictly concave down in $x$; in $D_{xy2}$, $\im v_1(x,y)$ is strictly concave up in $x$;
          \item For $(x,y)\in D_{x,y}$,
            \begin{enumerate}[1]
              \item If $y\in[0,\frac{\pi}{2})$, $\im v_1(x,y)\le \im v_1(0,y)$
              \item If $y=\frac{\pi}{2}$, $\im v_1(x,y)=0$
              \item If $y\in(\frac{\pi}{2},\pi)$, $\im v_1(x,y)\le \im v_1(\frac{\pi}{2},y)=\im v_1(-\frac{\pi}{2},y)$
            \end{enumerate}
        \end{enumerate}
       \end{cl}


       (a) is because of
       $$\frac{\D\im v_1(x,y)}{\D x}=\im (v_1)'_x(x,y)=2\log|\sin(y-x)|-2\log|\sin(y+x)|$$
       $$\frac{\D^2\im v_1(x,y)}{\D x^2}=-2(\cot(y-x)+\cot(y+x))=-\frac{2\sin2y}{\sin(y-x)\sin(y+x)}$$

       we have $\frac{\D\im v_1(x,y)}{\D x}=0\Leftrightarrow x\in\frac{\pi\mathbb{Z}}{2}$ or $y\in\frac{\pi\mathbb{Z}}{2}$. Combine this with (a), we get (b).

      \item Let's study the imaginary part of the function
        $$v_2(x,y)=\li(e^{2\sqrt{-1}(-y+z)})+\li(e^{2\sqrt{-1}(-y-z)})$$
        on $D_{yz}=\{(y,z)\in\mathbb{R}^2|z\in(0,\pi-y)\}$. We can divide $D_{yz}$ into three smaller regions (see Figure \ref{Dyz})
        $$D_{yz}=D_{yz1}\cup D_{yz2} \cup D_{yz3}$$

       \begin{figure}[H]
         \centering
         \includegraphics[width=0.3\textwidth]{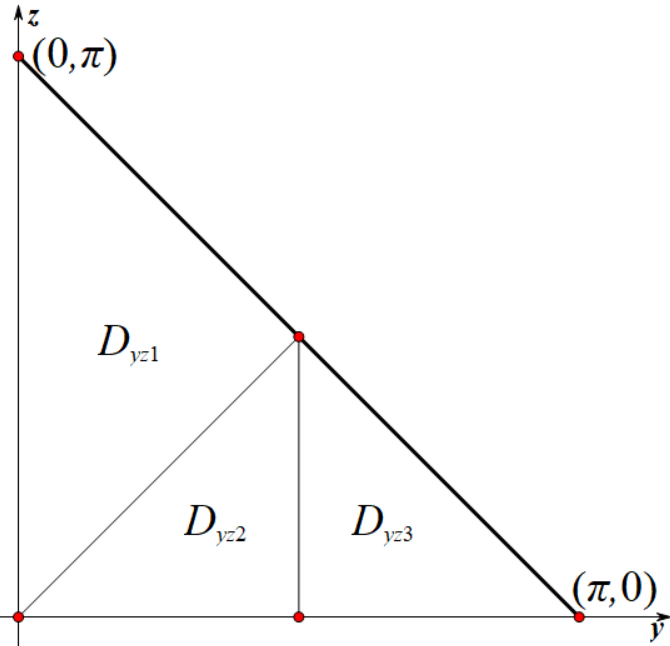}
         \label{Dyz}
       \end{figure}

       where
       \begin{eqnarray*}
         D_{yz1}&=&\{(y,z)\in\mathbb{R}^2|y\in(0,|z-\frac{\pi}{2}|+\frac{\pi}{2})\}\\
         D_{yz2}&=&\{(y,z)\in\mathbb{R}^2|y\le\frac{\pi}{2},0<z\le y)\}\\
         D_{yz3}&=&\{(y,z)\in\mathbb{R}^2|y\ge\frac{\pi}{2},0<z\le\pi-y)\}
       \end{eqnarray*}

        \begin{figure}[H]
         \centering
         \includegraphics[width=0.4\textwidth]{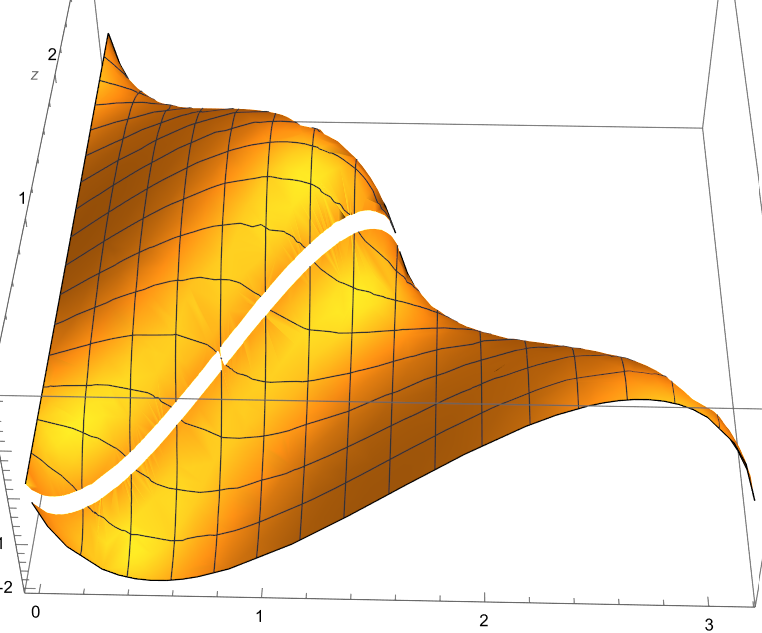}
         \caption{The graph of $\im v_2$ on $D_{yz}$. Due to the lack of precision in Mathematica, white ribbon gaps appear in the image. The graph is actually a continuous surface.}\label{v2}
       \end{figure}

       \begin{cl}\label{cl2}
        \begin{enumerate}
          \item In $D_{yz1}$ and $D_{yz3}$, $\im v_2(y,z)$ is strictly concave down in $z$; in $D_{yz2}$, $\im v_2(y,z)$ is an increasing function of $z$;
          \item For $(y,z)\in D_{y,z}$,
            \begin{enumerate}[1]
              \item If $y\in(0,\frac{\pi}{2})$, $\im v_2(y,z)\le \im v_2(y,\frac{\pi}{2})$
              \item If $y=\frac{\pi}{2}$, $\im v_1(x,y)=0$
              \item If $y\in(\frac{\pi}{2},\pi)$, $\im v_2(y,z)\le \im v_2(y,0)$
            \end{enumerate}
        \end{enumerate}
       \end{cl}

       (a) is because of
       $$\frac{\D\im v_2(y,z)}{\D z}=\im (v_1)'_x(x,y)=2\log|\sin(z+y)|-2\log|\sin(z-y)|$$
       $$\frac{\D^2\im v_2(y,z)}{\D z^2}=2(\cot(z+y)-\cot(z-y))=-\frac{2\sin2y}{\sin(z-y)\sin(z+y)}$$

       We have $\frac{\D\im v_2(y,z)}{\D z}=0\Leftrightarrow y\in\frac{\pi\mathbb{Z}}{2}$ or $z\in\frac{\pi\mathbb{Z}}{2}$. $\frac{\D\im v_2(y,z)}{\D z}$ is positive on $D_{yz2}$ and negative on $D_{yz2}$. Combine this with (a), we get (b).

      \item We divide $D$ into three smaller regions
        $$D=D'\cup D''\cup D'''$$
        where
        \begin{eqnarray}
          D'&=&\{(x,y,z)\in \mathbb{R}^3|y\in[|x|,\frac{\pi}{2}],z\in(y,\pi-y)\}\\
          D''&=&\{(x,y,z)\in \mathbb{R}^3|y\in[|x|,\frac{\pi}{2}],z\in(0,y]\}\\
          D'''&=&\{(x,y,z)\in D|y\in[\frac{\pi}{2},\pi)\}
        \end{eqnarray}

        we have $D_0\in D'$

        \begin{figure}[H]
          \centering
          \includegraphics[width=0.5\textwidth]{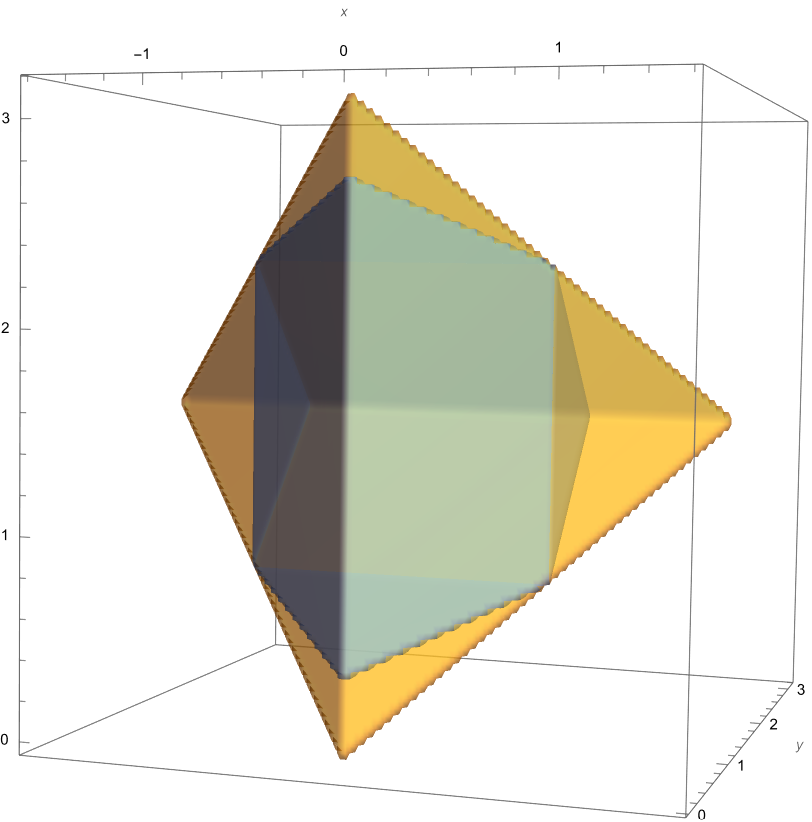}
          \caption{The orange region is $D'$ and the blue region is $D_0$}\label{D0D'}
        \end{figure}

        Combine (1) and (2), we have
        \begin{enumerate}
          \item If $(x,y,z)\in D'$, $\im v(x,y,z)$ is  strictly concave down in $(x,z)$ and $\im v(x,y,z)\le\im v(0,y,\frac{\pi}{2})$.

            By the convexity of $(x,z)$, it suffices to prove that for $(x,y,z)\in \partial D_0$, $\im v(x,y,z)<3.5$.

            There are $8$ faces of $D_0$. Because of the symmetry of $v$, it suffices to prove for $(x,y,z)\in F_i(i=1,2,3,4)$, $\im v(x,y,z)<3.5$, where
            \begin{eqnarray*}
              F_1&=&\{(x,y,z)\in\mathbb{R}^3|x=y,0<y<\frac{\pi}{4},2z+y<\frac{7\pi}{4},2z-y>\frac{\pi}{4}\}\\
              F_2&=&\{(x,y,z)\in\mathbb{R}^3|x=\frac{\pi}{4},y>\frac{\pi}{4},z+y<\pi,z-y>0\}\\
              F_3&=&\{(x,y,z)\in\mathbb{R}^3||x|<y,y=2z-\frac{\pi}{4},\frac{\pi}{8}<z<\frac{\pi}{4}\}\\
              F_4&=&\{(x,y,z)\in\mathbb{R}^3||x|<\frac{\pi}{4},\frac{\pi}{4}<y=z<\frac{\pi}{2}\}
            \end{eqnarray*}

            \begin{enumerate}[1]
              \item For $(x,y,z)\in F_1$, we have $\im v(x,y,z)=\im v(y,y,z)\le\im v(y,y,\frac{\pi}{2})$
                $$\frac{\D \im v(y,y,\frac{\pi}{2})}{\D y}=\im\frac{\D v(y,y,\frac{\pi}{2})}{\D y}=2(2\log2|\cos y|-\log2|\sin y|-2\log2|\sin2y|)=-6\log2|\sin y|$$

                $\frac{\D \im v(y,y,\frac{\pi}{2})}{\D y}$ is a decreasing function with the unique root $y_0=\frac{\pi}{6}$, so we have
                $$\im v(y,y,\frac{\pi}{2})\le\im v(\frac{\pi}{6},\frac{\pi}{6},\frac{\pi}{2})(=3v_3)=3.0448\cdots<3.5$$

               \item For $(x,y,z)\in F_2$, we have $\im v(x,y,z)=\im v(\frac{\pi}{4},y,z)\le\im v(\frac{\pi}{4},y,\frac{\pi}{2})$

                \begin{eqnarray*}
                  &&\frac{\D \im v(\frac{\pi}{4},y,\frac{\pi}{2})}{\D y}=\im\frac{\D v(\frac{\pi}{4},y,\frac{\pi}{2})}{\D y}\\
                  &=&-2(\log2|\sin y|-2\log2|\cos y|+\log2|\sin(y+\frac{\pi}{4})|+\log2|\cos(y+\frac{\pi}{4})|)\\
                  &=&-2(\log2|\sin y|-2\log2|\cos y|+\log2|\cos2y|)\\
                  &=&-2\log\frac{\sin y(2\sin^2y-1)}{1-\sin^2y}
                \end{eqnarray*}

                $\frac{\D \im v(\frac{\pi}{4},y,\frac{\pi}{2})}{\D y}$ is a decreasing function with the unique root $y_0=\arcsin a_0=0.978\cdots$, where $a_0$ is the root of $2a^3+a^2-a-1=0$ which is in $(\frac{\sqrt{2}}{2},1)$. So we have
                $$\im v(\frac{\pi}{4},y,\frac{\pi}{2})<\im v(\frac{\pi}{4},y_0,\frac{\pi}{2})=3.2527\cdots<3.5$$

              \item For $(x,y,z)\in F_3$, we have $\im v(x,y,z)=\im v(x,2z-\frac{\pi}{4},z)\le\im v(0,2z-\frac{\pi}{4},z)$
                \begin{eqnarray*}
                  &&\frac{\D \im v(0,2z-\frac{\pi}{4},z)}{\D z}=\im\frac{\D v(0,2z-\frac{\pi}{4},z)}{\D z}\\
                  &=&2(\log2|\sin(\frac{\pi}{4}-z)|-6\log2|\sin(2z-\frac{\pi}{4})|+3\log2|\sin(3z-\frac{\pi}{4})|)
                \end{eqnarray*}
                \begin{eqnarray*}
                  &&\frac{\D^2 \im v(0,2z-\frac{\pi}{4},z)}{\D z^2}\\
                  &=&-2\tan(z+\frac{\pi}{4})+24\tan(2z+\frac{\pi}{4})-18\tan(3z+\frac{\pi}{4})\\
                  &=&-2\tan(z+\frac{\pi}{4})-24\tan(\frac{3\pi}{4}-2z)+18\tan(\frac{3\pi}{4}-3z)\\
                  &<&-2\tan(z+\frac{\pi}{4})-6\tan(\frac{3\pi}{4}-2z)<0
                \end{eqnarray*}

                So $\frac{\D \im v(0,2z-\frac{\pi}{4},z)}{\D z}$ is a decreasing function with the unique root $z_0=0.674695\cdots$. So we have
                $$\im v(0,2z-\frac{\pi}{4},z)\le \im v(0,2z_0-\frac{\pi}{4},z_0)=3.1439\cdots<3.5$$

              \item For $(x,y,z)\in F_4$, we have $\im v(x,y,z)=\im v(x,y,y)\le\im v(0,y,y)$

              \begin{eqnarray*}
                &&\frac{\D \im v(0,y,y)}{\D y}=\im\frac{\D v(0,y,y)}{\D y}\\
                &=&-2(3\log2|\sin y|-2\log2|\sin2y|)\\
                &=&-2\log\frac{\sin y}{2(1-\sin^2y)}
              \end{eqnarray*}
               So $\frac{\D \im v(0,y,y)}{\D y}$ is a decreasing function with the unique root $y_0=\arcsin\frac{\sqrt{17}-1}{4}$. So we have
               $$\im v(0,y,y)\le\im v(0,y_0,y_0)=2.7868\cdots<3.5$$

            \end{enumerate}

          \item If $(x,y,z)\in D''$, $\im v(x,y,z)\le \im v(x,y,y)$, where $(x,y,y)\in\partial D'$, so it suffices to study the case (a);
          \item If $(x,y,z)\in D'''$, $\im v(x,y,z)\le \im v(\frac{\pi}{2},y,0)=\im v(-\frac{\pi}{2},y,0)$.

            Define $w(y)=v(\frac{\pi}{2},y,0)=2\li(-e^{2\sqrt{-1}y})+\li(e^{-2\sqrt{-1}y})$ for $y\in [\frac{\pi}{2},\pi)$. We have
            $$\frac{\D\im w(y)}{\D y}=\im w'(y)=2\log|2\sin y|-4\log|2\cos y|=2\log\frac{\sin y}{2(1-\sin^2y)}$$
            So $\frac{\D\im w(y)}{\D y}$ is a decreasing function with the unique root $y_0=\pi-\arcsin\frac{\sqrt{17}-1}{4}$. So we have $\im w(y)\le\im w(y_0)=2.7868\cdots<2.8$.
        \end{enumerate}

    \end{enumerate}

\end{proof}

\begin{proof}[Proof of Lemma \ref{imV>3.5}]
    It suffices to prove that for $(x,y,z)\in D_0\backslash D_H$
    $$\im V(x,y+\frac{\log5}{4}\sqrt{-1},z,s,m,n,l)<3.5$$

    Define $$v(x,y,z)=\li(e^{2\sqrt{-1}(y+x)})+\li(e^{2\sqrt{-1}(y-x)})+\li(e^{2\sqrt{-1}(-y+z)})+\li(e^{2\sqrt{-1}(-y-z)})\\
      -\li(e^{-2\sqrt{-1}y})-4y^2$$
    then for $(x,y,z)\in D_0$, $\im v(x,y+\frac{\log5}{4}\sqrt{-1},z)=\im V(x,y+\frac{\log5}{4}\sqrt{-1},z,s,m,n,l)$.
    \begin{enumerate}[(1)]
      \item First let's study the imaginary part of the function
      $$v_1(x,y)=\li(e^{2\sqrt{-1}(y+x)})+\li(e^{2\sqrt{-1}(y-x)})$$
       on $D_{xy}+(0,\frac{\log5}{4}\sqrt{-1})$ where $D_{xy}=\{(x,y)\in\mathbb{R}^2|y\in(|x|,\frac{\pi}{2}),|x|<\frac{\pi}{4}\}$. We have
       \begin{eqnarray*}
         &&\frac{\D\im v_1(x,y+\frac{\log5}{4}\sqrt{-1})}{\D x}=\im (v_1)'_x(x,y+\frac{\log5}{4}\sqrt{-1})\\
         &=&2\log|\frac{1-e^{2\sqrt{-1}(y+\frac{\log5}{4}\sqrt{-1}-x)}}{1-e^{2\sqrt{-1}(y+\frac{\log5}{4}\sqrt{-1}+x)}}|\\
         &=&\log\frac{3-\sqrt{5}\cos2(y-x)}{3-\sqrt{5}\cos2(y+x)}\\
         &=&\log(1-\frac{2\sin2x\sin2y}{\frac{3}{\sqrt{5}}-\cos2(y+x)})
       \end{eqnarray*}
       so $\frac{\D\im v_1(x,y+\frac{\log5}{4}\sqrt{-1})}{\D x}\begin{cases}
                                                                 >0, & \mbox{if } x\in(-y,0) \\
                                                                 <0, & \mbox{if } x\in(0,y)
                                                               \end{cases}$.

       We have for $(x,y,z)\in D_0+(0,\frac{\log5}{4}\sqrt{-1},0)$, $\im v(x,y,z)$ is an increasing function of $x$ when $x<0$, and an decreasing function of $x$ when $x>0$.

      \item  Let's study the imaginary part of the function
        $$v_2(x,y)=\li(e^{2\sqrt{-1}(-y+z)})+\li(e^{2\sqrt{-1}(-y-z)})$$
        on $D_{yz}+(\frac{\log5}{4}\sqrt{-1},0)$ where
        $$D_{yz}=\{(y,z)\in\mathbb{R}^2|y>0,z>y,z+y<\pi,y+2z<\frac{7}{4}\pi,2z-y>\frac{\pi}{4}\}$$

        We have
       \begin{eqnarray*}
         &&\frac{\D\im v_2(y+\frac{\log5}{4}\sqrt{-1},z)}{\D z}=\im (v_2)'_z(y+\frac{\log5}{4}\sqrt{-1},z)\\
         &=&2\log|\frac{1-e^{2\sqrt{-1}(-y-\frac{\log5}{4}\sqrt{-1}-z)}}{1-e^{2\sqrt{-1}(-y-\frac{\log5}{4}\sqrt{-1}+z)}}|\\
         &=&\log\frac{3-\sqrt{5}\cos2(-y-z)}{3-\sqrt{5}\cos2(-y+z)}\\
         &=&\log(1+\frac{2\sin2z\sin2y}{\frac{3}{\sqrt{5}}-\cos2(z-y)})
       \end{eqnarray*}
       so $\frac{\D\im v_2(y+\frac{\log5}{4}\sqrt{-1},z)}{\D z}\begin{cases}
                                                                 >0, & \mbox{if } z<\frac{\pi}{2} \\
                                                                 <0, & \mbox{if } z>\frac{\pi}{2}
                                                               \end{cases}$.

       We have for $(x,y,z)\in D_0+(0,\frac{\log5}{4}\sqrt{-1},0)$, $\im v(x,y,z)$ is an increasing function of $z$ when $z<\frac{\pi}{2}$, and an decreasing function of $z$ when $z>\frac{\pi}{2}$.

      \item  Combine (1) and (2), it suffices to prove that for $(x,y,z)\in \partial D_H+(0,\frac{\log5}{4}\sqrt{-1},0)$, $\im v(x,y,z)<3.5$.

            There are $8$ faces of $D_H+(0,\frac{\log5}{4}\sqrt{-1},0)$. Because of the symmetry of $v$, it suffices to prove for $(x,y,z)\in F_i+(0,\frac{\log5}{4}\sqrt{-1},0)$, $\im v(x,y,z)<3.5$, where
            \begin{eqnarray*}
              F_1&=&\{(x,y,z)\in\mathbb{R}^3|x=y,0<y<\frac{\pi}{4},z\in(\frac{\pi}{2}-y,\frac{\pi}{2}+y)\}\\
              F_2&=&\{(x,y,z)\in\mathbb{R}^3|x=\frac{\pi}{2}-y,\frac{\pi}{4}<y<\frac{\pi}{2},z\in(y,\pi-y)\}\\
              F_3&=&\{(x,y,z)\in\mathbb{R}^3|x\in(-y,y),0<y<\frac{\pi}{4},z=\frac{\pi}{2}-y\}\\
              F_4&=&\{(x,y,z)\in\mathbb{R}^3|x\in(y-\frac{\pi}{2},\frac{\pi}{2}-y),\frac{\pi}{4}<y<\frac{\pi}{2},z=y\}
            \end{eqnarray*}

            \begin{enumerate}[1]
              \item For $(x,y,z)\in F_1$, we have
                $$\im v(x,y+\frac{\log5}{4}\sqrt{-1},z)=\im v(y,y+\frac{\log5}{4}\sqrt{-1},z)\le\im v(y,y+\frac{\log5}{4}\sqrt{-1},\frac{\pi}{2})$$

                \begin{eqnarray*}
                  &&\frac{\D \im v(y,y+\frac{\log5}{4}\sqrt{-1},\frac{\pi}{2})}{\D y}=\im\frac{\D v(y,y+\frac{\log5}{4}\sqrt{-1},\frac{\pi}{2})}{\D y}\\
                  &=&-\log(6-2\sqrt{5}\cos2y)+2\log(6+2\sqrt{5}\cos2y)-2\log(6-2\sqrt{5}\cos4y)
                \end{eqnarray*}

                $\frac{\D \im v(y,y+\frac{\log5}{4}\sqrt{-1},\frac{\pi}{2})}{\D y}$ is a decreasing function with the unique root $y_0=0.372498\cdots$, so we have
                $$\im v(y,y+\frac{\log5}{4}\sqrt{-1},\frac{\pi}{2})\le\im v(y_0,y_0+\frac{\log5}{4}\sqrt{-1},\frac{\pi}{2})=3.2543\cdots<3.5$$

              \item For $(x,y,z)\in F_2$, we have
                $$\im v(x,y+\frac{\log5}{4}\sqrt{-1},z)=\im v(\frac{\pi}{2}-y,y+\frac{\log5}{4}\sqrt{-1},z)\le\im v(\frac{\pi}{2}-y,y+\frac{\log5}{4}\sqrt{-1},\frac{\pi}{2})$$

                \begin{eqnarray*}
                  &&\frac{\D \im v(\frac{\pi}{2}-y,y+\frac{\log5}{4}\sqrt{-1},\frac{\pi}{2})}{\D y}=\im\frac{\D v(\frac{\pi}{2}-y,y+\frac{\log5}{4}\sqrt{-1},\frac{\pi}{2})}{\D y}\\
                  &=&-\log(6-2\sqrt{5}\cos2y)+2\log(6+2\sqrt{5}\cos2y)-2\log(6+2\sqrt{5}\cos4y)
                \end{eqnarray*}

                $\frac{\D \im v(\frac{\pi}{2}-y,y+\frac{\log5}{4}\sqrt{-1},\frac{\pi}{2})}{\D y}$ is a decreasing function with the unique root $y_0=0.891\cdots$, so we have
                $$\im v(\frac{\pi}{2}-y,y+\frac{\log5}{4}\sqrt{-1},\frac{\pi}{2})\le\im v(\frac{\pi}{2}-y_0,y_0+\frac{\log5}{4}\sqrt{-1},\frac{\pi}{2})=2.635\cdots<3.5$$

              \item For $(x,y,z)\in F_3$, we have
                $$\im v(x,y+\frac{\log5}{4}\sqrt{-1},z)=\im v(x,y+\frac{\log5}{4}\sqrt{-1},\frac{\pi}{2}-y)\le\im v(0,y+\frac{\log5}{4}\sqrt{-1},\frac{\pi}{2}-y)$$

                \begin{eqnarray*}
                  &&\frac{\D \im v(0,y+\frac{\log5}{4}\sqrt{-1},\frac{\pi}{2}-y)}{\D y}=\im\frac{\D v(0,y+\frac{\log5}{4}\sqrt{-1},\frac{\pi}{2}-y)}{\D y}\\
                  &=&-3\log(6-2\sqrt{5}\cos2y)+2\log(6+2\sqrt{5}\cos4y)
                \end{eqnarray*}

                $\frac{\D \im v(0,y+\frac{\log5}{4}\sqrt{-1},\frac{\pi}{2}-y)}{\D y}$ is a decreasing function with the unique root $y_0=0.4278594\cdots$, so we have
                $$\im v(0,y+\frac{\log5}{4}\sqrt{-1},\frac{\pi}{2}-y)\le\im v(0,y_0+\frac{\log5}{4}\sqrt{-1},\frac{\pi}{2}-y_0)=3.4595\cdots<3.5$$

              \item For $(x,y,z)\in F_4$, we have
                $$\im v(x,y+\frac{\log5}{4}\sqrt{-1},z)=\im v(x,y+\frac{\log5}{4}\sqrt{-1},y+0)\footnote{$"+0"$ means that take the limit from the right side: $v(x,y+\frac{\log5}{4}\sqrt{-1},y+0)=\lim_{z\to y+0}\limits v(x,y+\frac{\log5}{4}\sqrt{-1},z)$. The reason why we have to take the limit here is that the $\im v$ is a multi-valued function, and it is well-defined by taking the limit.}\le\im v(0,y+\frac{\log5}{4}\sqrt{-1},y+0)$$

                \begin{eqnarray*}
                  &&\frac{\D \im v(0,y+\frac{\log5}{4}\sqrt{-1},y+0)}{\D y}=\im\frac{\D v(0,y+\frac{\log5}{4}\sqrt{-1},y+0)}{\D y}\\
                  &=&-3\log(6-2\sqrt{5}\cos2y)+2\log(6-2\sqrt{5}\cos4y)
                \end{eqnarray*}

                $\frac{\D \im v(0,y+\frac{\log5}{4}\sqrt{-1},y+0)}{\D y}$ is a decreasing function
                $$\frac{\D \im v(0,y+\frac{\log5}{4}\sqrt{-1},y)}{\D y}<\frac{\D \im v(0,\frac{\pi}{4}+\frac{\log5}{4}\sqrt{-1},\frac{\pi}{4}+0)}{\D y}=\log\frac{(6+2\sqrt{5})^2}{6^3}<0$$
                so $\im v(0,y+\frac{\log5}{4}\sqrt{-1},y+0)$ is a decreasing function, we have
                $$\im v(0,y+\frac{\log5}{4}\sqrt{-1},y+0)<\im v(0,\frac{\pi}{4}+\frac{\log5}{4}\sqrt{-1},\frac{\pi}{4}+0)=2.5778\cdots<3.5$$

            \end{enumerate}

    \end{enumerate}

\end{proof}

\begin{lem}\label{f(y0)}
    Let $f(y)=2\li(e^{2\sqrt{-1}y})+2\li(-e^{-2\sqrt{-1}y})-\li(e^{-2\sqrt{-1}y})-4y^2$. Then we have
    \begin{enumerate}[(1)]
      \item $f(\frac{\arctan2}{2}+\frac{\log5}{4}\sqrt{-1})=-\frac{\pi^2}{4}+4\dd(\sqrt{-1})\sqrt{-1}$
      \item $f(-\frac{\arctan2}{2}+\frac{\log5}{4}\sqrt{-1})=-\frac{\pi^2}{4}-4\dd(\sqrt{-1})\sqrt{-1}$
    \end{enumerate}
\end{lem}

\begin{proof}[Proof of Lemma \ref{f(y0)}]
    Because of
    $$\overline{f(y)}=f(-\overline{y})$$
    it suffices to prove (1). We have the following steps.

    \begin{enumerate}
      \item     We have $e^{2\sqrt{-1}y}|_{y=\frac{\arctan2}{2}+\frac{\log5}{4}\sqrt{-1}}=\frac{1+2\sqrt{-1}}{5}=\frac{1}{1-2\sqrt{-1}}$,
        \begin{eqnarray}
          &&f(\frac{\arctan2}{2}+\frac{\log5}{4}\sqrt{-1})\nonumber\\
          &=&2\li(\frac{1}{1-2\sqrt{-1}})+2\li(2\sqrt{-1}-1)-\li(1-2\sqrt{-1})+\log^2(1-2\sqrt{-1})\nonumber\\
          &\stackrel{(\ref{e1})}{=}&2\li(2\sqrt{-1}-1)-3\li(1-2\sqrt{-1})-\frac{\pi^2}{3}-\log^2(2\sqrt{-1}-1)+\log^2(1-2\sqrt{-1})\label{eap1}
        \end{eqnarray}

      \item We have the following identities of dilogarithm function.
        \begin{equation}\label{edi2}
          \li(z)+\li(1-z)=\frac{\pi^2}{6}-\log z\log(1-z)
        \end{equation}

        By (\ref{e1}) and (\ref{edi2}), we have
        \begin{equation}\label{edi3}
          \li(z)-\li(\frac{1}{1-z})=\frac{\pi^2}{3}+\frac{\log^2(z-1)}{2}-\log z\log(1-z)
        \end{equation}

      \item We have $\li(-1)=-\frac{\pi^2}{12}$, $\im\li(\sqrt{-1})=\dd(\sqrt{-1})$ and
        $$2\re\li(\sqrt{-1})=\li(\sqrt{-1})+\li(-\sqrt{-1})\stackrel{\ref{e1}}{=}-\frac{\pi^2}{6}-\frac{\log^2(-\sqrt{-1})}{2}=-\frac{\pi^2}{24}$$

        So we have
        \begin{eqnarray}
          \li(\sqrt{-1})&=&-\frac{\pi^2}{48}+\dd(\sqrt{-1})\sqrt{-1}\label{edi4}\\
          \li(-\sqrt{-1})&=&-\frac{\pi^2}{48}-\dd(\sqrt{-1})\sqrt{-1}\label{edi5}
        \end{eqnarray}

        \begin{eqnarray}
          &&\li(1+\sqrt{-1})+\li(-\sqrt{-1})\stackrel{(\ref{edi2})}{=}\frac{\pi^2}{6}-\log(-\sqrt{-1})\log(1+\sqrt{-1})\nonumber\\
          &\stackrel{(\ref{edi5})}{\Rightarrow}&\li(1+\sqrt{-1})=\frac{\pi^2}{48}+\dd(\sqrt{-1})\sqrt{-1}+\frac{\pi^2}{6}-\log(-\sqrt{-1})\log(1+\sqrt{-1})\label{edi6}\\
          &&\li(\frac{1+\sqrt{-1}}{2})=\li(\frac{1}{1-\sqrt{-1}})\nonumber\\
          &\stackrel{(\ref{edi3})}{=}&\li(\sqrt{-1})-\frac{\pi^2}{3}-\frac{\log^2(-1+\sqrt{-1})}{2}+\log\sqrt{-1}\log(1-\sqrt{-1})\nonumber\\
          &\stackrel{(\ref{edi4})}{=}&-\frac{\pi^2}{48}+\dd(\sqrt{-1})\sqrt{-1}-\frac{\pi^2}{3}-\frac{\log^2(-1+\sqrt{-1})}{2}+\log\sqrt{-1}\log(1-\sqrt{-1})\label{edi7}\\
          &&\li(-1+2\sqrt{-1})+\li(2-2\sqrt{-1})\stackrel{(\ref{edi2})}{=}\frac{\pi^2}{6}-\log(-1+2\sqrt{-1})\log(2-2\sqrt{-1})\label{edi8}
        \end{eqnarray}

      \item By Proposition \ref{5terms} we have
        \begin{eqnarray}
          &&\li(1-2\sqrt{-1})+\li(2-2\sqrt{-1})-\li(-1)-\li(-1-\sqrt{-1})+\li(1+\sqrt{-1})\nonumber\\
          &=&-\log(2\sqrt{-1})\log(-1+2\sqrt{-1})\label{e5term1}\\
          &&\li(1-2\sqrt{-1})+\li(-1-\sqrt{-1})-\li(-\sqrt{-1})-\li(\frac{-1+\sqrt{-1}}{2})+\li(\frac{1+\sqrt{-1}}{2})\nonumber\\
          &=&-\log(2\sqrt{-1})\log(2+\sqrt{-1})\label{e5term2}
        \end{eqnarray}

      \item Take (\ref{e5term1})$\times2+$(\ref{e5term2}), we have
        \begin{eqnarray*}
          &&3\li(1-2\sqrt{-1})+2\li(2-2\sqrt{-1})+\frac{\pi^2}{6}+(2\li(1+\sqrt{-1})-\li(-\sqrt{-1})+\li(\frac{1+\sqrt{-1}}{2}))\\
          &&-(\li(\frac{-1+\sqrt{-1}}{2})+\li(-1-\sqrt{-1}))=-\log(2\sqrt{-1})(2\log(-1+2\sqrt{-1})+\log(2+\sqrt{-1}))\\
          &\Rightarrow&3\li(1-2\sqrt{-1})+2(-\li(-1+2\sqrt{-1})+\frac{\pi^2}{6}-\log(-1+2\sqrt{-1})\log(2-2\sqrt{-1}))\footnote{By (\ref{edi8})}+\frac{\pi^2}{6}\\
          &&+(4\dd(\sqrt{-1})\sqrt{-1}+\frac{\pi^2}{24}-2\log(-\sqrt{-1})\log(1+\sqrt{-1})-\frac{\log^2(-1+\sqrt{-1})}{2}+\log\sqrt{-1}\log(1-\sqrt{-1}))\footnote{By (\ref{edi5}), (\ref{edi6}), (\ref{edi7})}\\
          &&-(-\frac{\pi^2}{6}-\frac{\log^2(1+\sqrt{-1})}{2})\footnote{By (\ref{e1})}=-\log(2\sqrt{-1})(2\log(-1+2\sqrt{-1})+\log(2+\sqrt{-1}))\\
          &\Rightarrow&2\li(2\sqrt{-1}-1)-3\li(1-2\sqrt{-1})=\frac{\pi^2}{2}-2\log(-1+2\sqrt{-1})\log(2-2\sqrt{-1})+4\dd(\sqrt{-1})\sqrt{-1}\\
          &&+\frac{\pi^2}{24}-2\log(-\sqrt{-1})\log(1+\sqrt{-1})-\frac{\log^2(-1+\sqrt{-1})}{2}+\log\sqrt{-1}\log(1-\sqrt{-1})+\frac{\pi^2}{6}\\
          &&+\frac{\log^2(1+\sqrt{-1})}{2}+\log(2\sqrt{-1})(2\log(-1+2\sqrt{-1})+\log(2+\sqrt{-1}))\\
          &\Rightarrow&(\ref{eap1})=\frac{\pi^2}{2}-2\log(-1+2\sqrt{-1})\log(2-2\sqrt{-1})+4\dd(\sqrt{-1})\sqrt{-1}+\frac{\pi^2}{24}-2\log(-\sqrt{-1})\log(1+\sqrt{-1})\\
          &&-\frac{\log^2(-1+\sqrt{-1})}{2}+\log\sqrt{-1}\log(1-\sqrt{-1})+\frac{\pi^2}{6}+\frac{\log^2(1+\sqrt{-1})}{2}\\
          &&+\log(2\sqrt{-1})(2\log(-1+2\sqrt{-1})+\log(2+\sqrt{-1}))-\frac{\pi^2}{3}-\log^2(2\sqrt{-1}-1)+\log^2(1-2\sqrt{-1})\\
          &&=-\frac{\pi^2}{4}+4\dd(\sqrt{-1})\sqrt{-1}
        \end{eqnarray*}

    \end{enumerate}

\end{proof}

 \bibliographystyle{plain}
  \bibliography{bibtex}

~

~

\noindent Huabin Ge, hbge@ruc.edu.cn\\[2pt]
\emph{School of Mathematics, Renmin University of China, Beijing 100872, P. R. China}\\

\noindent Yunpeng Meng, 2230501005@cnu.edu.cn\\[2pt]
\emph{School of Mathematical Sciences, Capital Normal University, Beijing 100048, P. R. China}\\

\noindent Chuwen Wang, chuwenwang@ruc.edu.cn\\[2pt]
\emph{School of Mathematics, Renmin University of China, Beijing 100872, P. R. China}\\

\noindent Yuxuan Yang, 2001110022@pku.edu.cn\\[2pt]
\emph{Beijing International Center for Mathematical Research, Peking University, Beijing 100871, P. R. China}

\end{document}